%% file: DeWit99c.tex
\documentclass[fleqn,10pt]{article}

\usepackage{latexsym}
\usepackage{amssymb}

\def\ket#1{\left\vert #1 \right\rangle}
\def\bra#1{\left\langle #1 \right\vert}
\def\ketbra#1#2{\left\vert #1\right\rangle\!\left\langle #2\right\vert}
\def\ketket#1#2{\left\vert #1 \! \otimes \! #2 \right\rangle}
\def\braket#1#2{\left\langle #1 \vert #2 \right\rangle}
\def\brabra#1#2{\left\langle #1 \! \otimes \! #2 \right\vert}

\newtheorem{lemma}{Lemma}

\title{
  Automatic Construction of Explicit R Matrices for the One-Parameter
  Families of Irreducible Typical Highest
  Weight $(\dot{0}_m|\dot{\alpha}_n)$ Representations of
  {$U_q[gl(m|n)]$}
}

\author{
  David~~De Wit
  \\
  RIMS, Kyoto University, JAPAN
}

\begin{document}

\maketitle

\begin{abstract}
  \noindent
   We detail the automatic construction of R matrices corresponding to
   (the tensor products of) the $(\dot{0}_m|\dot{\alpha}_n)$ families
   of highest-weight representations of the quantum superalgebras
   $U_q[gl(m|n)]$.  These representations are irreducible, contain
   a free complex parameter $\alpha$, and are $2^{mn}$ dimensional.
   Our R matrices are actually (sparse) rank $4$ tensors, containing a
   total of $2^{4mn}$ components, each of which is in general an
   algebraic expression in the two complex variables $q$ and $\alpha$.

   Although the constructions are straightforward, we describe them in
   full here, to fill a perceived gap in the literature.  As the
   algorithms are generally impracticable for manual calculation, we
   have implemented the entire process in \textsc{Mathematica};
   illustrating our results with $U_q[gl(3|1)]$.

\end{abstract}



\section{Introduction}

Broadly, R matrices are solutions to the various versions of the
Yang--Baxter equation, and as such, are of great interest in
mathematical physics and knot theory (see, e.g.  \cite{Kauffman:93}),
both in their algebraic (i.e. ``universal'') forms, and in their
(matrix) representations (i.e. ``quantum'' forms), useful for explicit
computations.  Here, we will be specifically concerned with quantum R
matrices associated with the quantum superalgebras $U_q[gl(m|n)]$.

Although much is known about the origin and properties of quantum
superalgebra R matrices (e.g. \cite{KhoroshkinTolstoy:91} provides
universal R matrices), explicit examples of their quantum R matrices
are rare in the literature, due largely to the computational effort
involved in obtaining them.  This paper describes the automation of an
algorithm to generate a suite of explicit quantum R matrices for
$U_q[gl(m|n)]$.  As readers of this organ may not be familiar with
these algebraic structures, we provide a full description of their
details.

Specifically, we construct \emph{trigonometric} R matrices
$\check{R}^{m,n}(u)$ corresponding to the $\alpha$-parametric highest
weight minimal representations labeled $(\dot{0}_m|\dot{\alpha}_n)$ of
the $U_q[gl(m|n)]$.  These irreducible representations are $2^{mn}$
dimensional, and contain free complex parameters $q$ and $\alpha$; the
real variable $u$ is a `spectral' parameter.  Quantum R matrices
$\check{R}^{m,n}$ are immediately obtainable as the spectral limits
$u\to\infty$ of the $\check{R}^{m,n}(u)$.

Our R matrices are in fact \emph{graded}, as they are based on graded
vector spaces, hence they actually satisfy \emph{graded} Yang--Baxter
equations.  However, it is a simple matter to remove this grading and
transform them into objects that satisfy the usual Yang--Baxter
equations.

The constructions have been implemented in \textsc{Mathematica}, and
results obtained for $n=1$ and $m=1,2,3,4$; we illustrate the
algorithms using $U_q[gl(3|1)]$. Full listings of all our R matrices
have been announced in \cite{DeWit:99d}.

As they are solutions to Yang--Baxter equations, our R matrices are
of immediate practical interest.  Firstly, they are of physical
interest in that they are applicable to the construction of exactly
solvable models of interacting fermions.  Corresponding to
$\check{R}^{m,1}(u)$, we may construct an integrable $2^m$ state
fermionic model on a lattice.  Models associated with $U_q[gl(2|1)]$
and $U_q[gl(3|1)]$ have been discussed in
\cite{GouldHibberdLinksZhang:96} and \cite{GeGouldZhangZhou:98a},
respectively. The $U_q[gl(4|1)]$ case has an elegant interpretation in
terms of a $2$-leg ladder model for interacting electrons: a discussion
of this is provided in \cite{DeWit:99d}.

Furthermore, corresponding to each $\check{R}^{m,n}$, we may obtain a
polynomial `Links--Gould' link invariant $LG^{m,n}$
\cite{LinksGould:92b}, cf.  the celebrated Jones polynomial.  These
$LG^{m,n}$ are two-variable, integer-coefficient Laurent polynomials,
and are generally substantially more powerful than the Jones polynomial
in distinguishing knots.  ($LG^{1,1}$ degenerates to the
well-known Alexander--Conway polynomial in the single variable
$q^{2\alpha}$ (cf.  \cite{BrackenGouldZhangDelius:94b}).) A fuller
documentation of the suite $LG^{m,n}$ has been provided by myself in
collaboration with Louis Kauffman and Jon Links in
\cite{DeWit:98,DeWit:99a,DeWit:99e,DeWitKauffmanLinks:99a}.  Although
the $LG^{m,n}$ are far from being complete invariants, as they can
distinguish neither mutants nor inversion
\cite{DeWit:99a,DeWitKauffmanLinks:99a}, it turns out that even
$LG^{2,1}$ is in fact more powerful than the well-known two-variable
HOMFLY and Kauffman invariants, being able to distinguish (including
chirality) all prime knots of up to $10$ crossings \cite{DeWit:99a}.
Their evaluation also involves automatic symbolic computation, but the
computational aspects are comparatively pedestrian.

Lastly, we mention explicitly that this paper contains no new theorems,
although it does contain two new technical lemmas, proven in Appendix
\ref{app:Lemmas}.  It is primarily intended to provide a proper
foundation for the results presented in \cite{DeWit:99d,DeWit:99e},
although it also serves as a tutorial on an application of symbolic
computation. Whilst it specifically pertains to representations of
$U_q[gl(m|n)]$, many of the algorithms have a much broader
application.

The following subsections provide a synopsis of the paper.

\pagebreak


\subsection{Algebraic overview}

Fixing $m$ and $n$, we are initially interested in a $2^{mn}$
dimensional vector space $V$ that is a module for the $U_q[gl(m|n)]$
minimal typical highest weight representation
$\Lambda=(\dot{0}_m|\dot{\alpha}_n)$.  The algebra contains a free
complex variable $q$, whilst the representation $\pi_{\Lambda}$ acting
on $V$ contains a free complex variable $\alpha$.  Our $V$ is actually
($\mathbb{Z}_2$) \emph{graded}; this ensures compatibility with the
($\mathbb{Z}_2$) grading of $U_q[gl(m|n)]$.

Using the properties of $U_q[gl(m|n)]$, we apply a version of the
\emph{Kac induced module construction} (KIMC) \cite{Kac:77b,Kac:78} to
establish a (\emph{weight}) basis $\{\ket{i}\}_{i=1}^{2^{mn}}$ for $V$.
This involves postulating $\ket{1}$ as a highest weight vector, and
recursively acting on $\ket{1}$ with all possible distinct products of
simple lowering generators ${E^{a+1}}_a$ to define the other basis
vectors, normalising as we go.  This construction requires a
`Poincar\'{e}--Birkhoff--Witt (PBW) lemma' for $U_q[gl(m|n)]$
\cite{DeWit:2000,Zhang:93}, i.e. a set of commutations sufficient to
transform any product of algebra generators into a normal form (see
\S\ref{sec:Commutations}), together with a statement that the
algebra is spanned by the set of all such normal forms.

Where $V$ has a graded weight basis
$\{\ket{i}\}_{i=1}^{2^{mn}}$, the tensor product module
$V\otimes V$ has a natural $2^{2mn}$ dimensional basis
$\{\ket{i}\otimes\ket{j}\}_{i,j=1}^{2^{mn}}$, which inherits a weight
system and a grading from $V$.  For our particular representation, the
orthogonal decomposition of $V\otimes V$ is known
\cite{DeliusGouldLinksZhang:95b}, and contains no multiplicities, viz:
\begin{eqnarray*}
  V \otimes V
  =
  \bigoplus_{k}
    V_k,
\end{eqnarray*}
where the submodule $V_k$ has highest weight $\lambda_k$, and these
$\lambda_k$ are known, and all distinct.  To build R matrices acting on
$V\otimes V$, we require an alternative, orthonormal weight basis
$\mathfrak{B}=\bigcup_{k} \mathfrak{B}_k$ for $V \otimes V$,
corresponding to this decomposition, viz $\mathfrak{B}_k$ is a basis
for $V_k$. Again using the KIMC, the basis vectors of each
$\mathfrak{B}_k$ are derived as linear combinations of the form
$\theta_{ij}(\ket{i}\otimes\ket{j})$, where the coefficients
$\theta_{ij}$ are algebraic expressions in $q$ and $\alpha$.  This
process initially yields a basis $\overline{\mathfrak{B}}_k$ that is
not necessarily orthonormal, so we also apply a Gram--Schmidt process
to orthonormalise $\overline{\mathfrak{B}}_k$ into $\mathfrak{B}_k$.
The desired R matrix is then a weighted sum of projectors onto these
$V_k$, where the weights are eigenvalues of the appropriate second
order Casimir invariants.

The algebraic structure of $U_q[gl(m|n)]$ is detailed in
\S\ref{sec:UqglmnStructure}, and an introduction to its highest weight
representations is provided in
\S\ref{sec:HighestweightUqglmnrepresentations}.  In
\S\ref{sec:UqglmnPBW}, we provide a normal ordering and a PBW lemma for
$U_q[gl(m|n)]$.  The construction of our particular
$(\dot{0}_m|\dot{\alpha}_n)$ representations is detailed in
\S\ref{sec:Uqglmnalphareps}.  \S\ref{sec:ThesubmodulesVk} describes the
construction of the bases $\mathfrak{B}_k$, and
\S\ref{sec:ProjectorsRMatrices} describes the construction of
projectors and R matrices.

\pagebreak


\subsection{Implementation and results}

Explicit computations within the representation theory of quantum
superalgebras are tedious and error-prone when performed manually.  The
dimensions of representations are generally large, and in our case, we
have the presence of the two variables $q$ and $\alpha$; these
generally manifest themselves in complicated rational algebraic
expressions, whose symmetries must be continually identified and
exploited to avoid the arising of intractable messes of algebra.

The construction of the basis $\{\ket{i}\}_{i=1}^{2^{mn}}$ involves
many applications of the PBW lemma to simplify long strings of algebra
generators.  This is computationally expensive; firstly as the
simplification involves a minimally-efficient sorting process, and
second as it involves a geometric explosion in the number of terms
being sorted.

The construction of the weight space bases $\mathfrak{B}_k$ is
nontrivial, as each basis vector of each $\mathfrak{B}_k$ generally
contains many terms of the form $\theta_{ij}(\ket{i}\otimes \ket{j})$,
where the coefficients $\theta_{ij}$ are generally complicated rational
algebraic expressions in $q$ and $\alpha$.  (That said, we have avoided
the more theoretically difficult situation of computing weight space
decompositions in cases where there are weight multiplicities in the
underlying carrier space $V$.) Although the R matrices have $2^{4mn}$
components, Nature is kind to us in that most of these components are
zero, and those that are not are generally simpler than the
$\theta_{ij}$.

To the best of our knowledge, computer implementation of the algebraic
structures and algorithms described herein has not previously been
achieved.  We have implemented the entire process as a suite of
\textsc{Mathematica} functions; the thousands of lines of code perform
algebraic computations that a human being could not ever realistically
expect to perform correctly.

From \S\ref{sec:Uqglmnrootsystem} onwards, we use $U_q[gl(3|1)]$ to
illustrate our results. These are summarised in Appendix
\ref{app:Data}; where we list the explicit matrix elements for the
generators of the underlying $8$ dimensional representation,
orthonormal bases $\mathfrak{B}_k$ for the $4$ submodules
$V_k\subset V\otimes V$, the components of the associated $4$
projectors $P_k$ onto the $V_k$, and finally, the trigonometric and
quantum R matrices, $\check{R}^{3,1}(u)$ and $\check{R}^{3,1}$,
respectively.

Whilst there are no theoretical limits to $m$ and $n$, a current
practical limit for computation is $mn\leqslant 4$.  This is
convenient, as an immediate application \cite{DeWit:99d} of the
material critically requires $\check{R}^{4,1}(u)$.  Although translation
of the interpreted \textsc{Mathematica} code into a compiled language
would increase the speed of the computations, storage requirements
would still limit $mn$ to perhaps $7$ in the general case.

Further discussion of implementational issues and results is provided
in \S\ref{sec:ImplementationResults}.


\section{The quantum superalgebras $U_q[gl(m|n)]$}
\label{sec:UqglmnStructure}

The algebraic structures labeled $U_q[gl(m|n)]$ are
\emph{quantum superalgebras},%
\footnote{%
  These structures are sometimes called ``quantum super\emph{groups}'',
  but they \emph{are} actually (associative, noncommutative) algebras.
}
described in many places, e.g.
\cite{%
 BrackenGouldZhang:90,%
 DeliusGouldLinksZhang:95a,%
 ScheunertNahmRittenberg:77a,%
 Yamane:91,%
 Yamane:94%
},
and in the book \cite[see~\S6.5]{ChariPressley:94}.

For our purposes, $m$ and $n$ are positive integers, to be regarded as
\emph{fixed}, and $q$ is to be regarded as a nonzero complex variable.
As $U_q[gl(m|n)]$ may be unfamiliar to the readers of this organ, in
\S\ref{sec:UqglmnPhylogeny} we introduce its phylogeny, and in
\S\ref{sec:UqglmnFullDescription}, we provide a full description of its
structure in terms of generators and relations. Beyond that, in
\S\ref{sec:Uqglmnrootsystem} we describe its root system, and in
\S\ref{sec:UqglmnasaHopfsuperalgebra} we show how it may be regarded as
a Hopf (super)algebra.


\subsection{The phylogeny of $U_q[gl(m|n)]$}
\label{sec:UqglmnPhylogeny}

\begin{enumerate}
\item
  Where $n$ is a positive integer, recall that the \emph{Lie algebra}
  $gl(n)$ is equivalent to the usual (complex) vector space of
  $n\times n$ (complex) matrices augmented by a `vector multiplication'
  operation which is the usual matrix multiplication. $gl(n)$ is of
  course a unital algebra, and is of dimension $n^2$ and rank $n-1$.
  The $n^2$ generators ${\{{e^a}_b\}}_{a,b=1}^{n}$ of $gl(n)$ satisfy a
  \emph{commutation relation}:
  \begin{eqnarray*}
    [{e^a}_b, {e^c}_d]
    =
    \delta^c_b {e^a}_d
    -
    \delta^a_d {e^c}_b,
  \end{eqnarray*}
  where $[\cdot,\cdot]$ is the usual \emph{commutator (bracket)},
  defined for $X,Y\in gl(n)$ by:
  \begin{eqnarray*}
    [X, Y]
    \triangleq
    XY-YX.
  \end{eqnarray*}

\item
  Letting both $m$ and $n$ be positive integers, the
  \emph{Lie superalgebra} $gl(m|n)$ may be obtained from $gl(m+n)$ by
  retaining the generators ${\{{e^a}_b\}}_{a,b=1}^{m+n}$, but
  modifying the definition of the commutator bracket and commutation
  relations to include some `parity factors' of $\pm 1$. Specifically,
  we have the commutation relation:
  \begin{equation}
    [{e^{a}}_{b}, {e^{c}}_{d} ]
    =
    \delta^c_b {e^{a}}_{d}
    -
    {(-)}^{[{e^{a}}_{b}][{e^{c}}_{d}]}
    \delta^a_d {e^{c}}_{b},
    \label{eq:glmnCommutationRelations}
  \end{equation}
  where $[\cdot,\cdot]$ is now the \emph{graded commutator (bracket)},
  defined for homogeneous (see below) $X,Y\in gl(m|n)$ by:
  \begin{equation}
    [X, Y]
    \;
    \triangleq
    \;
    XY
    -
    (-)^{[X][Y]}
    YX,
    \label{eq:glmnCommutatorBracket}
  \end{equation}
  and extended by linearity.  In both
  (\ref{eq:glmnCommutationRelations}) and
  (\ref{eq:glmnCommutatorBracket}), $[X]\in\{0,1\}$ refers to the
  \emph{grading} of the homogeneous element $X$.  For this reason, Lie
  superalgebras are sometimes called ``graded Lie algebras''.
  From the $gl(m+n)$ case, we see that $gl(m|n)$ is of dimension
  $(m+n)^2$ and rank $m+n-1$.

\item
  $U[gl(m|n)]$ is then the usual \emph{universal enveloping algebra}
  obtained from $gl(m|n)$ by regarding the $gl(m|n)$ generators as
  letters in words contained in $U[gl(m|n)]$, where the (graded)
  commutator bracket acts as a relation to reduce the algebra
  somewhat.  $U[gl(m|n)]$ is infinite dimensional, although of
  finite rank, viz, again $m+n-1$.

\item
  The \emph{quantum superalgebra} $U_q[gl(m|n)]$ is then
  a so-called `$q$-deformation'%
  \footnote{%
    We might say `quantum deformation' here, but the relation to
    quantum mechanics is more of analogy than of rigor.
  }
  of $U[gl(m|n)]$, which maintains its viability as a Hopf
  (super)algebra structure (see below) \cite{DeliusZhang:95}.  Roughly
  speaking, the deformation amounts to `exponentiation by $q$'; indeed
  $U[gl(m|n)]$ may be recovered as the limit $q\to 1$ of $U_q[gl(m|n)]$.
  $U_q[gl(m|n)]$ is of course also
  infinite dimensional, and again of rank $m+n-1$.

\end{enumerate}


\subsection{Generators and relations for $U_q[gl(m|n)]$}
\label{sec:UqglmnFullDescription}

Following Zhang \cite[pp1237-1238]{Zhang:93}, we provide a full
description of $U_q[gl(m|n)]$ in terms of generators and relations.
For various invertible $X$, we will repeatedly use the notation
$\overline{X}\triangleq X^{-1}$.


\subsubsection{$U_q[gl(m|n)]$ generators}

\enlargethispage{3\baselineskip}
Where $\mathcal{I}\triangleq\{1,\dots,m+n\}$ is the set of
the $gl(m|n)$ \emph{indices}, we define a $\mathbb{Z}_2$ grading
$[\cdot]:\mathcal{I}\to\mathbb{Z}_2$:
\begin{eqnarray*}
  [a]
  \triangleq
  \left\{
  \begin{array}{lll}
    0 \quad & \mathrm{if~} a \leqslant m \qquad\qquad
            & \mathrm{even~indices}
    \\
    1       & \mathrm{else}
            & \mathrm{odd~indices}.
  \end{array}
  \right.
\end{eqnarray*}
Throughout, we shall use dummy indices
$a,b\in\mathcal{I}$ where meaningful.
A set of ${(m+n)}^2$ generators for $U_q[gl(m|n)]$ is then:
\begin{eqnarray}
  \left\{
  \! \!
    \begin{array}{rlll}
      K_{a},       &       & \quad \qquad \qquad m+n & \mathrm{Cartan}
      \\[0.5mm]
      {E^{b}}_{a}, & a < b & \qquad \qquad
        \frac{1}{2}(m+n)(m+n-1) & \mathrm{lowering}
      \\[0.5mm]
      {E^{a}}_{b}, & a < b & \qquad \qquad
        \frac{1}{2}(m+n)(m+n-1) & \mathrm{raising}
    \end{array}
  \! \!
  \right\}.
  \label{eq:GeneratorSet}
\end{eqnarray}
Let us now introduce the notation, for any $a\in\mathcal{I}$:
\begin{eqnarray*}
  q_a
  \triangleq
  q^{{(-)}^{[a]}},
\end{eqnarray*}
For any power $N$, replacing $q$ with $q^N$ immediately shows that
$(q_a)^N=(q^N)_a$, so we may write $q_a^N$ with impunity, specifically,
we will write $\overline{q}_a\equiv q_a^{-1}$.  Next, an equivalent
notation for $K_{a}$ is $q_a^{{E^a}_a}$; where the exponential is
defined in the usual manner as an infinite sum, thus powers $K_{a}^N$
are meaningful; specifically, we will often be working with
$N\in\frac{1}{2}\mathbb{Z}$.  Thus, under the mapping $q\mapsto
\overline{q}$, $K_{a}$ is mapped to $\overline{K}_a$, where we intend
$\overline{K}_a\equiv K_a^{-1}$.  As expected, for arbitrary powers
$M,N$, we have:
\begin{eqnarray*}
  K_{a}^{M}
  K_{a}^{N}
  =
  K_{a}^{M+N}
  \quad
  \mathrm{where}
  \quad
  K_{a}^0
  \equiv
  \mathrm{Id},
  \qquad
\end{eqnarray*}
where $\mathrm{Id}$ is the $U_q[gl(m|n)]$ identity element.
Apart from $N\in \mathbb{N}$, powers (i.e. products) of the non-Cartan
generators $({E^{a}}_{b})^N$ for $a\neq b$, are not meaningful.

On the generators we define a natural $\mathbb{Z}_2$ grading in terms
of the grading on the indices:
\begin{eqnarray}
  [K_{a}^N]
  \triangleq
  0,
  \qquad \qquad
  [{E^a}_b]
  \triangleq
  [a] + [b]
  \quad
  (\mathrm{mod}\;2),
  \label{eq:Gradingofgenerators}
\end{eqnarray}
where the former may be seen as a special case of the latter by setting
$a=b$ and making the identification
$K_a\equiv q^{(-)^{[a]}{E^{a}}_{a}}$.
We use the terms ``even'' and ``odd'' for generators in the same manner
as we do for indices.  Elements of $U_q[gl(m|n)]$ are said to be
\emph{homogeneous} if they are linear combinations of generators of the
same grading.  The product $XY$ of homogeneous $X,Y\in U_q[gl(m|n)]$
has grading:
\begin{eqnarray}
  [ X Y ]
  \triangleq
  [X] + [Y]
  \quad
  (\mathrm{mod}\;2).
  \label{eq:Gradingofproduct}
\end{eqnarray}
Thus, for example, inspection of (\ref{eq:Gradingofgenerators}) and
(\ref{eq:Gradingofproduct}) shows that we may cheerfully substitute
$[{E^{a}}_{c}]$ for $[{E^{a}}_{b}{E^{b}}_{c}]$. Further, we also
have the following useful results for $a<b<c<d$:
\begin{eqnarray*}
  [{E^{a}}_{d}]
  [{E^{b}}_{c}]
  =
  [{E^{b}}_{c}]
  \qquad
  \mathrm{and}
  \qquad
  \left[{E^{a}}_{b}\right]
  [{E^{b}}_{c}]
  =
  [{E^{a}}_{b}]
  [{E^{c}}_{d}]
  =
  0.
\end{eqnarray*}


\subsubsection{$U_q[gl(m|n)]$ simple generators}

The full set of generators (\ref{eq:GeneratorSet}) includes some
redundancy; some can be regarded as \emph{simple} in that the rest may
be expressed in terms of them.  We shall call the following subset of
$3(m+n)-2$ generators the $U_q[gl(m|n)]$ \emph{simple} generators:
\begin{eqnarray*}
  \left\{
    \! \! \!
    \begin{array}{rl@{\hspace{4pt}}r@{\hspace{10pt}}l}
      K_{a},         & &
        \qquad \qquad m+n   & \mathrm{Cartan}
      \\
      {E^{a+1}}_{a}, & a < m+n &
        \qquad \qquad m+n-1 & \mathrm{simple~lowering}
      \\
      {E^{a}}_{a+1}, & a < m+n & \quad
        \qquad \qquad m+n-1 & \mathrm{simple~raising}
    \end{array}
    \! \! \!
  \right\}.
\end{eqnarray*}
The fact that there are $m+n-1$ simple lowering generators indicates
that $U_q[gl(m|n)]$ has rank $m+n-1$.  Note that there are only two
\emph{odd} simple generators: ${E^{m+1}}_{m}$ (lowering) and
${E^{m}}_{m+1}$ (raising).


\subsubsection{$U_q[gl(m|n)]$ nonsimple generators}

In the $gl(m|n)$ case, the remaining nonsimple (non-Cartan) generators
satisfy the same commutation relations as the simple generators.  The
situation is very different for $U_q[gl(m|n)]$; the nonsimple
generators do \emph{not} satisfy the same commutation relations as do
the simple generators. Instead, they are recursively defined in terms
of sums of products of the simple generators (see
\cite[p1971,~(3)]{Zhang:92} and \cite[p1238,~(2)]{Zhang:93}).  Strictly
speaking, they are not explicitly required for the definition of the
algebra; their use can help simplify otherwise large expressions.

To whit, a set of $U_q[gl(m|n)]$ nonsimple generators may be
defined recursively for $a<b$ by:
\begin{equation}
  \left.
  \begin{array}{lrcll}
    (\mathrm{a}) &
    {E^b}_a
    & \triangleq &
      {E^b}_c
      {E^c}_a
      -
      q_c
      {E^c}_a
      {E^b}_c
    &
    \mathrm{nonsimple~lowering}
    \\
    (\mathrm{b}) &
    {E^a}_b
    & \triangleq &
      {E^a}_c
      {E^c}_b
      -
      \overline{q}_c
      {E^c}_b
      {E^a}_c
    \qquad \hspace{-1pt}
    & \mathrm{nonsimple~raising},
  \end{array}
  \hspace{21pt}
  \right\}
  \label{eq:UqglmnNonSimpleGenerators}
\end{equation}
where $a<c<b$; viz $c$ is an \emph{arbitrary} index, we do \emph{not}
intend a sum here.

In \S\ref{sec:Uqglmnalphareps}, we will have use
for an alternative set of nonsimple generators, again defined
recursively for $a<b$ by:
\begin{equation}
  \hspace{-30pt}
  \left.
  \begin{array}{lrcll}
    (\mathrm{a}) &
    {\mathbf{E}^b}_a
    & \triangleq &
      {\mathbf{E}^b}_c
      {\mathbf{E}^c}_a
      -
      \overline{q}_c
      {\mathbf{E}^c}_a
      {\mathbf{E}^b}_c
    &
    \mathrm{alternative~nonsimple~lowering}
    \\
    (\mathrm{b}) &
    {\mathbf{E}^a}_b
    & \triangleq &
      {\mathbf{E}^a}_c
      {\mathbf{E}^c}_b
      -
      q_c
      {\mathbf{E}^c}_b
      {\mathbf{E}^a}_c
    \qquad
    \hspace{-13pt}
    & \mathrm{alternative~nonsimple~raising},
  \end{array}
  \hspace{15pt}
  \right\}
  \label{eq:UqglmnAlternativeNonSimpleGenerators}
\end{equation}
where we intend ${\mathbf{E}^a}_b\triangleq {E^a}_b$ when ${E^a}_b$ is
simple, viz, for any $|a-b|=1$.  Note that we use a boldface
${\mathbf{E}^a}_b$ where the original source \cite{Zhang:92} uses an
overline ${\overline{E}^a}_b$; the use of the boldface notation
saves the overline for indicating inverses.

These definitions may be written more concisely with some more notation.
Writing $S^a_b\triangleq\mathrm{sign}(a-b)$, we may replace
(\ref{eq:UqglmnNonSimpleGenerators}) and
(\ref{eq:UqglmnAlternativeNonSimpleGenerators}), for all $a\neq b$, by:
\begin{eqnarray}
  \left.
  \begin{array}{lrcl}
    (\mathrm{a}) &
    {E^a}_b
    & = &
    {E^a}_c
    {E^c}_b
    -
    q_c^{S^a_b}
    {E^c}_b
    {E^a}_c
    \\
    (\mathrm{b}) &
    {\mathbf{E}^a}_b
    & = &
    {\mathbf{E}^a}_c
    {\mathbf{E}^c}_b
    -
    \overline{q}_c^{S^a_b}
    {\mathbf{E}^c}_b
    {\mathbf{E}^a}_c.
  \end{array}
  \hspace{118pt}
  \right\}
  \label{eq:UqglmnNonSimpleGeneratorsBetterNotation}
\end{eqnarray}

The two different sets of generators are in fact Hermitian conjugates.
For all meaningful indices $a,b$, we have:
\begin{eqnarray*}
  ({E^{a}}_{b})^\dagger
  =
  {\mathbf{E}^{b}}_{a},
  \qquad
  ({\mathbf{E}^{a}}_{b})^\dagger
  =
  {E^{b}}_{a},
  \qquad
  {(K_a^N)}^\dagger
  =
  K_a^N,
\end{eqnarray*}
and these definitions ensure that $(X^\dagger)^\dagger=X$ for all
$U_q[gl(m|n)]$ generators $X$.  Note that these are ordinary, not
$\mathbb{Z}_2$ graded Hermitian conjugates, meaning that we have
$(XY)^\dagger=Y^\dagger X^\dagger$; expressly \emph{not}
$(XY)^\dagger=(-)^{[X][Y]}Y^\dagger X^\dagger$.

Lastly, we mention a result of Zhang \cite[Lemma~3]{Zhang:92}, which
gives us a more efficient formula than
(\ref{eq:UqglmnNonSimpleGeneratorsBetterNotation}b) for expanding the
alternative nonsimple generators:
\begin{eqnarray}
  {\mathbf{E}^a}_b
  =
  {E^a}_b
  +
  S^a_b
  \sum_c
    \Delta_c
    {E^c}_b
    {\mathbf{E}^{a}}_{c},
  \label{eq:Zhang92Lemma3}
\end{eqnarray}
for \emph{any} indices $a\neq b$, where the sum is over all $c$
strictly between $a$ and $b$.
(If $|a-b|=1$, then the sum is ignored, and the result is trivial.)

Note that in (\ref{eq:Zhang92Lemma3}), we have introduced the following
handy notation:
\begin{eqnarray}
  \left.
  \begin{array}{rcl}
    \Delta
    & \triangleq &
    q - \overline{q},
    \qquad
    \Delta_{a}
    \triangleq
    q_a-\overline{q}_a
    =
    {(-)}^{[a]}
    (q-\overline{q})
    =
    {(-)}^{[a]}
    \Delta
    \\
    \overline{\Delta}
    & \triangleq &
    (\Delta)^{-1}
    \qquad
    \overline{\Delta}_{a}
    \triangleq
    (\Delta_{a})^{-1}.
  \end{array}
  \hspace{26pt}
  \right\}
  \label{eq:DefnofDeltaa}
\end{eqnarray}


\subsubsection{The graded commutator}

The \emph{graded commutator}
$
  [\cdot,\cdot]
  :
  U_q[gl(m|n)] \times U_q[gl(m|n)]
  \to
  U_q[gl(m|n)]
$,
is defined for homogeneous
$ X, Y \in U_q[gl(m|n)] $ by (\ref{eq:glmnCommutatorBracket}), viz:
\begin{eqnarray*}
  [ X, Y ]
  \triangleq
  X Y - {(-)}^{[X][Y]} Y X,
\end{eqnarray*}
and extended by linearity.  For completeness, we mention that for
\emph{associative} superalgebras, of which $U_q[gl(m|n)]$ is certainly
an example, we have the following useful graded commutator identities:
\begin{eqnarray}
  \left.
  \begin{array}{lrcl}
    (\mathrm{a}) &
    [XY,Z]
    & = &
    X[Y,Z]
    +
    {(-)}^{[Y][Z]}
    [X,Z]Y
    \\
    (\mathrm{b}) &
    [X,YZ]
    & = &
    [X,Y]Z
    +
    {(-)}^{[X][Y]}
    Y[X,Z].
  \end{array}
  \hspace{73pt}
  \right\}
  \label{eq:AssociativeSAIdentity}
\end{eqnarray}


\subsubsection{$U_q[gl(m|n)]$ relations}
\label{sec:Uqglmnrelations}

With this notation, we have the following $U_q[gl(m|n)]$
\emph{relations}:

\begin{enumerate}
\item
  The Cartan generators all commute; for any powers $M,N$:
  \begin{eqnarray}
    K_{a}^{M}
    K_{b}^{N}
    =
    K_{b}^{N}
    K_{a}^{M}.
    \label{eq:CartanGensCommute}
  \end{eqnarray}

\item
  The Cartan generators commute with the simple raising and lowering
  generators in the following manner:
  \begin{eqnarray}
    K_{a} {E^b}_{b\pm1} \overline{K}_a
    =
    q_a^{(\delta^a_b - \delta^a_{b\pm1} )}
    {E^b}_{b\pm1}.
    \label{eq:UqglmnCartanRaisingCommutation}
  \end{eqnarray}
  From (\ref{eq:UqglmnCartanRaisingCommutation}), we have the following
  useful interchange:
  \begin{eqnarray}
    K_{a}
    {E^b}_{b\pm1}
    =
    q_a^{( \delta^a_b - \delta^a_{b\pm1} )}
    {E^b}_{b\pm1}
    K_{a}.
    \label{eq:Donkey}
  \end{eqnarray}
  In Lemma \ref{lem:CommutatorLemma} (proved in Appendix
  \ref{app:CommutatorLemma}), we show that (\ref{eq:Donkey}) may be
  much strengthened to:
  \begin{eqnarray}
    K_{a}^N
    {E^b}_{c}
    =
    q_a^{N (\delta^a_b-\delta^a_{c})}
    {E^b}_{c}
    K_{a}^N,
    \label{eq:FullCNCCommutator}
  \end{eqnarray}
  for \emph{any} meaningful indices $b,c$ (viz $b<c$, $b>c$, and even
  $b=c$), and any power $N$.

\item
  The non-Cartan simple generators satisfy the following commutation
  relations (this is the really interesting part!):
  \begin{eqnarray}
    [ {E^a}_{a+1}, {E^{b+1}}_b ]
    =
    \delta^a_b
    \frac{
      K_{a} \overline{K}_{a+1} - \overline{K}_a K_{a+1}
    }{
      q_a - \overline{q}_a
    }.
    \label{eq:SimpleRLSymmetricCommutator}
  \end{eqnarray}
  Alternatively, again employing the notation of
  (\ref{eq:DefnofDeltaa}), we may write this:
  \begin{eqnarray}
    [ {E^a}_{a+1}, {E^{b+1}}_b ]
    =
    \delta^a_b
    {(-)}^{[a]}
    \left[
      {(-)}^{[a]} {E^{a}}_{a}
      -
      {(-)}^{[a+1]} {E^{a+1}}_{a+1}
    \right]_q,
    \label{eq:ConvenientSimpleRLSymmetricCommutator}
  \end{eqnarray}
  where we have introduced the $q$-bracket,
  defined for various invertible $X\in U_q[gl(m|n)]$, including scalars
  (well, scalar multiples of $\mathrm{Id}$):
  \begin{eqnarray}
    {[X]}_q
    \triangleq
    \frac{q^{X} - \overline{q}^{X}}{q - \overline{q}},
    \qquad
    \mathrm{observe~that}
    \qquad
    \lim_{q \to 1}{[X]}_q
    =
    X.
    \label{eq:qbracketdefintion}
  \end{eqnarray}
  Note that in (\ref{eq:ConvenientSimpleRLSymmetricCommutator}),
  Zhang \cite{Zhang:92} replaces
  ${(-)}^{[a]} {E^{a}}_{a}-{(-)}^{[a+1]} {E^{a+1}}_{a+1}$
  with the more convenient expression $h_a$.

\pagebreak

  In fact, (\ref{eq:SimpleRLSymmetricCommutator}) generalises to a more
  useful result (proven in \cite{DeWit:2000}):
  \begin{eqnarray}
    [ {E^a}_{b}, {E^b}_{a} ]
    =
    \frac{
      K_{a} \overline{K}_{b} - \overline{K}_a K_{b}
    }{
      q_a - \overline{q}_a
    }
    =
    \overline{\Delta}_a
    (K_{a} \overline{K}_{b} - \overline{K}_a K_{b}),
    \label{eq:NonSimpleRLSymmetricCommutator}
  \end{eqnarray}
  viz:
  \begin{eqnarray*}
    [ {E^a}_{b}, {E^b}_{a} ]
    =
    {(-)}^{[a]}
    \left[
      {(-)}^{[a]} {E^{a}}_{a}
      -
      {(-)}^{[b]} {E^{b}}_{b}
    \right]_q.
  \end{eqnarray*}

  We also have, for $|a-b| > 1$, the commutations:
  \begin{eqnarray}
    \hspace{-27pt}
    {E^{a+1}}_{a}
    {E^{b+1}}_{b}
    =
    {E^{b+1}}_{b}
    {E^{a+1}}_{a}
    \quad
    \mathrm{and}
    \quad
    {E^{a}}_{a+1}
    {E^{b}}_{b+1}
    =
    {E^{b}}_{b+1}
    {E^{a}}_{a+1}.
    \label{eq:DistantSimpleGeneratorsCommute}
  \end{eqnarray}

\item
  The squares of the odd simple generators are zero:
  \begin{eqnarray*}
    {( {E^{m}}_{m+1} )}^2
    =
    {( {E^{m+1}}_{m} )}^2
    =
    0.
  \end{eqnarray*}
  In fact, we may show that this implies that the squares of
  \emph{nonsimple} odd generators are also zero:
  \begin{eqnarray}
    {( {E^{a}}_{b} )}^2
    =
    0,
    \qquad
    [a]\neq [b].
    \label{eq:SquaresofOddGeneratorsareZero}
  \end{eqnarray}

\item
  Lastly, we have the $U_q[gl(m|n)]$ \emph{Serre relations}; their
  inclusion ensures that the algebra is reduced enough to be
  \emph{simple}.  For $a\neq m$:
  \small
  \begin{eqnarray*}
    \hspace{-30pt}
    \begin{array}{l}
      {( {E^{a+1}}_{a} )}^2 {E^{a\pm1+1}}_{\hspace{-5pt}a\pm1}
      -
      \nabla
      {E^{a+1}}_{a} {E^{a\pm1+1}}_{\hspace{-5pt}a\pm1} {E^{a+1}}_{a}
      +
      {E^{a\pm1+1}}_{\hspace{-5pt}a\pm1} {( {E^{a+1}}_{a} )}^2
      =
      0
      \\
      {( {E^{a}}_{a+1} )}^2 {E^{a\pm1}}_{\hspace{-5pt}a\pm1+1}
      -
      \nabla
      {E^{a}}_{a+1} {E^{a\pm1}}_{\hspace{-5pt}a\pm1+1} {E^{a}}_{a+1}
      +
      {E^{a\pm1}}_{\hspace{-5pt}a\pm1+1} {( {E^{a}}_{a+1} )}^2
      =
      0.
    \end{array}
  \end{eqnarray*}
  \normalsize
  where to save space, we have introduced the notation:
  $\nabla\triangleq q + \overline{q}$.
  These may be more succinctly expressed using nonsimple generators.
  Noting that for $a\neq m$, we have $q_{a}=q_{a+1}$, the above
  become, again for $a\neq m$:
  \begin{eqnarray}
    \left.
    \begin{array}{crcl}
      (\mathrm{a}) &
      {E^{a+1}}_{a}
      {E^{a+2}}_{a}
      &=&
      q_{a}
      {E^{a+2}}_{a}
      {E^{a+1}}_{a}
      \\[0.5mm]
      (\mathrm{b}) &
      {E^{a}}_{a+1}
      {E^{a}}_{a+2}
      &=&
      q_{a}
      {E^{a}}_{a+2}
      {E^{a}}_{a+1}
      \\[0.5mm]
      (\mathrm{c}) &
      {E^{a+1}}_{a-1}
      {E^{a+1}}_{a}
      &=&
      q_a
      {E^{a+1}}_{a}
      {E^{a+1}}_{a-1}
      \\[0.5mm]
      (\mathrm{d}) &
      {E^{a-1}}_{a+1}
      {E^{a}}_{a+1}
      &=&
      q_a
      {E^{a}}_{a+1}
      {E^{a-1}}_{a+1}.
    \end{array}
    \hspace{60pt}
    \right\}
    \label{eq:SerreRelationsaneqm}
  \end{eqnarray}

\pagebreak

  Alternatively, if we define a \emph{$q$ graded commutator}:
  \begin{eqnarray*}
    \left[ X, Y \right]_q
    \triangleq
    X Y
    -
    (-)^{[X][Y]}
    q
    Y X,
  \end{eqnarray*}
  and, equivalently, a \emph{$\overline{q}$ graded commutator} by
  replacing $q$ with $\overline{q}$, then
  (\ref{eq:SerreRelationsaneqm}) may be more elegantly expressed:
  \begin{eqnarray}
    [X, [ X, Y]_q]_{\overline{q}} = 0,
    \label{eq:XYqgradedcommutator}
  \end{eqnarray}
  where the pair $(X,Y)$ represents the four pairs
  $({E^{a}}_{a+1}, {E^{a\pm 1}}_{a\pm 1 + 1})$
  and
  $({E^{a+1}}_{a}, {E^{a\pm 1 + 1}}_{a\pm 1})$.
  Equivalently, we may exchange $q$ and $\overline{q}$ in
  (\ref{eq:XYqgradedcommutator}).  Note that in these cases, the parity
  factors $(-)^{[X][Y]}$ in (\ref{eq:XYqgradedcommutator}) are always
  $+1$ as $a\neq m$ in ${E^{a}}_{a+1}$ and ${E^{a+1}}_{a}$, viz the
  graded commutators degenerate to ordinary commutators.

  The $a\neq m$ Serre relations are complemented by a pair dealing
  with the case $a=m$:
  \begin{eqnarray*}
    {E^{m+1}}_{m}
    {E^{m+2}}_{m-1}
    & = &
    -
    {E^{m+2}}_{m-1}
    {E^{m+1}}_{m}
    \\
    {E^{m}}_{m+1}
    {E^{m-1}}_{m+2}
    & = &
    -
    {E^{m-1}}_{m+2}
    {E^{m}}_{m+1},
  \end{eqnarray*}
  more succinctly expressed as:
  \begin{eqnarray*}
    \left[ {E^{m+1}}_{m}, {E^{m+2}}_{m-1} \right]
    =
    \left[ {E^{m}}_{m+1}, {E^{m-1}}_{m+2} \right]
    =
    0,
  \end{eqnarray*}
  where, as the generators are all odd, the graded commutators
  are read as anticommutators.

  Observe that if either $m$ or $n$ is $1$, there are actually \emph{no}
  $U_q[gl(m|n)]$ Serre relations; making life a little simpler.
\end{enumerate}

These relations tell us, in principle, how to reexpress products of
simple generators.  In general, to reexpress a product containing
\emph{non}simple generators, those nonsimple generators must first be
recursively expanded using
(\ref{eq:UqglmnNonSimpleGeneratorsBetterNotation}a), with any graded
commutators expanded by linearity, before the above relations can be
invoked.

The above description of the $U_q[gl(m|n)]$ relations should convince
the gentle reader that $U_q[gl(m|n)]$ has a formidable structure.  To
facilitate examination of its representation theory, in
\S\ref{sec:UqglmnPBW}, we will rewrite the $U_q[gl(m|n)]$ relations
into a PBW basis formulation, which is suitable for implementation on a
computer.

\pagebreak


\subsection{$U_q[gl(m|n)]$ root system}
\label{sec:Uqglmnrootsystem}

We next introduce the $U_q[gl(m|n)]$ root system, which is identical to
that of $gl(m|n)$. We will have use for it in
\S\ref{sec:TheKacinducedmoduleconstruction}, \S\ref{sec:NormalOrdering}
and \S\ref{sec:AnorthonrmalbasisB}.

Where the $gl(m|n)$ Cartan subalgebra is denoted by $H$, its dual, the
$gl(m|n)$ weight space $H^{*}$, has a basis given by the $gl(m|n)$
\emph{fundamental weights} $\{\varepsilon_a\}_{a\in\mathcal{I}}$, which
are lists of zeros of length $m+n$, with a $1$ in position $a$.
The $\varepsilon_a$ inherit a grading from that on the indices.
As $H$ and $H^{*}$ are dual, where
$\{{e^{b}}_{b}\}_{b\in\,\mathcal{I}}$ are the $gl(m|n)$ Cartan
generators, we have the form:
$
  {e^{b}}_{b} (\varepsilon_a)
  \triangleq
  \delta_{a b}
$.
On $H^*$, we have an invariant symmetric bilinear form
$(\cdot,\cdot):H^{*} \times H^{*} \to \mathbb{C}\,$, defined by:
\begin{eqnarray*}
  (\varepsilon_a, \varepsilon_b)
  \triangleq
  (-)^{[a]} \delta_{a b},
\end{eqnarray*}
and extended by linearity.

Next, to each non-Cartan $gl(m|n)$ generator ${e^{a}}_{b}$, there
corresponds a $gl(m|n)$ \emph{root}
$\alpha_{ab}\triangleq\varepsilon_a-\varepsilon_b$, which is the weight
of ${e^{a}}_{b}$ in the adjoint representation.%
\footnote{%
  We apologise for the overloading of $\alpha$. The notation in this
  subsection will go no further.
}
For our purposes, it is convenient to bastardise the notation. Also
permitting $a=b$, we will refer to $\alpha_{ab}$ as the `weight' of
${e^{a}}_{b}$:
\begin{eqnarray}
  \mathrm{wt} ({e^{a}}_{b})
  \triangleq
  \varepsilon_{a} - \varepsilon_{b}
  =
  \alpha_{ab},
  \label{eq:AdjointActionWeight}
\end{eqnarray}
indicating that within a $gl(m|n)$ weight module (see
\S\ref{sec:HighestweightUqglmnrepresentations}), the action of a
generator ${e^{a}}_{b}$ sends a vector of weight $\gamma$ to another of
weight $\gamma+\alpha_{ab}$.

The roots inherit a grading from the indices:
$[\alpha_{ab}]\equiv [{e^{a}}_{b}]$.  Further, we assign signs to them
in accordance with those of these generators, viz that if ${e^{a}}_{b}$
is a lowering generator, the corresponding root $\alpha_{ab}$ is said
to be \emph{negative}, written $\alpha_{ab}\prec 0$, and if
${e^{a}}_{b}$ is a raising generator, then $\alpha_{ab}$ is said to be
\emph{positive}, written $\alpha_{ab}\succ 0$.  To illustrate, weights
for the $U_q[gl(3|1)]$ lowering generators are presented in Table
\ref{tab:31generatorweights}.

\begin{table}[ht]
  \begin{centering}
  \begin{tabular}{
      c
      l@{\hspace{0pt}}r@{\hspace{2pt}}r@{\hspace{2pt}}r
        @{\hspace{4pt}}c@{\hspace{4pt}}r@{\hspace{0pt}}r
      c
    }
    $X$ & \multicolumn{7}{c}{$\mathrm{wt}(X)$}
    \\[1mm]
    \cline{1-8}
    \\[-3mm]
    ${E^{4}}_{3}$ & $($ & $ 0$, & $ 0$, & $-1$ & $|$ & $+1$ & $)$\\
    ${E^{4}}_{2}$ & $($ & $ 0$, & $-1$, & $ 0$ & $|$ & $+1$ & $)$&
      odd\\
    ${E^{4}}_{1}$ & $($ & $-1$, & $ 0$, & $ 0$ & $|$ & $+1$ & $)$\\[1mm]
    \cline{1-8}
    \\[-3mm]
    ${E^{3}}_{2}$ & $($ & $ 0$, & $-1$, & $+1$ & $|$ & $ 0$ & $)$\\
    ${E^{3}}_{1}$ & $($ & $-1$, & $ 0$, & $+1$ & $|$ & $ 0$ & $)$&
      even\\
    ${E^{2}}_{1}$ & $($ & $-1$, & $+1$, & $ 0$ & $|$ & $ 0$ & $)$
  \end{tabular}
  \caption{
    Weights (all negative) of the $U_q[gl(3|1)]$ lowering generators.
  }
  \label{tab:31generatorweights}
  \end{centering}
\end{table}

Using this notation, $gl(m|n)$ has the following \emph{simple},
positive roots:
\begin{eqnarray*}
  \alpha_{a,a+1}
  \triangleq
  \varepsilon_{a}
  -
  \varepsilon_{a+1},
  \qquad
  a = 1, \dots, m+n-1.
\end{eqnarray*}
Apart from the single odd simple root $\alpha_{m,m+1}$, the simple
positive roots are all even.  (Of various choices for Lie superalgebra
root systems, this \emph{distinguished root system} is unique in
containing only one odd simple root.)

Then, we define $\Delta^+_{i}$ to be the set of $gl(m|n)$ positive roots
of grading $i$, and $\Delta^+$ to be the union of the $\Delta^+_i$,
viz:%
\footnote{%
  We apologise for the overloading of $\Delta$.  In practice, this
  $\Delta$ will only appear with a positive superscript, so it is
  easily distinguishable.
}
\begin{eqnarray*}
  \Delta^+_{i}
  =
  \{
    \gamma
    :
    [\gamma]=i
    \mathrm{~and~}
    \gamma\succ 0
  \},
  \qquad
  \qquad
  \qquad
  \Delta^+
  =
  \Delta^+_{0}
  \cup
  \Delta^+_{1},
\end{eqnarray*}
where $[\gamma]$ denotes the grading of the root $\gamma$.

In terms of these, we define the half sums of all positive even and odd
$gl(m|n)$ roots, and their graded sum $\rho$:
\begin{eqnarray*}
  \rho_{i}
  =
  {\textstyle \frac{1}{2}}
  \sum_{\gamma\in\Delta^+_i}
    \gamma,
  \qquad
  \qquad
  \qquad
  \rho
  =
  {\textstyle \frac{1}{2}}
  \sum_{\gamma\in\Delta^+}
    {(-)}^{[\gamma]}
    \gamma
  =
  \sum_{i=0,1}
    {(-)}^{i}
    \rho_{i},
\end{eqnarray*}
viz $\rho=\rho_0-\rho_1$.
Specifically, for $gl(m|n)$ (and hence for $U_q[gl(m|n)]$), we find
\cite[p6207]{DeliusGouldLinksZhang:95b}:%
\footnote{%
  We make a correction to
  \cite{DeliusGouldLinksZhang:95b}, which appears to cite an error
  reproduced several times before and after, e.g. it appears in
  Zhang \cite{Zhang:93}. To whit, the term ``$2m$'' in the formula
  for $\rho$ repeatedly appears as ``$m$''.
}
\begin{eqnarray*}
  \rho_{0}
  &=&
  {\textstyle \frac{1}{2}}
  \sum_{a=1}^{m}
    \left(
      m-2a+1
    \right)
    \varepsilon_{a}
  +
  {\textstyle \frac{1}{2}}
  \sum_{a=m+1}^{m+n}
    \left(
      m+n-2a+1
    \right)
    \varepsilon_{a}
  \\
  \rho_{1}
  &=&
  {\textstyle \frac{1}{2}}
  \sum_{a=1}^{m}
    n
    \varepsilon_{a}
  -
  {\textstyle \frac{1}{2}}
  \sum_{a=m+1}^{m+n}
    m
    \varepsilon_{a}
  \\
  \rho
  &=&
  {\textstyle \frac{1}{2}}
  \sum_{a=1}^{m}
    \left(
      m-n-2a+1
    \right)
    \varepsilon_{a}
  +
  {\textstyle \frac{1}{2}}
  \sum_{a=m+1}^{m+n}
    \left(
      2m+n-2a+1
    \right)
    \varepsilon_{a}.
\end{eqnarray*}


\subsection{$U_q[gl(m|n)]$ as a Hopf superalgebra}
\label{sec:UqglmnasaHopfsuperalgebra}

$U_q[gl(m|n)]$ may be regarded as a Hopf superalgebra when equipped
with the following (compatible!) coproduct $\Delta$, counit
$\varepsilon$ and antipode $S$ structures. The material is taken from
\cite[p1238]{Zhang:93}, except that we have modified the definition of
the coproduct and antipode so that they have increased symmetry.  This
material is included for completeness; in \S\ref{sec:ThesubmodulesVk},
we will only have need for the coproduct.  $U_q[gl(m|n)]$ is in fact a
\emph{quasitriangular Hopf superalgebra}, i.e. it possesses a
(universal) R matrix.

We first introduce some notation. A homomorphism or antihomomorphism
$H$ on $U_q[gl(m|n)]$ is described as
\emph{($\mathbb{Z}_2$) graded} if it is compatible with the graded
commutator, viz:
\begin{eqnarray*}
  H([X,Y])
  =
  [H(X),H(Y)],
\end{eqnarray*}
where the latter graded commutator may even exist on an ungraded space,
e.g. $\mathbb{C}\,$, where it is actually trivial.

\enlargethispage{\baselineskip}
This means that a graded homomorphism $H$ and a graded
antihomomorphism $A$ must necessarily satisfy:
\begin{eqnarray*}
  H(XY)
  =
  H(X) H(Y)
  \qquad
  \mathrm{and}
  \qquad
  A(XY)
  =
  {(-)}^{[X][Y]}
  A(Y) A(X),
\end{eqnarray*}
of which only the latter varies from the usual, ungraded situation.


\subsubsection{Coproduct $\Delta$}

The coproduct (a.k.a. comultiplication) is a
$\mathbb{Z}_2$ graded algebra homomorphism
$
  \Delta : U_q[gl(m|n)] \to U_q[gl(m|n)] \otimes U_q[gl(m|n)]
$,
defined by:
\begin{eqnarray}
  \left.
    \begin{array}{lrcl}
      (\mathrm{a}) &
      \Delta ({E^{a+1}}_{a})
      & = &
      {E^{a+1}}_{a}
      \otimes
      \overline{K}_a^{\frac{1}{2}} K_{a+1}^{\frac{1}{2}}
      +
      K_{a}^{\frac{1}{2}} \overline{K}_{a+1}^{\frac{1}{2}}
      \otimes
      {E^{a+1}}_{a}
      \\[1mm]
      (\mathrm{b}) &
      \Delta ({E^{a}}_{a+1})
      & = &
      {E^{a}}_{a+1}
      \otimes
      \overline{K}_a^{\frac{1}{2}} K_{a+1}^{\frac{1}{2}}
      +
      K_{a}^{\frac{1}{2}}
      \overline{K}_{a+1}^{\frac{1}{2}}
      \otimes {E^{a}}_{a+1}
      \\[1mm]
      (\mathrm{c}) &
      \Delta (K_{a})
      & = &
      K_{a} \otimes K_{a},
    \end{array}
    \hspace{14pt}
  \right\}
  \label{eq:UqglmnCoproduct}
\end{eqnarray}
and extended to all of $U_q[gl(m|n)]$ by:
\begin{eqnarray}
  \Delta(X Y)
  =
  \Delta(X)
  \Delta(Y),
  \qquad
  \mathrm{for~all~}
  X,Y \in U_q[gl(m|n)].
  \label{eq:Deltaisahomo}
\end{eqnarray}
Observe that $\Delta$ preserves grading, viz that $[\Delta(X)]=[X]$ for
homogeneous $X \in U_q[gl(m|n)]$, where we have
$[X \otimes Y]\triangleq [X]+[Y]$.

By substitution of $q^{N}$ for $q$ in (\ref{eq:UqglmnCoproduct}c), we
discover:
\begin{eqnarray*}
  \Delta (K_{a}^{N})
  =
  K_{a}^{N} \otimes K_{a}^{N},
\end{eqnarray*}
and, setting $N=0$, hence
$\Delta(\mathrm{Id})=\mathrm{Id}\otimes\mathrm{Id}$, as expected.

Before proceeding, we mention that our definition in
(\ref{eq:UqglmnCoproduct}) is only one of various possibilities; we
have chosen it for its symmetry.  In fact, in comparison with the
literature, our $\Delta$ agrees with that of
\cite{GeGouldZhangZhou:98a}, and differs from that of Zhang
\cite{Zhang:93} and that of my PhD thesis \cite{DeWit:98}.  We
mention that given \emph{any} coproduct, it is possible to
write down another coproduct structure, the
``opposite coproduct'': $\Delta^T\triangleq T\cdot\Delta$.
Here, the \emph{twist map} $T$ is an operator on the tensor product
$U_q[gl(m|n)]\otimes U_q[gl(m|n)]$, defined for homogeneous
$X,Y\in U_q[gl(m|n)]$ by:
\begin{eqnarray*}
  T ( X \otimes Y )
  =
  {(-)}^{[X][Y]}
  ( Y \otimes X ).
\end{eqnarray*}

More relevant to our purposes here, we may extend the expression for
the coproduct for simple generators to that for nonsimple generators.
Firstly, as in (\ref{eq:UqglmnNonSimpleGeneratorsBetterNotation}),
writing $S^a_b\triangleq\mathrm{sign}(a-b)$, we may cheerfully rewrite
(\ref{eq:UqglmnCoproduct}) for the simple generators ${E^{a}}_{b}$, for
any $|a-b|=1$:
\begin{eqnarray}
  \Delta({E^a}_b)
  =
  {E^a}_b
  \otimes
  K_{a}^{\frac{1}{2}S^a_b}
  \overline{K}_b^{\frac{1}{2}S^a_b}
  +
  \overline{K}_{a}^{\frac{1}{2}S^a_b}
  K_{b}^{\frac{1}{2}S^a_b}
  \otimes
  {E^a}_b.
  \label{eq:SimpleFullCoproduct}
\end{eqnarray}
Using this notation, we prove in Lemma \ref{lem:FullCoproduct} in
Appendix \ref{app:CoproductLemma} the following more general statement,
for \emph{any} valid indices $a,b$:
\begin{eqnarray}
  \hspace{-20pt}
  \Delta({E^a}_b)
  & = &
  {E^a}_b
  \otimes
  K_{a}^{\frac{1}{2}S^a_b}
  \overline{K}_b^{\frac{1}{2}S^a_b}
  +
  \overline{K}_a^{\frac{1}{2}S^a_b}
  K_{b}^{\frac{1}{2}S^a_b}
  \otimes
  {E^a}_b
  -
  \nonumber
  \\
  & &
  \qquad
  S^a_b
  {\displaystyle \sum_{c}}
    \Delta_c
    \left(
      \overline{K}_c^{\frac{1}{2}S^a_b}
      K_{b}^{\frac{1}{2}S^a_b}
      {E^a}_c
      \otimes
      {E^c}_b
      K_{a}^{\frac{1}{2}S^a_b}
      \overline{K}_c^{\frac{1}{2}S^a_b}
    \right),
  \label{eq:FullCoproduct}
\end{eqnarray}
where the sum ranges over all $c$ strictly between $a$ and $b$, and is
simply ignored if $|a-b|\leqslant 1$. Where $a=b$, the statement is
also true; this is made clear when the equivalence
$K_a=q^{(-)^{[a]}{E^a}_a}$ is noted.

Lastly, we apologise for even further overloading the definition of
$\Delta$; to be sure, the coproduct will only appear with parentheses
enclosing its argument.


\subsubsection{Counit $\varepsilon$}

The counit $\varepsilon:U_q[gl(m|n)]\to\mathbb{C}\,$,
is also a $\mathbb{Z}_2$ graded algebra homomorphism, defined by:
\begin{eqnarray*}
  \varepsilon ( {E^{a\pm1}}_{a} )
  =
  \varepsilon ( {E^{a}}_{a\pm1} )
  =
  0,
  \qquad \qquad
  \varepsilon ( K_{a} )
  =
  1,
\end{eqnarray*}
and extended to all of $U_q[gl(m|n)]$ by
$\varepsilon(XY)=\varepsilon(X)\varepsilon(Y)$.  Again, we have
$\varepsilon (K_{a}^{N})=1$, and, setting $N=0$, thus
$\varepsilon(\mathrm{Id})=1$, as expected.  We apologise for
overloading the definition of $\varepsilon$ as the counit with the
$gl(m|n)$ fundamental weights (see \S\ref{sec:Uqglmnrootsystem}). As we
shall have no further use for the counit, we are safe.


\subsubsection{Antipode $S$}

Lastly, the antipode $S:U_q[gl(m|n)]\to U_q[gl(m|n)]$, is a
$\mathbb{Z}_2$ graded algebra \emph{antiauto}morphism, defined by:
\begin{eqnarray*}
  S ({E^{a+1}}_{a})
  & = &
  -
  K_{a}^{\frac{1}{2}}
  K_{a+1}^{\frac{1}{2}}
  {E^{a+1}}_{{a}}
  \\
  S ({E^{a}}_{{a+1}})
  & = &
  -
  \overline{K}_{a}^{\frac{1}{2}}
  \overline{K}_{a+1}^{\frac{1}{2}}
  {E^{a}}_{{a+1}}
  \\
  S (K_{a})
  & = &
  \overline{K}_a,
\end{eqnarray*}
and extended to all of $U_q[gl(m|n)]$ by:
\begin{eqnarray*}
  S ( X Y )
  =
  {(-)}^{[X] [Y]}
  S ( Y ) S ( X ),
\end{eqnarray*}
for homogeneous $X,Y\in U_q[gl(m|n)]$.  Again, immediately
$S(K_{a}^{N})=\overline{K}_a^{N}$, and thus
$S(\mathrm{Id})=\mathrm{Id}$, as expected.

$S$ is perhaps better expressed in terms of the notation introduced for
the coproduct $\Delta$. We have, for simple generators ${E^{a}}_{b}$,
where $|a-b|=1$:
\begin{eqnarray}
  S ( {E^{a}}_{b} )
  =
  -
  K_{a}^{\frac{1}{2} S^a_b}
  K_{b}^{\frac{1}{2} S^a_b}
  {E^{a}}_{{b}}.
  \label{eq:Antipode}
\end{eqnarray}
This result is also valid for the case $a=b$ (where
$S^a_a=\mathrm{sign}(a-a)=0$), and the formula degenerates to
$S({E^{a}}_{a})=-K_a^{0}K_a^{0}{E^{a}}_{a}=-{E^{a}}_{a}$, which is
equivalent to $S(K_a)=\overline{K}_a$.  Furthermore, a direct
inductive proof%
\footnote{%
  An example proof, albeit for a different definition of $S$, is
  provided in Lemma 3 of \cite{Zhang:92}.
}
shows that (\ref{eq:Antipode}) generalises to the case of nonsimple
generators ${E^{a}}_{b}$, so that (\ref{eq:Antipode}) is, in a sense,
the most general expression of $S$.


\pagebreak


\section{Highest weight $U_q[gl(m|n)]$ representations $\pi_{\lambda}$}
\label{sec:HighestweightUqglmnrepresentations}

\subsection{Introduction}

The construction of highest weight representations for $U_q[gl(m|n)]$
involves initially postulating a \emph{highest weight vector}, which we
shall call $\ket{1}$. The action of the $U_q[gl(m|n)]$ Cartan
generators on $\ket{1}$ is that of scalar multipliers; the details of
these multiplications are encoded in the \emph{weight} of $\ket{1}$.
Thus, if we are dealing with a representation labeled:
\begin{eqnarray}
  \lambda
  \equiv
  (
    \lambda_1, \dots, \lambda_m
    \; | \;
    \lambda_{m+1}, \dots, \lambda_{m+n}
  )
  =
  \sum_{i=1}^{m+n}
    \lambda_i
     \varepsilon_i,
  \label{eq:glmnFundamentalWeights}
\end{eqnarray}
where the $\varepsilon_i$ are $gl(m|n)$ fundamental weights (see
\S\ref{sec:Uqglmnrootsystem}), then we intend $\ket{1}$ to have weight
$\lambda$, that is the action of $K_{a}$ on $\ket{1}$ is:%
\footnote{
  For consistency of the weight notation between $gl(m|n)$ and
  $U_q[gl(m|n)]$, (\ref{eq:Kaactiononket1}) tells us that $\lambda_a$
  is actually the weight of $\ket{1}$ in terms of the $U_q[gl(m|n)]$
  `generators' ${E^a}_a$.
}
\begin{eqnarray}
  K_{a}
  \cdot
  \ket{1}
  \equiv
  \pi_{\lambda} (K_{a})
  \cdot
  \ket{1}
  \triangleq
  q_a^{\lambda_{a}}
  \ket{1}.
  \label{eq:Kaactiononket1}
\end{eqnarray}
Substituting $q^N$ for $q$, we immediately have that:
\begin{eqnarray}
  K_{a}^{N}
  \cdot
  \ket{1}
  =
  q_a^{N \lambda_{a}}
  \ket{1}.
  \label{eq:Cartanactiononket1}
\end{eqnarray}
We implement the notion that $\ket{1}$ is a \emph{highest} weight
vector by declaring that it be annihilated by the actions of
\emph{all raising} generators:
\begin{eqnarray}
  {E^{a}}_{b}
  \cdot
  \ket{1}
  \triangleq
  0,
  \qquad
  a < b.
  \label{eq:raisingactiononket1}
\end{eqnarray}
The module $V_{\lambda}$ is then defined by the action of
all possible products of the $U_q[gl(m|n)]$ lowering generators on
$\ket{1}$. We may determine a basis $B_{\lambda}$ for
$V_{\lambda}$, with elements $\ket{i}$ defined by:
\begin{eqnarray*}
  \ket{i}
  \triangleq
  \beta_{i}
  \;
  X_{i_1} X_{i_2} \cdots X_{i_p}
  \cdot
  \ket{1},
\end{eqnarray*}
where the $X_{i_j}$ are $U_q[gl(m|n)]$ generators, $p$ is the number of
generators in the product, and $\beta_{i}$ is a normalisation
constant.

We call $B_{\lambda}$ a \emph{graded weight} basis, meaning that we may
assign to it (i.e. to $V_{\lambda}$) a grading consistent with that of
$U_q[gl(m|n)]$, and a system of \emph{weights} (both, see
\S\ref{sec:Uqglmnalphareps}). For our specific choices of $\lambda$
(again, see \S\ref{sec:Uqglmnalphareps}), $V_{\lambda}$ is
finite-dimensional.


\subsection{The Kac induced module construction}
\label{sec:TheKacinducedmoduleconstruction}

The \emph{Kac induced module construction} (KIMC) is a two-stage
process which efficiently implements the construction of
$B_{\lambda}$.

\begin{itemize}
\item
  Firstly, we construct a basis $B^0_{\lambda}$ for the so-called `even
  subalgebra submodule' $V^0_{\lambda}\subset V_{\lambda}$; this
  being the module of highest weight $\lambda$ of the $U_q[gl(m|n)]$
  `even subalgebra' $U_q[gl(m)\oplus gl(n)]$, viz the algebra generated
  by the even generators of $U_q[gl(m|n)]$. That is, $V^0_{\lambda}$ is
  defined by the action of all possible combinations of \emph{even
  lowering} generators on $\ket{1}$, where we have declared that
  $\ket{1}$ is annihilated by the action of all \emph{even raising}
  generators.

\item
  Secondly, $V_{\lambda}$ is induced from $V_{\lambda}^0$ by the
  repeated action of the odd lowering generators on $V_{\lambda}^0$,
  subject to the proviso that $V_{\lambda}^0$ is annihilated by the
  (unique) odd raising generator, to whit:
  ${E^{m}}_{m+1} \cdot V_{\lambda}^0\triangleq\{0\}$. This implies that
  ${E^a}_b\cdot V_{\lambda}^0 = \{0\} $ for \emph{all} odd raising
  generators ${E^a}_b$, $a<b$. Thus, we construct $B_{\lambda}$ from
  $B^0_{\lambda}$.

  In this process, there is a subtlety: the resultant ``Kac module''
  (i.e. $V_{\lambda}$) may not be irreducible. However, we shall choose
  $\lambda$ such that $\pi_{\lambda}$ is a so-called \emph{typical}
  representation \cite{Kac:78,Zhang:93}, ensuring that $V_{\lambda}$
  \emph{is} irreducible.

\end{itemize}


\subsection{Dimension of $V_{\lambda}$}
\label{sec:DimensionofVlambda}

For arbitrary typical highest weight $U_q[gl(m|n)]$ representations
$V_{\lambda}$, we have the following \emph{Kac--Weyl dimension formula}
\cite{Kac:78}:
\begin{eqnarray}
  \mathrm{dim}(V_{\lambda})
  =
  2^{mn}
  \cdot
  \mathrm{dim}(V^0_{\lambda}),
  \qquad
  \mathrm{where}
  \quad
  \mathrm{dim}(V^0_{\lambda})
  =
  \prod_{\gamma \in \Delta^+_0}
    \frac{
      (\lambda+\rho_0,\gamma)
    }{
      (\rho_0,\gamma)
    },
  \label{eq:KacWeylDimension}
\end{eqnarray}
where $\Delta^+_0$, $\rho_0$ and the inner product
$(\cdot,\cdot)$ on the $gl(m|n)$ fundamental weights are presented in
\S\ref{sec:Uqglmnrootsystem}.

For the specific choice $\lambda=\Lambda=(\dot{0}_m|\dot{\alpha}_n)$,
for even positive roots $\gamma$, we have:
\begin{eqnarray*}
  (\Lambda,\gamma)
  =
  {\textstyle \frac{1}{2}}
  \sum_{a<b,[a]=[b]}
    (\Lambda, \varepsilon_a-\varepsilon_b)
  =
  {\textstyle \frac{1}{2}}
  \sum_{a<b,[a]=[b]}
    (-)^{[a]}
    (\Lambda_a-\Lambda_b)
  =
  0,
\end{eqnarray*}
as $\Lambda_a=\Lambda_b$ for $[a]=[b]$, thus
$\mathrm{dim}(V^0_{\Lambda})=1$, hence
$\mathrm{dim}(V_{\Lambda})=2^{mn}$, which simplifies things.  Details
of the KIMC for this case are presented in
\S\ref{sec:Uqglmnalphareps}.


\subsection{Matrix elements}

To construct explicit \emph{matrix elements} $\pi_{\Lambda} (X)$ for a
particular $U_q[gl(m|n)]$ generator $X$, the action of $X$ on
\emph{each} of the basis vectors of $B_{\lambda}$ must be determined.

Whilst the action of the generators on $\ket{1}$ is predefined, more
generally, the determination of the action of $X$ on an arbitrary
vector $\ket{i}$ requires the rendering of a string of generators into
a normal ordering and the application of the `KIMC rules' to simplify
that normally ordered expression into (a multiple of) a basis vector.

Thus, we must first determine an appropriate ordering (see
\S\ref{sec:NormalOrdering}), and then describe an appropriate set of
generator commutations to implement that ordering (viz the PBW lemma of
\S\ref{sec:Commutations}). By the latter, we mean that we intend not
to use all the commutators of \S\ref{sec:Uqglmnrelations} directly.
Instead, we shall use expressions taken from lemmas in
\cite{Zhang:92,Zhang:93} for the commutations between non-Cartan
generators.

With these tools, we proceed to build bases and explicit matrix elements
for our particular representations in \S\ref{sec:Uqglmnalphareps}.


\section{A normal ordering and a PBW lemma}
\label{sec:UqglmnPBW}

Finding a normal ordering for a string of $U_q[gl(m|n)]$ algebra
generators involves the recursive use of commutation relations to
rewrite the string as a sum of strings, with respect to some chosen
(hopefully natural) ordering. Both the initial string and the resultant
may contain initial scalar multipliers, which for $U_q[gl(m|n)]$ are
typically algebraic expressions in $q$. When we speak of the length
of a string, we shall ignore these scalars.

A PBW lemma describes the appropriate commutations, but we must
determine an ordering ourselves.  Perhaps the most natural ordering is
purely by weight (see (\ref{eq:AdjointActionWeight}) in
\S\ref{sec:Uqglmnrootsystem}), but there are reasons for choosing other
orderings.


\subsection{A normal ordering for $U_q[gl(m|n)]$}
\label{sec:NormalOrdering}

We begin with the convention that if generators $G_1$ and $G_2$ are
ordered, viz $G_1\leqslant G_2$, then the string $G_1G_2$ is ordered.
With this, the ordering we choose is based on the following
principles:

\begin{enumerate}
\item
  Our string will often be regarded as (right) acting on the highest
  weight vector $\ket{1}\in V_{\lambda}$, and the KIMC directs us to
  first build an even subalgebra submodule $V_{\lambda}^{0}$ (see
  \S\ref{sec:TheKacinducedmoduleconstruction} and
  \S\ref{sec:AnorthonrmalbasisB}) based on this $\ket{1}$, i.e. to
  define basis vectors of $V_{\lambda}^0$ in terms of the right actions
  of strings of even lowering generators on $\ket{1}$.  Thus, we
  require even generators to be greater than odd generators, i.e.
  within normally-ordered strings, even generators lie to the right of
  odd ones.

\item
  Within the even generators, $\ket{1}$ is always annihilated by the
  (right) action of raising generators, so these must be greatest, i.e.
  rightmost. By symmetry, we then demand that the least amongst the
  even generators are the lowering generators, so the Cartan must be
  lie between the even lowering and the even raising.

\item
  Within the odd generators, the (right) action of the raising
  generators always annihilates any vectors from $V_{\lambda}^{0}$, so
  the odd lowering must be lesser than the odd raising.

\item
  Within the five equivalence classes created by these considerations,
  non-Cartan generators are ordered by increasing weight.
  Doing this ensures that squares of odd generators can be
  systematically identified and annihilated; it also facilitates a
  systematic way of defining basis vectors for $V_{\lambda}$ (see
  \S\ref{sec:Uqglmnalphareps}).

  Furthermore, (powers of) Cartan generators are ordered by index, viz
  $K_{a}^{M}\leqslant K_{b}^{N}$ if $a\leqslant b$.  Doing this ensures
  that powers of the same generator may be combined.

\end{enumerate}

\pagebreak

We call this ordering ``$OL<OR<EL<C<ER$''.  It differs slightly from
that (implicitly) described in \cite[p1240]{Zhang:93}, viz
$OR<OL<EL<C<ER$. To implement it, we say that distinct \emph{weights}
$\gamma_1$ and $\gamma_2$ are ordered (viz $\gamma_1\prec\gamma_2$) if
the first nonzero component of $\gamma_1-\gamma_2$ is positive. Then,
say that we are comparing generators $G_1$ and $G_2$; where $\gamma_i$
is the weight of $G_i$ (see (\ref{eq:AdjointActionWeight}) in
\S\ref{sec:Uqglmnrootsystem}); $L_i$ is the `lifting' of $G_i$, being
$-1$, $0$ or $+1$ if $G_i$ is a lowering, Cartan or raising generator,
respectively; and, if $G_i$ is Cartan, then let $a_i$ be its index (the
exponent is unimportant).  Then:
\begin{eqnarray*}
  A
  \vee
  (
    B
    \wedge
    (
      C
      \vee
      (
        D
        \wedge
        (
          (E\wedge F)
          \vee
          G
        )
      )
    )
  )
  \;
  \Longleftrightarrow
  \;
  G_1 \leqslant G_2,
\end{eqnarray*}
where:
\begin{eqnarray*}
  \hspace{-11pt}
  \begin{array}{rcll}
    A
    & \mathrm{is} &
    [G_1]>[G_2]
    \quad
    &
    G_1 \mathrm{~is~odd~and~} G_2 \mathrm{~even}
    \\
    B
    & \mathrm{is} &
    [G_1]=[G_2]
    &
    \mathrm{both~odd~or~both~even}
    \\
    C
    & \mathrm{is} &
    L_1<L_2
    &
    \mathrm{ordered~liftings}
    \\
    D
    & \mathrm{is} &
    L_1=L_2
    &
    \mathrm{same~liftings}
    \\
    E
    & \mathrm{is} &
    L_1=0
    &
    G_1 \mathrm{~is~Cartan}
    \\
    F
    & \mathrm{is} &
    a_1\leqslant a_2
    &
    \mathrm{(implicitly)~both~Cartan~and~ordered}
    \\
    G
    & \mathrm{is} &
    \gamma_1\preccurlyeq\gamma_2
    &
    \mbox{ordered~within~(implicitly~non-Cartan)~class}.
  \end{array}
\end{eqnarray*}


To illustrate the ordering, for $U_q[gl(3|1)]$, with reference to
Table \ref{tab:31generatorweights}, we have:
\begin{eqnarray*}
  \hspace{-18pt}
  \begin{array}{l}
    \underbrace{{E^4}_3 < {E^4}_2 < {E^4}_1}_\mathrm{Odd~Lowering}
    <
    \underbrace{{E^1}_4 < {E^2}_4 < {E^3}_4}_\mathrm{Odd~Raising}
    <
    \\
    \hspace{35pt}
    \underbrace{{E^3}_2 < {E^2}_1 < {E^3}_1}_\mathrm{Even~Lowering}
    <
    \underbrace{
      K_{1}^{N_1} < K_{2}^{N_2} < K_{3}^{N_3}
    }_\mathrm{Cartan}
    <
    \underbrace{{E^1}_3 < {E^1}_2 < {E^2}_3}_\mathrm{Even~Raising}.
  \end{array}
\end{eqnarray*}

This ordering ensures that in the KIMC, the action of a
normally-ordered string of generators on a highest weight vector
$\ket{1}$ may be evaluated by the following \emph{ordered} steps.

\begin{enumerate}
\item
  If there are any terminal even raising generators, then the string
  evaluates to $0$.

\item
  If there are terminal Cartan generators, then these may be replaced
  by their known scalar actions on $\ket{1}$, and the string is reduced
  in length.

\item
  Next, the action of any even lowering generators is considered. In
  the general situation, these map $\ket{1}$ to another basis vector of
  $V_{\lambda}^{0}$. For our particular modules $V_{\Lambda}$, where
  $\Lambda=(\dot{0}_m|\dot{\alpha}_n)$, as $V^0_{\Lambda}$ is one
  dimensional, $\ket{1}$ is in fact \emph{annihilated} by the even
  lowering generators, so if there are \emph{any} even raising
  generators, then the string evaluates to $0$.

\item
  Next, the odd raising generators annihilate any vectors of
  $V_{\lambda}^0$, so if any are present, then the string evaluates to
  $0$.

\item
  Lastly, when the string is reduced to (a scalar multiple of) the
  action of some unrepeated odd lowering generators on a
  $V_{\lambda}^0$ basis vector, that residual string may be identified
  as (a scalar multiple of) a particular $V_{\lambda}$ basis vector.

\end{enumerate}

We mention that although our PBW lemma provides us with the
\emph{means} to normally order generator strings, the normal ordering
is a computationally expensive process.  Firstly, \emph{each} exchange
may generate up to two extra terms in a sum,%
\footnote{%
  Although the exchanges can in fact add up to two extra terms, in
  practice they add only one extra term, but they can be `sum-neutral'
  or even subtract a term.
}
so there is a geometric increase in the number of terms with
exchanges.  Secondly, implementation of the exchanges is really a
sorting procedure, but we have not been able to implement an efficient
algorithm -- we in fact use the dumbest possible opportunistic
exchange. This failure is partly due to the complexities in developing
a sorting algorithm in the presence of the continual creation of extra
terms.

Thus, the process to normally order generator strings requires time and
storage which both of which grow at least exponentially with string
length.  Using \textsc{Mathematica}, we currently get into serious
trouble beyond length $8$.


\subsection{Commutations implementing the normal ordering}
\label{sec:Commutations}

To implement the normal ordering described in
\S\ref{sec:NormalOrdering}, we describe here a set of
generator-exchanging commutations.  The material originates in
\cite{Zhang:92,Zhang:93}; we have modified the results a little in
light of (\ref{eq:FullCNCCommutator}), rearranged many things, and
corrected several minor mistakes. In what follows, we intend distinct
abstract indices to represent different concrete indices.


\subsubsection{A PBW commutator lemma}

The following result contains some corrections to the original
\cite{Zhang:93}. In it, we use the notation presented in
(\ref{eq:DefnofDeltaa}).

\pagebreak

\begin{lemma}
  ~\\
  We have the following commutations. Firstly,
  (\ref{eq:SimpleRLSymmetricCommutator})
  generalises to the case of nonsimple generators
  (\ref{eq:NonSimpleRLSymmetricCommutator}), viz:
  \begin{eqnarray}
    \left[{E^{a}}_{b},{E^{b}}_{a}\right]
    =
    \overline{\Delta}_a
    (
      K_{a} \overline{K}_{b}
      -
      \overline{K}_{a} K_{b}
    )
    \qquad
    \mathit{all~}
    a, b.
    \label{eq:TomWaits}
  \end{eqnarray}
  Secondly, where there are three distinct indices, we have:
  \begin{eqnarray}
    \left[{E^{a}}_{c},{E^{c}}_{b}\right]
    & = &
    \left\{
      \begin{array}{lll}
        (\mathrm{a}) & K_{c} \overline{K}_{b} {E^{a}}_{b} & c<b<a \\
        (\mathrm{b}) & {E^{a}}_{b} K_{a} \overline{K}_{c} & c<a<b \\
        (\mathrm{c}) & {E^{a}}_{b} K_{c} \overline{K}_{a} & b<a<c \\
        (\mathrm{d}) & K_{b} \overline{K}_{c} {E^{a}}_{b} & a<b<c
      \end{array}
    \right\}
    \label{eq:BrendanPerry}
    \\
    \left[{E^{c}}_{a},{E^{c}}_{b}\right]
    & = &
    \left[{E^{a}}_{c},{E^{b}}_{c}\right]
    =
    0
    \label{eq:AustralianCrawl}
    \qquad
    a<c<b
    \mathit{~~or~~}
    b<c<a
    \\
    {E^{c}}_{a}
    {E^{c}}_{b}
    & = &
    \left\{
      \begin{array}{lll}
        (\mathrm{a}) &
        {(-)}^{[{E^{c}}_{b}]} & a<b<c
        \\
        (\mathrm{b}) &
        {(-)}^{[{E^{c}}_{a}]} & c<a<b
      \end{array}
    \right\}
    q_{c}
    {E^{c}}_{b}
    {E^{c}}_{a}
    \label{eq:LeonardCohen}
    \\
    {E^{a}}_{c}
    {E^{b}}_{c}
    & = &
    \left\{
      \begin{array}{lll}
        (\mathrm{a}) &
        {(-)}^{[{E^{b}}_{c}]} & a<b<c
        \\
        (\mathrm{b}) &
        {(-)}^{[{E^{a}}_{c}]} & c<a<b
      \end{array}
    \right\}
    q_{c}
    {E^{b}}_{c}
    {E^{a}}_{c}.
    \label{eq:LisaGerrard}
  \end{eqnarray}
  Thirdly, we describe the situation where there are no common indices,
  and we have $a<b$ and $c<d$.  Let $S(\mathsf{x},\mathsf{y})$ denote
  the set of integers $\{\mathsf{x},\mathsf{x}+1,\dots,\mathsf{y}\}$.
  Then, if $S(a,b)$ and $S(c,d)$ are either disjoint or one is wholly
  contained within the other, viz $a<c<d<b$, $a<b<c<d$, $c<a<b<d$ or
  $c<d<a<b$, we have a total of $16$ cases:
  \begin{eqnarray}
    \left[{E^{a}}_{b},{E^{c}}_{d}\right]
    =
    \left[{E^{a}}_{b},{E^{d}}_{c}\right]
    =
    \left[{E^{b}}_{a},{E^{c}}_{d}\right]
    =
    \left[{E^{b}}_{a},{E^{d}}_{c}\right]
    =
    0.
    \label{eq:DrHooks16cases}
  \end{eqnarray}
  More interestingly, if there is some other overlap between the sets
  $S(a,b)$ and $S(c,d)$, viz $a<c<b<d$ or $c<a<d<b$, then we have the
  $8$ cases:
  \begin{eqnarray}
    \left[{E^{a}}_{b},{E^{c}}_{d}\right]
    & = &
    \left\{
      \begin{array}{lll}
        (\mathrm{a}) &
        +
        \Delta_{b}
        &
        a<c<b<d
        \\
        (\mathrm{b}) &
        -
        \Delta_{d}
        &
        c<a<d<b
      \end{array}
    \right\}
    {E^{a}}_{d}
    {E^{c}}_{b}
    \label{eq:RyCooder1}
    \\
    \left[{E^{b}}_{a},{E^{d}}_{c}\right]
    & = &
    \left\{
      \begin{array}{lll}
        (\mathrm{a}) &
        +
        \Delta_{b}
        &
        a<c<b<d
        \\
        (\mathrm{b}) &
        -
        \Delta_{d}
        &
        c<a<d<b
      \end{array}
    \right\}
    {E^{d}}_{a}
    {E^{b}}_{c}
    \label{eq:RyCooder2}
    \\
    \left[{E^{a}}_{b},{E^{d}}_{c}\right]
    & = &
    \left\{
      \begin{array}{lll}
        (\mathrm{a}) &
        -
        \Delta_{b}
        \overline{K}_{b}
        K_{c}
        {E^{a}}_{c}
        {E^{d}}_{b}
        &
        a<c<b<d
        \\
        (\mathrm{b}) &
        +
        \Delta_{d}
        {E^{d}}_{b}
        {E^{a}}_{c}
        \overline{K}_{a}
        K_{d}
        &
        c<a<d<b
      \end{array}
    \right\}
    \label{eq:VishwaMohanBhatt}
    \\
    \left[{E^{b}}_{a},{E^{c}}_{d}\right]
    & = &
    \left\{
      \begin{array}{lll}
        (\mathrm{a}) &
        -
        \Delta_{c}
        {E^{b}}_{d}
        {E^{c}}_{a}
        \overline{K}_{c}
        K_{b}
        &
        a<c<b<d
        \\
        (\mathrm{b}) &
        +
        \Delta_{a}
        \overline{K}_{d}
        K_{a}
        {E^{c}}_{a}
        {E^{b}}_{d}
        &
        c<a<d<b
      \end{array}
    \right\}.
    \label{eq:RaviShankar}
  \end{eqnarray}
  \label{lem:Lemma2ofZhang93}
\end{lemma}

\pagebreak

Rearranging the indices in Lemma \ref{lem:Lemma2ofZhang93} gives us the
following simplified results.
\begin{itemize}
\item
  The entirety of (\ref{eq:BrendanPerry}) may be summarised by:
  \begin{eqnarray*}
    \left[{E^{a}}_{c},{E^{c}}_{b}\right]
    =
    \left\{
      \begin{array}{ll}
        K_{c}
        \overline{K}_{b}
        {E^{a}}_{b}
        &
        c<b<a
        \\
        {E^{a}}_{b}
        K_{a}
        \overline{K}_{c}
        &
        c<a<b
        \\
        {E^{a}}_{b}
        K_{c}
        \overline{K}_{a}
        &
        b<a<c
        \\
        K_{b}
        \overline{K}_{c}
        {E^{a}}_{b}
        &
        a<b<c.
      \end{array}
    \right.
  \end{eqnarray*}

\item
  The entirety of (\ref{eq:AustralianCrawl}) to (\ref{eq:LisaGerrard})
  may be summarised by:
  \begin{eqnarray*}
    {E^{a}}_{c}
    {E^{b}}_{c}
    =
    \kappa
    {E^{b}}_{c}
    {E^{a}}_{c}
    \qquad
    \mathrm{and}
    \qquad
    {E^{c}}_{a}
    {E^{c}}_{b}
    =
    \kappa
    {E^{c}}_{b}
    {E^{c}}_{a},
  \end{eqnarray*}
  where:
  \begin{eqnarray*}
    \kappa
    \triangleq
    \left\{
    \begin{array}{ll}
      1
      &
      \mathrm{if~} z(a,b,c) = c
      \\
      {(-)}^{[{E^{z(a,b,c)}}_c]}
      \overline{q}_{c}^{S^a_b} & \mathrm{else},
    \end{array}
    \right.
  \end{eqnarray*}
  and $z(a,b,c)$ picks out the middle element of $\{a,b,c\}$.
  (The $1$ factor follows as $[{E^{a}}_{c}][{E^{b}}_{c}]=0$ for $c$
  between $a$ and $b$.)

\item
  The entirety of (\ref{eq:DrHooks16cases}) to (\ref{eq:RaviShankar})
  may be summarised by:
  \begin{eqnarray*}
    \left[{E^{a}}_{b},{E^{c}}_{d}\right]
    =
    \left\{
      \begin{array}{l@{\hspace{0pt}}ll}
        +
        \Delta_{b} &
        {E^{a}}_{d}
        {E^{c}}_{b}
        &
        a<c<b<d
        \\
        -
        \Delta_{d} &
        {E^{a}}_{d}
        {E^{c}}_{b}
        &
        c<a<d<b
        \\
        +
        \Delta_{a} &
        {E^{c}}_{b}
        {E^{a}}_{d}
        &
        b<d<a<c
        \\
        -
        \Delta_{c} &
        {E^{c}}_{b}
        {E^{a}}_{d}
        &
        d<b<c<a
        \\
        -
        \Delta_{b} &
        \overline{K}_{b}
        K_{d}
        {E^{a}}_{d}
        {E^{c}}_{b}
        &
        a<d<b<c
        \\
        +
        \Delta_{c} &
        {E^{c}}_{b}
        {E^{a}}_{d}
        \overline{K}_{a}
        K_{c}
        &
        d<a<c<b
        \\
        -
        \Delta_{c} &
        {E^{a}}_{d}
        {E^{c}}_{b}
        \overline{K}_{c}
        K_{a}
        &
        b<c<a<d
        \\
        +
        \Delta_{b} &
        \overline{K}_{d}
        K_{b}
        {E^{c}}_{b}
        {E^{a}}_{d}
        &
        c<b<d<a
        \\
        & 0 &
        a\neq b\neq c\neq d \mathrm{~~else}.
      \end{array}
    \right.
  \end{eqnarray*}

\end{itemize}

From these, we deduce the following rules for exchanges:

\pagebreak

\begin{description}


\item[From (\ref{eq:CartanGensCommute}), replace $K_{a}^M K_{b}^N$ with]
      $K_{b}^N K_{a}^M$.

  If also $a=b$, then replace it with $K_{a}^{M+N}$.

  \qquad If also $M+N=0$ then replace it with $\mathrm{Id}$.


\item[From (\ref{eq:FullCNCCommutator}), replace $K_{a}^N {E^{b}}_{c}$
      with]
  $q_a^{N (\delta^a_b - \delta^a_c)} {E^{b}}_{c} K_{a}^N$, and

\hspace{29pt} \textbf{replace} ${E^{b}}_{c} K_{a}^N$ \textbf{with}
  $ q_a^{-N (\delta^a_b - \delta^a_c)} K_{a}^N {E^{b}}_c$.


\item[From (\ref{eq:NonSimpleRLSymmetricCommutator}),
      replace ${E^{a}}_{b} {E^{b}}_{a}$ with]
  $
    {(-)}^{[{E^{a}}_{b}]} {E^{b}}_{a} {E^{a}}_{b}
    +
    \overline{\Delta}_{a}
    (\textstyle K_{a} \overline{K}_{b} - \overline{K}_{a} K_{b})
  $.


\item[From (\ref{eq:SquaresofOddGeneratorsareZero}),
      replace ${E^{a}}_{b} {E^{a}}_{b}$ with]
  $0$ if $[{E^{a}}_{b}]=1$.


\item[Replace ${E^{a}}_{c} {E^{c}}_{b}$ with:]
  \begin{eqnarray*}
    \hspace{-20pt}
    \begin{array}{r@{\hspace{5pt}}c@{\hspace{5pt}}l@{\hspace{0pt}}l@{\hspace{20pt}}l}
      {E^a}_b & +
      & q_c^{S^a_b} & {E^c}_b {E^a}_c
      &
      a<c<b,\;\;\;b<c<a
      \\
        {(-)}^{[{E^c}_b]} {E^c}_b {E^a}_c & +
         & & K_{c} \overline{K}_b {E^{a}}_{b}
      &
      c<b<a
      \\
        {(-)}^{[{E^a}_c]} {E^c}_b {E^a}_c & +
         & & {E^{a}}_{b} K_{a} \overline{K}_c
      &
      c<a<b
      \\
        {(-)}^{[{E^a}_c]} {E^c}_b {E^a}_c & +
         & & {E^{a}}_{b} K_{c} \overline{K}_a
      &
      b<a<c
      \\
        {(-)}^{[{E^c}_b]} {E^c}_b {E^a}_c & +
         & & K_{b} \overline{K}_c {E^{a}}_{b}
      &
      a<b<c.
    \end{array}
  \end{eqnarray*}

\item[Replace ${E^c}_b {E^a}_c$ with:]
  \begin{eqnarray*}
    \hspace{-20pt}
    \begin{array}{r@{\hspace{5pt}}l@{\hspace{5pt}}c@{\hspace{5pt}}l@{\hspace{20pt}}l}
      q_c^{S^b_a} & ({E^a}_c {E^c}_b & - & {E^a}_b)
      &
      a<c<b,\;\;\;b<c<a
      \\
      {(-)}^{[{E^c}_b]}
      & ({E^a}_c {E^c}_b & - & K_{c} \overline{K}_b {E^a}_b)
      &
      c<b<a
      \\
      {(-)}^{[{E^a}_c]}
      & ({E^a}_c {E^c}_b & - & {E^a}_b K_{a} \overline{K}_c)
      &
      c<a<b
      \\
      {(-)}^{[{E^a}_c]}
      & ({E^a}_c {E^c}_b & - & {E^a}_b K_{c} \overline{K}_a)
      &
      b<a<c
      \\
      {(-)}^{[{E^c}_b]}
      & ({E^a}_c {E^c}_b & - & K_{b} \overline{K}_c {E^a}_b)
      &
      a<b<c.
    \end{array}
  \end{eqnarray*}

\item[Replace ${E^c}_a {E^c}_b$ with]
  $\kappa {E^c}_b {E^c}_a$, and
  ~\\
  \textbf{replace ${E^a}_c {E^b}_c$ with}
  $\kappa {E^b}_c {E^a}_c$, where:
  \begin{eqnarray*}
    \kappa
    \triangleq
    \left\{
    \begin{array}{l@{\hspace{20pt}}l}
      1
      &
      \mathrm{if~} z(a,b,c) = c
      \\
      {(-)}^{[{E^{z(a,b,c)}}_c]}
      \overline{q}_{c}^{S^a_b} & \mathrm{else},
    \end{array}
    \right.
  \end{eqnarray*}
  and where $z(a,b,c)$ picks out the middle element of $\{a,b,c\}$.


\item[Replace ${E^a}_b {E^c}_d$ with]
  $(-)^{[{E^a}_b] [{E^c}_d]} {E^c}_d {E^a}_b + T$, where:
  \begin{eqnarray*}
    T
    =
    \left\{
      \begin{array}{l@{\hspace{0pt}}l@{\hspace{20pt}}l}
        +
        \Delta_{b} &
        {E^{a}}_{d}
        {E^{c}}_{b}
        &
        a<c<b<d
        \\
        -
        \Delta_{d} &
        {E^{a}}_{d}
        {E^{c}}_{b}
        &
        c<a<d<b
        \\
        +
        \Delta_{a} &
        {E^{c}}_{b}
        {E^{a}}_{d}
        &
        b<d<a<c
        \\
        -
        \Delta_{c} &
        {E^{c}}_{b}
        {E^{a}}_{d}
        &
        d<b<c<a
        \\
        -
        \Delta_{b} &
        \overline{K}_{b}
        K_{d}
        {E^{a}}_{d}
        {E^{c}}_{b}
        &
        a<d<b<c
        \\
        +
        \Delta_{c} &
        {E^{c}}_{b}
        {E^{a}}_{d}
        \overline{K}_{a}
        K_{c}
        &
        d<a<c<b
        \\
        -
        \Delta_{c} &
        {E^{a}}_{d}
        {E^{c}}_{b}
        \overline{K}_{c}
        K_{a}
        &
        b<c<a<d
        \\
        +
        \Delta_{b} &
        \overline{K}_{d}
        K_{b}
        {E^{c}}_{b}
        {E^{a}}_{d}
        &
        c<b<d<a
        \\
        & 0 &
        a\neq b\neq c\neq d \mathrm{~~else}.
      \end{array}
     \right.
  \end{eqnarray*}

\end{description}

\pagebreak


\section{The $U_q[gl(m|n)]$ representations
         $(\dot{0}_m|\dot{\alpha}_n)$}
\label{sec:Uqglmnalphareps}

Fixing $m$ and $n$, in this section we describe the use of a version of
the Kac induced module construction (KIMC, see
\S\ref{sec:TheKacinducedmoduleconstruction}) in the brute-force
construction of the $U_q[gl(m|n)]$ representation
$\Lambda\triangleq(\dot{0}_m|\dot{\alpha}_n)$.

Alternatively, we might have implemented the results presented in
\cite{PalevStoilovaVanderJeugt:94,PalevTolstoy:91}, which describe the
use of a Gel'fand--Tsetlin basis to explicitly construct the actions
for \emph{essentially typical} representations%
\footnote{%
  It turns out that these results sometimes hold for other
  representations, when various limits are evaluated using
  L'H\^{o}pital formulae.
}
(this class includes our representation, which is actually
\emph{typical}).  We avoid those fine results because we wish our code
to be more general, but we pay a price for this in the currency of
computational expense.

Strictly, this material applies only to \emph{generic} $q$, that is $q$
not a root of unity (in which case the representation theory changes
drastically). Also, our representations are \emph{unitary} only under
some constraints on $\alpha$ (viz that $\alpha$ is real and either
$\alpha>n-1$ or $\alpha<1-m$, see \cite{DeliusGouldLinksZhang:95b}),
and so we shall implicitly select these. In the application of
our results to the computation of link invariants \cite{DeWit:99e}, the
representation of the braid generator based on our quantum R matrix
$\check{R}^{m,n}$ again contains the variables $q$ and $\alpha$.
However, it turns out that $\check{R}^{m,n}$ is actually a valid braid
generator for \emph{any} $q$ and $\alpha$.


\subsection{An orthonormal basis $B$ for $V\equiv V_{\Lambda}$}
\label{sec:AnorthonrmalbasisB}

Recall from \S\ref{sec:TheKacinducedmoduleconstruction} and
\S\ref{sec:DimensionofVlambda} that $V\equiv V_{\Lambda}$ is of
dimension $2^{mn}$, and may be equipped with a grading compatible with
that of $U_q[gl(m|n)]$.  It is known that $V$ contains no
\emph{weight multiplicities} (that is, $V$ contains no constant weight
subspaces of dimension greater than $1$), so that a weight basis for
$V$ will contain no distinct vectors of the \emph{same} weight, and
this makes our task a little simpler.

Here, we use a version of the KIMC to construct a weight basis
$B=\{\ket{i}\}_{i=1}^{2^{mn}}$ for $V$.  That is, the $2^{mn}$ basis
vectors $\ket{i}$ are defined in terms of the actions of all $2^{mn}$
possible nonrepeated, ordered combinations of the $U_q[gl(m|n)]$
simple lowering generators (there are $mn$ of them) on a postulated
highest weight vector $\ket{1}$.  This $\ket{1}$ is further defined to
be of unit length and annihilated by \emph{all} $U_q[gl(m|n)]$ raising
generators.  Each of the vectors defined in this manner will be
orthogonal to all other such vectors, and it is a straightforward
matter to select constants to orthonormalise them.
The resulting $B$ is thus a graded orthonormal weight basis for $V$.
Using it, in \S\ref{sec:Matrixelements}, we construct matrix elements
for the $U_q[gl(m|n)]$ generators. In each subsection, we shall
illustrate our results using $U_q[gl(3|1)]$.

\pagebreak


\subsubsection{Details of the KIMC}

The KIMC firstly instructs us to construct a basis $B^0$ for $V^0$,
the submodule of $V$ determined by the action of the $U_q[gl(m|n)]$
even subalgebra on $\ket{1}$.  For our choice of $\Lambda$, in fact
$V^0$ has dimension $1$ (see \S\ref{sec:DimensionofVlambda}), hence
$B^0=\{\ket{1}\}$.

This means that $\ket{1}$ is annihilated not only by all \emph{raising}
generators, by also by all \emph{even} lowering generators.  More
generally, for other representations, we have to work harder to
construct a basis for $V^0$; that process has a similar appearance
to the following.

Secondly, $B$ is induced from $B^0$ by the actions of all
possible products of \emph{odd lowering} generators on $B^0$.
Thus, we must consider the set of all possible products of odd
lowering generators.  The PBW lemma allows us to reduce this set to
that of all possible \emph{ordered} products, and the knowledge that
the square of odd generators is zero allows us to reduce it to the set
of all possible \emph{nonrepeated} ordered products, a \emph{finite}
set, ensuring that $B$ is finite.  Recalling that $U_q[gl(m|n)]$ has
$mn$ (simple and nonsimple) odd lowering generators:
\begin{eqnarray}
  OLGS
  =
  \left\{
    {E^{a}}_{b}:
     a = m+1, \dots, m+n,
     b = 1, \dots, m
  \right\},
  \label{eq:OLGS}
\end{eqnarray}
thus $V$ is spanned by a set of vectors obtained by the actions of all
possible ordered products (that is, strings of length $0$ to $mn$)
of $OLGS$ generators on $\ket{1}$.  Indeed, this is the source of the
factor $2^{mn}$ in the dimension formula (\ref{eq:KacWeylDimension}).
As $V$ is known to have no weight multiplicities, this spanning set is
itself the desired weight basis $B$.  The set of ordered products of
generators may be obtained from the power set $\mathcal{P}(OLGS)$, by
replacing its elements with respective products -- \textsc{Mathematica}
is well-suited to this.

To illustrate, for the $U_q[gl(3|1)]$ case, we have
$
  OLGS
  =
  \{
    {E^{4}}_{3},
    {E^{4}}_{2},
    {E^{4}}_{1}
  \},
$
and $V$ has the following basis $B$:
\begin{eqnarray}
  \left(
  \begin{array}{c}
      \ket{1}
      \\
      \ket{2}
      \\
      \ket{3}
      \\
      \ket{4}
      \\
      \ket{5}
      \\
      \ket{6}
      \\
      \ket{7}
      \\
      \ket{8}
    \end{array}
  \right)
  =
  \left(
    \begin{array}{r}
      \beta_{1} \mathrm{Id}
      \\
      \beta_{2} {E^{4}}_{3}
      \\
      \beta_{3} {E^{4}}_{2}
      \\
      \beta_{4} {E^{4}}_{1}
      \\
      \beta_{5} {E^{4}}_{3} {E^{4}}_{2}
      \\
      \beta_{6} {E^{4}}_{3} {E^{4}}_{1}
      \\
      \beta_{7} {E^{4}}_{2} {E^{4}}_{1}
      \\
      \beta_{8} {E^{4}}_{3} {E^{4}}_{2} {E^{4}}_{1}
    \end{array}
  \right)
  \cdot
  \ket{1},
  \label{eq:3ketoriginlist}
\end{eqnarray}
where the $\beta_{i}$ are scaling factors, which we shall select to
normalise the $\ket{i}$; the redundant $\beta_{1}$ is implicitly $1$.

\pagebreak


\subsubsection{Weights, gradings and an ordering for $B$}
\label{sec:WeightsandgradingsforB}

Our $B$ is a \emph{weight} basis, in that each of the $\ket{i}$ is
of a definite \emph{weight}. If we have, for some $i$:
\begin{eqnarray}
  \ket{i}
  =
  \beta_{i}
  {E^{a_1}}_{b_1}
  {E^{a_2}}_{b_2}
  \cdots
  {E^{a_p}}_{b_p}
  \cdot
  \ket{1},
  \label{eq:ketdefinition}
\end{eqnarray}
for some $p\leqslant mn$, then we may define:
\begin{eqnarray*}
  \mathrm{wt} ( \ket{i} )
  =
  \Lambda
  +
  \sum_{j=1}^p
    \mathrm{wt} ({E^{a_j}}_{b_j}),
\end{eqnarray*}
where $\mathrm{wt} ({E^{a}}_{b})$ is the weight of the generator
${E^a}_b$ (see (\ref{eq:AdjointActionWeight}) in
\S\ref{sec:Uqglmnrootsystem}).  As we defined $\ket{1}$ to be a highest
weight vector, clearly we intend $\mathrm{wt}(\ket{1})=\Lambda$.

Further, our $B$ is $\mathbb{Z}_2$ \emph{graded}, in that its
elements are formed from the actions of products of graded
$U_q[gl(m|n)]$ generators on the (zero) graded $\ket{1}$. The
$\mathbb{Z}_2$ grading of $\ket{i}$ (as defined in
(\ref{eq:ketdefinition})) is defined by:
\begin{eqnarray*}
  [\ket{1}]
  \triangleq
  0,
  \qquad
  [\ket{i}]
  \triangleq
  \sum_{j=1}^p
    [{E^{a_j}}_{b_j}]
  =
  p
  \quad
  (\mathrm{mod}\;2),
\end{eqnarray*}
where the latter result holds as the ${E^{a_j}}_{b_j}$ are all odd.
This $\mathbb{Z}_2$ grading on $V$ is compatible with a notion of
\emph{$\mathbb{Z}$ graded level}, this being $p$, the number of factors
in the product forming $\ket{i}$. This notion is relevant in the
calculation -- we recursively form the basis vectors in level $l$ by
the action of the $OLGS$ on the basis vectors of level $l-1$, for
$l=1,\dots,mn$. We number our vectors $\ket{i}$ by decreasing weight
within increasing $\mathbb{Z}$ graded levels.  This ordering is
important in that it simplifies the process of identifying an arbitrary
string acting on $\ket{1}$, which is required in
\S\ref{sec:Matrixelements}.

The weights and gradings of the basis vectors for our
$U_q[gl(3|1)]$ example are supplied in Table \ref{tab:3wtsgradings}
(cf. Table \ref{tab:31generatorweights}).

\begin{table}[ht]
  \begin{centering}
  \begin{tabular}{
      r
      l@{\hspace{0pt}}r@{\hspace{2pt}}r@{\hspace{2pt}}r
        @{\hspace{4pt}}c@{\hspace{4pt}}l@{\hspace{0pt}}r
      c
      c
    }
    $i$ & \multicolumn{7}{c}{$\mathrm{wt}(\ket{i})$} & $[\ket{i}]$ &
      ${[\ket{i}]}_{\mathbb{Z}}$
    \\[1mm]
    \cline{1-9}
    \\[-3mm]
    $1$ & $($ & $ 0$, & $ 0$, & $ 0$ & $|$ & $\alpha  $ & $)$& $0$&$0$
    \\[1mm]
    \cline{1-9}
    \\[-3mm]
    $2$ & $($ & $ 0$, & $ 0$, & $-1$ & $|$ & $\alpha+1$ & $)$& $1$\\
    $3$ & $($ & $ 0$, & $-1$, & $ 0$ & $|$ & $\alpha+1$ & $)$& $1$&$1$\\
    $4$ & $($ & $-1$, & $ 0$, & $ 0$ & $|$ & $\alpha+1$ & $)$& $1$
    \\[1mm]
    \cline{1-9}
    \\[-3mm]
    $5$ & $($ & $ 0$, & $-1$, & $-1$ & $|$ & $\alpha+2$ & $)$& $0$\\
    $6$ & $($ & $-1$, & $ 0$, & $-1$ & $|$ & $\alpha+2$ & $)$& $0$&$2$\\
    $7$ & $($ & $-1$, & $-1$, & $ 0$ & $|$ & $\alpha+2$ & $)$& $0$
    \\[1mm]
    \cline{1-9}
    \\[-3mm]
    $8$ & $($ & $-1$, & $-1$, & $-1$ & $|$ & $\alpha+3$ & $)$& $1$&$3$
  \end{tabular}
  \caption{
    Weights and gradings of the basis vectors $\ket{i}$ of $B$, for
    $U_q[gl(3|1)]$, ordered by decreasing weight
    within increasing $\mathbb{Z}$ graded levels
    ${[\ket{i}]}_{\mathbb{Z}}$.
  }
  \label{tab:3wtsgradings}
  \end{centering}
\end{table}


\subsubsection{Normalisation of $B$}

To investigate questions of orthogonality, we require an inner product
on $V$.  To whit, we introduce a basis $B^*=\{\bra{i}\}_{i=1}^{2^{mn}}$
of $V^*$ (the dual of $V$), by:
\begin{eqnarray}
  \bra{i}
  \triangleq
  \beta_{i}^{*}
  \bra{1}
  \cdot
  {\mathbf{E}^{b_p}}_{a_p}
  \cdots
  {\mathbf{E}^{b_2}}_{a_2}
  {\mathbf{E}^{b_1}}_{a_1},
  \label{eq:bradefinition}
\end{eqnarray}
where $\beta_{i}^*$ is the complex conjugate of $\beta_{i}$, and
$\ket{i}$ is as supplied in (\ref{eq:ketdefinition}).  Note that here,
we explicitly intend the Hermitian conjugates ${\mathbf{E}^{b}}_{a}$
(see (\ref{eq:UqglmnNonSimpleGeneratorsBetterNotation})), and not
${E^{b}}_{a}$. This ensures that $B^*$ is conjugate to $B$, viz we have
$\ket{i}^\dagger=\bra{i}$, and that conjugate generators and conjugate
basis vectors remain conjugates in their matrix representations.  We
assign weights and $\mathbb{Z}_2$ (and also $\mathbb{Z}$) gradings to
the $\bra{i}$ such that
$\mathrm{wt}(\bra{i})\triangleq\mathrm{wt}(\ket{i})$ and
$[\bra{i}]\triangleq [\ket{i}]$.

Using this conjugate basis, we define an inner product on $V$:
\begin{eqnarray}
  (\ket{i}, \ket{j})
  \triangleq
  \bra{i} \cdot \ket{j}
  \equiv
  \braket{i}{j},
  \label{eq:innerproductonketsdefn}
\end{eqnarray}
where we implicitly have $\braket{1}{1}=1$. Next, where ${E^b}_a$ is
any raising generator (viz $b<a$), taking the conjugate of
${E^b}_a\cdot\ket{1}=0$ yields $\bra{1}\cdot{\mathbf{E}^a}_b=0$.
Expanding the ${\mathbf{E}^a}_b$ into simple generators using
(\ref{eq:UqglmnNonSimpleGeneratorsBetterNotation}b) shows that
$\bra{1}\cdot{E^{c+1}}_{c}=0$, for all \emph{simple} lowering
generators ${E^{c+1}}_{c}$, which in turn, recursively, yields
$\bra{1}\cdot{E^c}_d=0$ for \emph{all} lowering generators ${E^c}_d$
with $c>d$. In sum, the equivalent of (\ref{eq:raisingactiononket1})
is:
\begin{eqnarray}
  \bra{1}
  \cdot
  {E^a}_b
  =
  0,
  \qquad
  a > b.
  \label{eq:loweringactiononbra1}
\end{eqnarray}
For completeness, we mention the (left) action of Cartan generators on
$\bra{1}$, obtained by conjugating (\ref{eq:Cartanactiononket1}):
\begin{eqnarray}
  \bra{1}
  \cdot
  K_{a}^{N}
  =
  q_a^{N \lambda_{a}}
  \bra{1}.
  \label{eq:Cartanactiononbra1}
\end{eqnarray}

More generally, the value of an inner product $\braket{i}{j}$ may be
calculated by the following procedure.

\begin{enumerate}
\item
  We substitute the definitions of $\bra{i}$ and $\ket{j}$
  (viz the appropriate versions of (\ref{eq:bradefinition}) and
  (\ref{eq:ketdefinition})) into (\ref{eq:innerproductonketsdefn}),
  yielding a form
  $
    \beta_{i}^{*}
    \beta_{j}
    \bra{1}
    \cdot
    Z
    \cdot
    \ket{1}
  $,
  where $Z$ is a string of generators.  For our representations
  $\pi_{\Lambda}$, as $V^0$ is one-dimensional, $Z$ will contain no
  even generators.

\item
  We use the PBW lemma to normally order $Z$.  The resulting
  strings have their raising generators annihilating $\ket{1}$ pushed
  to their right hand ends and lowering generators annihilating
  $\bra{1}$ pushed to their left hand ends.

\item
  Implementing those annihilations, and evaluating the residual
  (scalar) Cartan generator actions on $\ket{1}$ and $\bra{1}$
  (viz (\ref{eq:Cartanactiononket1}) and (\ref{eq:Cartanactiononbra1})),
  we convert the resulting expression to a scalar.
\end{enumerate}

\pagebreak

It turns out that vectors $\ket{i}$ and $\ket{j}$ with distinct weights
satisfy $\braket{i}{j}=0$, viz
$\braket{i}{j}=\delta^i_j\braket{i}{i}$.  Recall that our
$(\dot{0}_m|\dot{\alpha}_n)$ representations have no
\emph{weight multiplicities}, thus distinct basis vectors have distinct
weights, hence our basis $B$ is orthogonal.  To make it
ortho\emph{normal}, we must select the $\beta_i$ appropriately. This
means that for each $i$, we must ensure that:
\begin{eqnarray*}
  \braket{i}{i}
  =
  \beta_{i}^{*}
  \beta_{i}
  \bra{1}
  \cdot
  {\mathbf{E}^{b_p}}_{a_p}
  \cdots
  {\mathbf{E}^{b_2}}_{a_2}
  {\mathbf{E}^{b_1}}_{a_1}
  \cdot
  {E^{a_1}}_{b_1}
  {E^{a_2}}_{b_2}
  \cdots
  {E^{a_p}}_{b_p}
  \cdot
  \ket{1}
  =
  1.
\end{eqnarray*}
Thus, for each $i$, we must use the commutations of the PBW lemma to
normal order the following string $Z_i$:
\begin{eqnarray*}
  Z_i
  =
  {\mathbf{E}^{b_p}}_{a_p}
  \cdots
  {\mathbf{E}^{b_2}}_{a_2}
  {\mathbf{E}^{b_1}}_{a_1}
  \cdot
  {E^{a_1}}_{b_1}
  {E^{a_2}}_{b_2}
  \cdots
  {E^{a_p}}_{b_p},
\end{eqnarray*}
into $\mathrm{NO}(Z_i)$, and then apply the algebra-module actions
(\ref{eq:Cartanactiononket1}),
(\ref{eq:raisingactiononket1}), (\ref{eq:loweringactiononbra1}),
(\ref{eq:Cartanactiononbra1}) and $\braket{1}{1}=1$ to
$\bra{1}\cdot\mathrm{NO}(Z_i)\cdot\ket{1}$ to yield, up to an arbitrary
complex constant (phase factor),
$\beta_i=\left(\bra{1}\cdot\mathrm{NO}(Z_i)\cdot\ket{1}\right)^{-1/2}$.
The phase factor is unimportant; different choices simply lead to bases
related by orthogonal transformations, and this will not affect our R
matrices. In practice, we let the internal machinery of
\textsc{Mathematica} decide on phase factors for us -- a human
calculator might make more elegant choices.

So, at this stage, we have determined the constants $\beta_{i}$ such
that we have an orthonormal basis $B$ for $V$. In general,
for arbitrary representations of $U_q[gl(m|n)]$, these constants will
be algebraic functions of the complex variable $q$. For our particular
$\Lambda$, these functions will also contain the variable $\alpha$.

To illustrate, for the $U_q[gl(3|1)]$ case, we have:
\begin{eqnarray}
  \left.
  \begin{array}{rcl@{\hspace{2pt}}l}
    \beta_{2}
    & = &
    & {[\alpha]}_q^{-\frac{1}{2}}
    \\
    \beta_{3}
    & = &
    \overline{q} &
    {[\alpha]}_q^{-\frac{1}{2}}
    \\
    \beta_{4}
    & = &
    \overline{q}^{2} &
    {[\alpha]}_q^{-\frac{1}{2}}
    \\
    \beta_{5}
    & = &
    \overline{q} &
    {[\alpha]}_q^{-\frac{1}{2}}
    {[\alpha+1]}_q^{-\frac{1}{2}}
    \\
    \beta_{6}
    & = &
    \overline{q}^{2} &
    {[\alpha]}_q^{-\frac{1}{2}}
    {[\alpha+1]}_q^{-\frac{1}{2}}
    \\
    \beta_{7}
    & = &
    \overline{q}^{3} &
    {[\alpha]}_q^{-\frac{1}{2}}
    {[\alpha+1]}_q^{-\frac{1}{2}}
    \\
    \beta_{8}
    & = &
    \overline{q}^{3} &
    {[\alpha]}_q^{-\frac{1}{2}}
    {[\alpha+1]}_q^{-\frac{1}{2}}
    {[\alpha+2]}_q^{-\frac{1}{2}}.
  \end{array}
  \hspace{119pt}
  \right\}
  \label{eq:3betalisting}
\end{eqnarray}
where we have use the $q$-bracket (see (\ref{eq:qbracketdefintion})) to
simplify the expressions. Thus, for example,
\begin{eqnarray*}
  \beta_{8}
  =
  \left(
    \frac{
      q^{6}
      (q^{\alpha}   - \overline{q}^{\alpha}  )
      (q^{\alpha+1} - \overline{q}^{\alpha+1})
      (q^{\alpha+2} - \overline{q}^{\alpha+2})
    }{
      (q-\overline{q})^3
    }
  \right)^{-\frac{1}{2}}.
\end{eqnarray*}

Our use of the $q$-bracket notation is more than cosmetic; the
$U_q[gl(m|n)]$ symmetries manifest themselves \emph{naturally} in
$\pi_{\Lambda}$ in these patterns, and if we do not recognise and
incorporate them into our notation, then expressions rapidly become
unreadable, and then intractable. Below, $q$-brackets will appear
at every point, and even in our R matrices.

\pagebreak

As mentioned in \S\ref{sec:NormalOrdering}, the normal ordering of
generator strings is a computationally expensive task. Here, the normal
ordering of $Z_{2^{mn}}$ typically dominates the computations as it
demands that we process a seriously disordered string of length $2mn$.
A theoretical insight would be valuable here -- for example an explicit
formula for the normal ordering of arbitrary ${\mathbf{E}^b}_a {E^a}_b$
would help speed the evaluation of the $\beta_i$.  The regularities
apparent in the above example suggest such the existence of such a
result, and (\ref{eq:Zhang92Lemma3}) may also be of use.
Alternatively, a more efficient computation of the $\beta_i$ should be
possible by the efficient reuse of previous calculations.

\pagebreak


\subsection{Matrix elements for $\pi \equiv \pi_{\Lambda}$}
\label{sec:Matrixelements}

Having established an orthonormal basis $B$ for the module
$V\equiv V_{\Lambda}$ corresponding to the representation
$\pi\equiv\pi_{\Lambda}$, we now use it to construct matrix elements
$\pi(X)$ for $U_q[gl(m|n)]$ generators $X$.

Where $I_{2^{mn}}$ is the identity transformation on $V$,
in the basis $B$ we have the identity:
$\sum_{i=1}^{2^{mn}}\ketbra{i}{i}=I_{2^{mn}}$, thus:
\begin{eqnarray}
  \pi (X)
  =
  \pi (X)
  \cdot
  \sum_{i=1}^{2^{mn}}
    \ketbra{i}{i}
  =
  \sum_{i=1}^{2^{mn}}
    (\pi (X) \cdot \ket{i}) \cdot \bra{i}
  =
  \sum_{i=1}^{2^{mn}}
    (X \cdot \ket{i}) \cdot \bra{i},
  \label{eq:howtobuildmatrixelements}
\end{eqnarray}
and the action $X\cdot\ket{i}$ may be computed knowing the expansion
of $\ket{i}$ in terms of generator products and $\ket{1}$. Thus,
if, as in (\ref{eq:ketdefinition}), we have:
\begin{eqnarray*}
  \ket{i}
  =
  \beta_{i}
  {E^{a_1}}_{b_1}
  {E^{a_2}}_{b_2}
  \cdots
  {E^{a_p}}_{b_p}
  \cdot
  \ket{1},
\end{eqnarray*}
then we may compute $X \cdot \ket{i}$ by the following process:

\begin{enumerate}
\item
  We again use the PBW basis commutators to normally order the string
  $
    Y
    \triangleq
    X
    {E^{a_1}}_{b_1}
    {E^{a_2}}_{b_2}
    \cdots
    {E^{a_p}}_{b_p}
  $,
  denoting the result by $\mathrm{NO}(Y)$.

\item
  We use the known actions of the raising and Cartan generators
  on $\ket{1}$ to reduce $\mathrm{NO}(Y) \cdot \ket{1}$ to an
  expression which is generally a \emph{sum} of scalar-multiplied,
  normally ordered products of (odd) lowering generators acting on
  $\ket{1}$.

\item
  Identifying the terms in the resulting products as scalar
  multiples of various $\ket{j}$, we obtain the result. To whit, if
  $\mathrm{NO}(Y)$ contains a term of the form:
  \begin{eqnarray*}
    T \cdot \ket{1}
    =
    \theta
    {E^{c_1}}_{d_1}
    {E^{c_2}}_{d_2}
    \cdots
    {E^{c_r}}_{d_r}
    \cdot
    \ket{1},
  \end{eqnarray*}
  for some scalar $\theta$, and $r$ odd lowering generators
  ${E^{c_j}}_{d_j}$, and we know that:
  \begin{eqnarray*}
    \ket{k}
    =
    \beta_{k}
    {E^{c_1}}_{d_1}
    {E^{c_2}}_{d_2}
    \cdots
    {E^{c_r}}_{d_r}
    \cdot
    \ket{1},
  \end{eqnarray*}
  then we may replace
  $T \cdot \ket{1}$
  with
  $\theta\overline{\beta}_{k}\ket{k}$.
\end{enumerate}

Repeating this procedure for all $2^{mn}$ basis vectors $\ket{i}$, and
substituting the results into (\ref{eq:howtobuildmatrixelements})
yields the required matrix element $\pi(X)$.  Again, as the ordering
chosen for the basis vectors $\ket{i}$ is compatible with the ordering
used in the PBW lemma, this process is robust.

We now divide the construction of matrix elements into two phases.
Firstly, in a direct implementation of the above, we build matrix
elements for the Cartan and simple raising generators.  We illustrate
this for the $U_q[gl(3|1)]$ case, for the generators $K_4$ (in
\S\ref{sec:MatrixelementsK4}) and ${E^{3}}_{4}$ (in
\S\ref{sec:MatrixelementsE34}).

Having done that, in \S\ref{sec:Theothermatrixelements}, we describe
the construction of the remaining matrix elements, as they may be
efficiently computed in terms of those for the simple lowering
generators.

Further illustrations are provided in my PhD thesis \cite{DeWit:98},
using the $U_q[gl(2|1)]$ case, although those results are somewhat less
formally explained.


\subsubsection{Matrix elements $\pi_{\Lambda}(K_4)$
               for the $U_q[gl(3|1)]$ case}
\label{sec:MatrixelementsK4}

Firstly, we must normal order a list of $8$ generator strings,
cf. (\ref{eq:3ketoriginlist}). We obtain:
\small
\begin{eqnarray}
  K_{4}
  \cdot
  \left(
    \begin{array}{l}
      \mathrm{Id}
      \\
      {E^{4}}_{3}
      \\
      {E^{4}}_{2}
      \\
      {E^{4}}_{1}
      \\
      {E^{4}}_{3} {E^{4}}_{2}
      \\
      {E^{4}}_{3} {E^{4}}_{1}
      \\
      {E^{4}}_{2} {E^{4}}_{1}
      \\
      {E^{4}}_{3} {E^{4}}_{2} {E^{4}}_{1}
    \end{array}
  \right)
  \stackrel{(\ref{eq:FullCNCCommutator})}{=}
  \left(
    \begin{array}{r}
      \mathrm{Id} ~~
      \\
      \overline{q}
      {E^{4}}_{3}
      \\
      \overline{q}
      {E^{4}}_{2}
      \\
      \overline{q}
      {E^{4}}_{1}
      \\
      \overline{q}^{2}
      {E^{4}}_{3} {E^{4}}_{2}
      \\
      \overline{q}^{2}
      {E^{4}}_{3} {E^{4}}_{1}
      \\
      \overline{q}^{2}
      {E^{4}}_{2} {E^{4}}_{1}
      \\
      \overline{q}^{3}
      {E^{4}}_{3} {E^{4}}_{2} {E^{4}}_{1}
    \end{array}
  \right)
  \cdot
  K_{4}.
  \label{eq:Kitami}
\end{eqnarray}
\normalsize
The action $K_{4}\cdot\ket{1}$ is known
explicitly from (\ref{eq:Cartanactiononket1}), that is, we have
$K_{4}\cdot\ket{1}=\overline{q}^{\alpha}\ket{1}$.
Thus, we have:
\small
\begin{eqnarray*}
  K_{4}
  \cdot
  \left(
    \begin{array}{r}
      \ket{1}
      \\
      \ket{2}
      \\
      \ket{3}
      \\
      \ket{4}
      \\
      \ket{5}
      \\
      \ket{6}
      \\
      \ket{7}
      \\
      \ket{8}
    \end{array}
  \right)
  & \stackrel{(\ref{eq:Kitami})}{=} &
  \overline{q}^{\alpha}
  \left(
    \begin{array}{r}
      1 ~~~~~
      \\
      \overline{q}
      \beta_{2}
      {E^{4}}_{3}
      \\
      \overline{q}
      \beta_{3}
      {E^{4}}_{2}
      \\
      \overline{q}
      \beta_{4}
      {E^{4}}_{1}
      \\
      \overline{q}^{2}
      \beta_{5}
      {E^{4}}_{3} {E^{4}}_{2}
      \\
      \overline{q}^{2}
      \beta_{6}
      {E^{4}}_{3} {E^{4}}_{1}
      \\
      \overline{q}^{2}
      \beta_{7}
      {E^{4}}_{2} {E^{4}}_{1}
      \\
      \overline{q}^{3}
      \beta_{8}
      {E^{4}}_{3} {E^{4}}_{2} {E^{4}}_{1}
    \end{array}
  \right)
  \cdot
  \ket{1}
  \stackrel{(\ref{eq:3ketoriginlist})}{=}
  \overline{q}^{\alpha}
  \left(
    \begin{array}{l@{\hspace{2pt}}r}
      & \ket{1}
      \\
      \overline{q}
      & \ket{2}
      \\
      \overline{q}
      & \ket{3}
      \\
      \overline{q}
      & \ket{4}
      \\
      \overline{q}^{2}
      & \ket{5}
      \\
      \overline{q}^{2}
      & \ket{6}
      \\
      \overline{q}^{2}
      & \ket{7}
      \\
      \overline{q}^{3}
      & \ket{8}
    \end{array}
  \right).
\end{eqnarray*}
\normalsize

Installing this information into (\ref{eq:howtobuildmatrixelements}),
we discover:
\begin{eqnarray*}
  \pi (K_{4})
  & = &
  \overline{q}^{\alpha}
  \left(
  \begin{array}{l}
    \ketbra{1}{1}
    +
    \overline{q}
    \left(
      \ketbra{2}{2}
      +
      \ketbra{3}{3}
      +
      \ketbra{4}{4}
    \right)
    \\
    \hspace{23pt}
    +
    \overline{q}^{2}
    \left(
      \ketbra{5}{5}
      +
      \ketbra{6}{6}
      +
      \ketbra{7}{7}
    \right)
    +
    \overline{q}^{3}
    \ketbra{8}{8}
  \end{array}
  \right).
\end{eqnarray*}

Note that this process \emph{doesn't} actually require explicit
knowledge of the $\beta_{i}$ presented in (\ref{eq:3betalisting}). This
is a feature of the evaluation of matrix elements of Cartan generators;
the more usual situation appears in the next example (i.e.
\S\ref{sec:MatrixelementsE34}).  Replacing $\ketbra{i}{j}$ with the
elementary matrix $e^i_j\in M_{2^{mn}}$, we have:
\begin{eqnarray*}
  \pi (K_{4})
  =
  \overline{q}^{\alpha}
  \left(
    e^1_1
    +
    \overline{q}
    \left(
      e^2_2
      +
      e^3_3
      +
      e^4_4
    \right)
    +
    \overline{q}^{2}
    \left(
      e^5_5
      +
      e^6_6
      +
      e^7_7
    \right)
    +
    \overline{q}^{3}
    e^8_8
  \right),
\end{eqnarray*}
that is:
\small
\begin{eqnarray*}
  \pi (K_{4})
  =
  \overline{q}^{\alpha}
  \left(
    \begin{array}{*8{c}}
      1 & 0 & 0 & 0 & 0 & 0 & 0 & 0 \\
      0 & \overline{q} & 0 & 0 & 0 & 0 & 0 & 0 \\
      0 & 0 & \overline{q} & 0 & 0 & 0 & 0 & 0 \\
      0 & 0 & 0 & \overline{q} & 0 & 0 & 0 & 0 \\
      0 & 0 & 0 & 0 & \overline{q}^{2} & 0 & 0 & 0 \\
      0 & 0 & 0 & 0 & 0 & \overline{q}^{2} & 0 & 0 \\
      0 & 0 & 0 & 0 & 0 & 0 & \overline{q}^{2} & 0 \\
      0 & 0 & 0 & 0 & 0 & 0 & 0 & \overline{q}^{3}
    \end{array}
  \right).
\end{eqnarray*}
\normalsize


\subsubsection{Matrix elements $\pi_{\Lambda}({E^{3}}_{4})$
                 for the $U_q[gl(3|1)]$ case}
\label{sec:MatrixelementsE34}

This is a more interesting case than that of $K_{4}$.  This time,
the normal ordering of the list of $8$ generator strings yields:
\small
\begin{eqnarray*}
  \hspace{-25pt}
  {E^{3}}_{4}
  \hspace{-2pt}
  \cdot
  \hspace{-2pt}
  \left(
    \begin{array}{@{\hspace{0pt}}l@{\hspace{0pt}}}
      \mathrm{Id}
      \\
      {E^{4}}_{3}
      \\
      {E^{4}}_{2}
      \\
      {E^{4}}_{1}
      \\
      {E^{4}}_{3} {E^{4}}_{2}
      \\
      {E^{4}}_{3} {E^{4}}_{1}
      \\
      {E^{4}}_{2} {E^{4}}_{1}
      \\
      {E^{4}}_{3} {E^{4}}_{2} {E^{4}}_{1}
    \end{array}
  \right)
  \hspace{-2pt}
  =
  \hspace{-2pt}
  \left(
    \begin{array}{@{\hspace{0pt}}l@{\hspace{0pt}}}
      {E^{3}}_{4}
      \\
      -
      {E^{4}}_{3} {E^{3}}_{4}
      +
      \overline{\Delta}
      (
        \overline{K}_{4} K_{3}
        -
        K_{4} \overline{K}_{3}
      )
      \\
      -
      {E^{4}}_{2} {E^{3}}_{4}
      +
      {E^{3}}_{2} K_{4} \overline{K}_{3}
      \\
      -
      {E^{4}}_{1} {E^{3}}_{4}
      +
      {E^{3}}_{1} K_{4} \overline{K}_{3}
      \\
      \overline{\Delta}
      {E^{4}}_{2}
      (
        q \overline{K}_{4} K_{3}
        -
        \overline{q} K_{4} \overline{K}_{3}
      )
      +
      {E^{4}}_{3}
      (
        {E^{4}}_{2} {E^{3}}_{4}
        -
        {E^{3}}_{2} K_{4} \overline{K}_{3}
      )
      \\
      \overline{\Delta}
      {E^{4}}_{1}
      (
        q \overline{K}_{4} K_{3}
        -
        \overline{q} K_{4} \overline{K}_{3}
      )
      +
      {E^{4}}_{3}
      (
        {E^{4}}_{1} {E^{3}}_{4}
        -
        {E^{3}}_{1} K_{4} \overline{K}_{3}
      )
      \\
      {E^{4}}_{2} {E^{4}}_{1} {E^{3}}_{4}
      +
      \overline{q}
      {E^{4}}_{1} {E^{3}}_{2} K_{4} \overline{K}_{3}
      -
      {E^{4}}_{2} {E^{3}}_{1} K_{4} \overline{K}_{3}
      \\
      -
      {E^{4}}_{3}
      (
        {E^{4}}_{2} {E^{4}}_{1} {E^{3}}_{4}
        +
        \overline{q}
        {E^{4}}_{1} {E^{3}}_{2} K_{4} \overline{K}_{3}
        -
        {E^{4}}_{2} {E^{3}}_{1} K_{4} \overline{K}_{3}
      )
      \\
      \hspace{76pt}
      +
      \overline{\Delta}
      {E^{4}}_{2} {E^{4}}_{1}
      (
        q^2
        \overline{K}_{4} K_{3}
        -
        \overline{q}^2
        K_{4} \overline{K}_{3}
      )
    \end{array}
  \hspace{-2pt}
  \right)
  \hspace{-2pt}
  ,
\end{eqnarray*}
\normalsize
where we have again written $\Delta=q-\overline{q}$ and
$\overline{\Delta}=\Delta^{-1}$ (cf. (\ref{eq:DefnofDeltaa})), and done
a little judicious factoring to improve readability.  Next, we know
that any terms which end with raising generators (e.g.  ${E^{3}}_{4}$)
will annihilate $\ket{1}$. Further, any term which contains an even
lowering generator immediately to the left of a terminal string of
Cartan generators%
\footnote{%
  For example, the last product contains
  ${E^{4}}_{3} {E^{4}}_{2} {E^{3}}_{1} K_{4} \overline{K}_{3}$, which
  contains the even lowering generator ${E^{3}}_{1}$ immediately to the
  left of a terminal string of Cartan generators.
}
will annihilate $\ket{1}$, as the action of the Cartan generators is
purely scalar. Omitting all such terms, we may thus write:
\small
\begin{eqnarray*}
  \hspace{0pt}
  {E^{3}}_{4}
  \cdot
  \left(
    \begin{array}{r}
      \mathrm{Id}
      \\
      {E^{4}}_{3}
      \\
      {E^{4}}_{2}
      \\
      {E^{4}}_{1}
      \\
      {E^{4}}_{3} {E^{4}}_{2}
      \\
      {E^{4}}_{3} {E^{4}}_{1}
      \\
      {E^{4}}_{2} {E^{4}}_{1}
      \\
      {E^{4}}_{3} {E^{4}}_{2} {E^{4}}_{1}
    \end{array}
  \right)
  \cdot
  \ket{1}
  =
  \left(
    \begin{array}{r}
      0
      \\
      \overline{\Delta}
      (
        \overline{K}_{4} K_{3}
        -
        K_{4} \overline{K}_{3}
      )
      \\
      0
      \\
      0
      \\
      \overline{\Delta}
      {E^{4}}_{2}
      (
        q \overline{K}_{4} K_{3}
        -
        \overline{q} K_{4} \overline{K}_{3}
      )
      \\
      \overline{\Delta}
      {E^{4}}_{1}
      (
        q \overline{K}_{4} K_{3}
        -
        \overline{q} K_{4} \overline{K}_{3}
      )
      \\
      0
      \\
      \overline{\Delta}
      {E^{4}}_{2} {E^{4}}_{1}
      (
        q^2 \overline{K}_{4} K_{3}
        -
        \overline{q}^2 K_{4} \overline{K}_{3}
      )
    \end{array}
  \right)
  \cdot
  \ket{1}
\end{eqnarray*}
\normalsize

Again, the action $K_{a}^N\cdot\ket{1}$ is known explicitly from
(\ref{eq:Cartanactiononket1}), that is, we have
$K_{3}^{\pm}\cdot\ket{1}=\ket{1}$ and, as above, that
$K_{4}^{\pm}\cdot\ket{1}=q^{\mp\alpha}\ket{1}$.  Thus, we have:
\small
\begin{eqnarray}
  {E^{3}}_{4}
  \cdot
  \left(
    \begin{array}{r}
      \mathrm{Id}
      \\
      {E^{4}}_{3}
      \\
      {E^{4}}_{2}
      \\
      {E^{4}}_{1}
      \\
      {E^{4}}_{3} {E^{4}}_{2}
      \\
      {E^{4}}_{3} {E^{4}}_{1}
      \\
      {E^{4}}_{2} {E^{4}}_{1}
      \\
      {E^{4}}_{3} {E^{4}}_{2} {E^{4}}_{1}
    \end{array}
  \right)
  \cdot
  \ket{1}
  =
  \left(
    \begin{array}{r}
      0
      \\
      {[\alpha]}_{q}
      \\
      0
      \\
      0
      \\
      {[\alpha+1]}_{q}
      {E^{4}}_{2}
      \\
      {[\alpha+1]}_{q}
      {E^{4}}_{1}
      \\
      0
      \\
      {[\alpha+2]}_{q}
      {E^{4}}_{2} {E^{4}}_{1}
    \end{array}
  \right)
  \cdot
  \ket{1}.
  \label{eq:Abashiri}
\end{eqnarray}
\normalsize
This time explicitly installing information of the $\beta_{i}$ from
(\ref{eq:3betalisting}), we have:
\small
\begin{eqnarray}
  \hspace{-25pt}
  {E^{3}}_{4}
  \cdot
  \left(
    \begin{array}{@{\hspace{0pt}}r@{\hspace{0pt}}}
      \ket{1}
      \\
      \ket{2}
      \\
      \ket{3}
      \\
      \ket{4}
      \\
      \ket{5}
      \\
      \ket{6}
      \\
      \ket{7}
      \\
      \ket{8}
    \end{array}
  \right)
  & \stackrel{(\ref{eq:3ketoriginlist})}{=} &
  {E^{3}}_{4}
  \cdot
  \left(
    \begin{array}{@{\hspace{0pt}}r@{\hspace{0pt}}}
      \beta_{1}
      \\
      \beta_{2} {E^{4}}_{3}
      \\
      \beta_{3} {E^{4}}_{2}
      \\
      \beta_{4} {E^{4}}_{1}
      \\
      \beta_{5} {E^{4}}_{3} {E^{4}}_{2}
      \\
      \beta_{6} {E^{4}}_{3} {E^{4}}_{1}
      \\
      \beta_{7} {E^{4}}_{2} {E^{4}}_{1}
      \\
      \beta_{8} {E^{4}}_{3} {E^{4}}_{2} {E^{4}}_{1}
    \end{array}
  \right)
  \cdot
  \ket{1}
  \stackrel{(\ref{eq:Abashiri})}{=}
  \left(
    \begin{array}{@{\hspace{0pt}}r@{\hspace{0pt}}}
      0
      \\
      \beta_{2}
      {[\alpha]}_{q}
      \ket{1}
      \\
      0
      \\
      0
      \\
      \beta_{5}
      {[\alpha+1]}_{q}
      {E^{4}}_{2}
      \cdot
      \ket{1}
      \\
      \beta_{6}
      {[\alpha+1]}_{q}
      {E^{4}}_{1}
      \cdot
      \ket{1}
      \\
      0
      \\
      \beta_{8}
      {[\alpha+2]}_{q}
      {E^{4}}_{2} {E^{4}}_{1}
      \cdot
      \ket{1}
    \end{array}
  \right)
  \nonumber
  \\
  & \stackrel{(\ref{eq:3ketoriginlist})}{=} &
  \left(
    \begin{array}{@{\hspace{0pt}}r@{\hspace{0pt}}}
      0
      \\
      \beta_{2}
      {[\alpha]}_{q}
      \ket{1}
      \\
      0
      \\
      0
      \\
      \beta_{5}
      \overline{\beta}_{3}
      {[\alpha+1]}_{q}
      \ket{3}
      \\
      \beta_{6}
      \overline{\beta}_{4}
      {[\alpha+1]}_{q}
      \ket{4}
      \\
      0
      \\
      \beta_{8}
      \overline{\beta}_{7}
      {[\alpha+2]}_{q}
      \ket{7}
    \end{array}
  \right)
  \stackrel{(\ref{eq:3betalisting})}{=}
  \left(
    \begin{array}{@{\hspace{0pt}}r@{\hspace{0pt}}}
      0
      \\
      {[\alpha]}_{q}^{\frac{1}{2}}
      \ket{1}
      \\
      0
      \\
      0
      \\
      {[\alpha+1]}_q^{\frac{1}{2}}
      \ket{3}
      \\
      {[\alpha+1]}_q^{\frac{1}{2}}
      \ket{4}
      \\
      0
      \\
      {[\alpha+2]}_q^{\frac{1}{2}}
      \ket{7}
    \end{array}
  \right).
  \label{eq:Barrow}
\end{eqnarray}
\normalsize

Note that there is a subtle point in the second application of
(\ref{eq:3ketoriginlist}), in that we are using it \emph{implicitly}:
this point wasn't so clear when we computed $\pi(K_4)$.  Whilst this is
simple enough for a human to perform, computer programs require
\emph{explicit} algorithms. In practice, what this means is that we
must invert (\ref{eq:3ketoriginlist}) to provide a list of
transformation rules for strings of odd lowering generators acting on
$\ket{1}$ that allow us to recover the $\ket{i}$. These rules must be
carefully coded, and applied in \emph{reverse} order to
(\ref{eq:3ketoriginlist}), to ensure that we collect the longest
strings first. Again, the ordering of the $\ket{i}$ facilitates this.
To illustrate, we want to apply, \emph{in order}, the following rules:
\begin{eqnarray*}
  \left\{
    \begin{array}{rcl}
      {E^{4}}_{3} {E^{4}}_{2} {E^{4}}_{1} \cdot \ket{1}
      & \mapsto &
      \overline{\beta}_{8}
      \ket{8}
      \\
      {E^{4}}_{2} {E^{4}}_{1} \cdot \ket{1}
      & \mapsto &
      \overline{\beta}_{7}
      \ket{7}
      \\
      {E^{4}}_{3} {E^{4}}_{1} \cdot \ket{1}
      & \mapsto &
      \overline{\beta}_{6}
      \ket{6}
      \\
      {E^{4}}_{3} {E^{4}}_{2} \cdot \ket{1}
      & \mapsto &
      \overline{\beta}_{5}
      \ket{5}
      \\
      {E^{4}}_{1} \cdot \ket{1}
      & \mapsto &
      \overline{\beta}_{4}
      \ket{4}
      \\
      {E^{4}}_{2} \cdot \ket{1}
      & \mapsto &
      \overline{\beta}_{3}
      \ket{3}
      \\
      {E^{4}}_{3} \cdot \ket{1}
      & \mapsto &
      \overline{\beta}_{2}
      \ket{2}
    \end{array}
  \right\}.
\end{eqnarray*}

Installing the information in (\ref{eq:Barrow})
into (\ref{eq:howtobuildmatrixelements}), we discover, again
substituting $e^i_j$ for $\ketbra{i}{j}$:
\begin{eqnarray*}
  \hspace{-11pt}
  \pi ({E^{3}}_{4})
  & = &
  {[\alpha]}_{q}^{\frac{1}{2}}
  \ketbra{1}{2}
  +
  {[\alpha+1]}_q^{\frac{1}{2}}
  \ketbra{3}{5}
  +
  {[\alpha+1]}_q^{\frac{1}{2}}
  \ketbra{4}{6}
  +
  {[\alpha+2]}_q^{\frac{1}{2}}
  \ketbra{7}{8}
  \\
  & = &
  {[\alpha]}_{q}^{\frac{1}{2}}
  e^1_2
  +
  {[\alpha+1]}_q^{\frac{1}{2}}
  e^3_5
  +
  {[\alpha+1]}_q^{\frac{1}{2}}
  e^4_6
  +
  {[\alpha+2]}_q^{\frac{1}{2}}
  e^7_8,
\end{eqnarray*}
and interpreting $e^i_j$ as an elementary matrix:
\small
\begin{eqnarray*}
  \pi ({E^{3}}_{4})
  =
  \left(
    \begin{array}{*{8}{c}}
      0 &
          {[\alpha]}_{q}^{\frac{1}{2}}
            & 0 & 0 & 0 & 0 & 0 & 0 \\
      0 & 0 & 0 & 0 & 0 & 0 & 0 & 0 \\
      0 & 0 & 0 & 0 &
                      {[\alpha+1]}_q^{\frac{1}{2}}
                        & 0 & 0 & 0 \\
      0 & 0 & 0 & 0 & 0 &
                          {[\alpha+1]}_q^{\frac{1}{2}}
                            & 0 & 0 \\
      0 & 0 & 0 & 0 & 0 & 0 & 0 & 0 \\
      0 & 0 & 0 & 0 & 0 & 0 & 0 & 0 \\
      0 & 0 & 0 & 0 & 0 & 0 & 0 &
                                  {[\alpha+2]}_q^{\frac{1}{2}} \\
      0 & 0 & 0 & 0 & 0 & 0 & 0 & 0
    \end{array}
  \right).
\end{eqnarray*}
\normalsize


\subsubsection{Matrix elements for the remaining generators}
\label{sec:Theothermatrixelements}

Having found the matrix elements for the \emph{simple raising}
generators $\pi({E^{a}}_{b})$, we may immediately write down the matrix
elements for the corresponding \emph{simple lowering} generators
$\pi({E^{b}}_{a})$, as these are necessarily the transposes of
$\pi({E^{a}}_{b})$.

To illustrate, for the $U_q[gl(3|1)]$ case, we have:
\begin{eqnarray*}
  \pi ({E^{4}}_{3})
  =
  {[\alpha]}_{q}^{\frac{1}{2}}
  e^2_1
  +
  {[\alpha+1]}_q^{\frac{1}{2}}
  e^5_3
  +
  {[\alpha+1]}_q^{\frac{1}{2}}
  e^6_4
  +
  {[\alpha+2]}_q^{\frac{1}{2}}
  e^8_7.
\end{eqnarray*}
Beyond this, we may construct matrix elements for the \emph{nonsimple}
generators from those for the simple ones by (recursively) applying
$\pi$ to (\ref{eq:UqglmnNonSimpleGeneratorsBetterNotation}a), viz:
\begin{eqnarray*}
  \pi({E^a}_b)
  =
  \pi({E^a}_c)
  \pi({E^c}_b)
  -
  q_c^{S^a_b}
  \pi({E^c}_b)
  \pi({E^a}_c),
\end{eqnarray*}
although this may not yield particularly \emph{useful} results. Note
that we can in fact write down the matrix elements for any (simple and
nonsimple) \emph{lowering} generators directly from their corresponding
raising generators, as we have, for $a>b$:
$\pi({E^a}_b)=\pi({E^b}_a)^{T_q}$, where the `$q$ transpose' $T_q$
indicates the combination of the transpose and the mapping
$q\mapsto\overline{q}$. This follows from inspection of the way that
(\ref{eq:UqglmnNonSimpleGeneratorsBetterNotation}a) depends on $q$;
indeed it is trivial for simple generators. Illustrations are visible
below (and in Appendix \ref{app:MatrixElements}); the reader
should note that the $q$ bracket is invariant under
$q\mapsto\overline{q}$.

For completeness, we list the matrix elements of the $U_q[gl(3|1)]$
\emph{simple} generators (the matrix elements of \emph{all}
the generators are listed in Appendix \ref{app:MatrixElements}):
\begin{eqnarray}
  \hspace{-30pt}
  \left.
  \begin{array}{l@{\hspace{3pt}}r@{\hspace{4pt}}c@{\hspace{4pt}}l}
    (\mathrm{a}) &
    \pi ( K_{1} )
    & = &
                    e^1_1
    +               e^2_2
    +               e^3_3
    + \overline{q}  e^4_4
    +               e^5_5
    + \overline{q}  e^6_6
    + \overline{q}  e^7_7
    + \overline{q}  e^8_8
    \\
    (\mathrm{b}) &
    \pi ( K_{2} )
    & = &
                    e^1_1
    +               e^2_2
    + \overline{q}  e^3_3
    +               e^4_4
    + \overline{q}  e^5_5
    +               e^6_6
    + \overline{q}  e^7_7
    + \overline{q}  e^8_8
    \\
    (\mathrm{c}) &
    \pi ( K_{3} )
    & = &
                    e^1_1
    + \overline{q}  e^2_2
    +               e^3_3
    +               e^4_4
    + \overline{q}  e^5_5
    + \overline{q}  e^6_6
    +               e^7_7
    + \overline{q}  e^8_8
    \\
    (\mathrm{d}) &
    \pi ( K_4 )
    & = &
      \overline{q}^{\alpha}      e^1_1
    + \overline{q}^{\alpha + 1}  (e^2_2 + e^3_3 + e^4_4)
    + \overline{q}^{\alpha + 2}  (e^5_5 + e^6_6 + e^7_7)
    + \overline{q}^{\alpha + 3}  e^8_8
    \\
    (\mathrm{e}) &
    \pi ( {E^1}_2 )
    & = &
    - e^3_4
    - e^5_6
    \\
    (\mathrm{f}) &
    \pi ( {E^2}_1 )
    & = &
    - e^4_3
    - e^6_5
    \\
    (\mathrm{g}) &
    \pi ( {E^2}_3 )
    & = &
    - e^2_3
    - e^6_7
    \\
    (\mathrm{h}) &
    \pi ( {E^3}_2 )
    & = &
    - e^3_2
    - e^7_6
    \\
    (\mathrm{i}) &
    \pi ({E^{3}}_{4})
    & = &
    {[\alpha]}_{q}^{\frac{1}{2}}
    e^1_2
    +
    {[\alpha+1]}_q^{\frac{1}{2}}
    e^3_5
    +
    {[\alpha+1]}_q^{\frac{1}{2}}
    e^4_6
    +
    {[\alpha+2]}_q^{\frac{1}{2}}
    e^7_8
    \\
    (\mathrm{j}) &
    \pi ({E^{4}}_{3})
    & = &
    {[\alpha]}_{q}^{\frac{1}{2}}
    e^2_1
    +
    {[\alpha+1]}_q^{\frac{1}{2}}
    e^5_3
    +
    {[\alpha+1]}_q^{\frac{1}{2}}
    e^6_4
    +
    {[\alpha+2]}_q^{\frac{1}{2}}
    e^8_7.
  \end{array}
  \right\}
  \label{eq:Uqgl31matrixelements}
\end{eqnarray}

\pagebreak


\section{The submodules $V_k \subset V\otimes V$}
\label{sec:ThesubmodulesVk}

We now turn our attention to the tensor product module $V\otimes V$.
Where $V$ has a basis $B=\{\ket{i}\}_{i=1}^{2^{mn}}$,
the $2^{2mn}$ dimensional $V\otimes V$ has a natural
basis $\{\ket{i}\otimes\ket{j}\}_{i,j=1}^{2^{mn}}$,
which inherits a weight system and a grading from $V$:
\begin{eqnarray}
  \mathrm{wt}(\ket{i}\otimes\ket{j})
  & \triangleq &
  \mathrm{wt}(\ket{i})
  +
  \mathrm{wt}(\ket{j})
  \label{eq:wtsofketket}
  \\
  \left[\ket{i}\otimes\ket{j}\right]
  & \triangleq &
  \left[\ket{i}\right]
  +
  \left[\ket{j}\right]
  \quad
  (\mathrm{mod}\;2).
  \nonumber
\end{eqnarray}

To build an R matrix acting on $V\otimes V$, we will use an alternative,
orthonormal weight basis $\mathfrak{B}$ for $V \otimes V$, which
corresponds to the (known) decomposition of $V\otimes V$ into
irreducible $U_q[gl(m|n)]$ submodules.  The basis vectors of
$\mathfrak{B}$ are expressed as linear combinations
$\theta_{ij} (\ket{i}\otimes \ket{j})$, where the coefficients
$\theta_{ij}$ are in general algebraic expressions in $q$ and $\alpha$.
Before proceeding, we introduce some machinery for dealing with the
tensor products.


\subsection{Tensor product representation tools}

We first introduce the action of $U_q[gl(m|n)]\otimes U_q[gl(m|n)]$ on
$V\otimes V$. Where $Y$ is an homogeneous $U_q[gl(m|n)]$ element, we
define:
\begin{eqnarray}
  (X \otimes Y)
  \cdot
  ( \ket{i} \otimes \ket{j} )
  \triangleq
  (-)^{[Y][i]}
  (X \cdot \ket{i} \otimes Y \cdot \ket{j}),
  \label{eq:XYxketketdefn}
\end{eqnarray}
and extend by linearity to all of $U_q[gl(m|n)]\otimes U_q[gl(m|n)]$;
note that we have written $[i]\equiv [\ket{i}]$ for readability.  This
action is compatible with the grading.

Next, we define a basis
$\{\bra{i}\otimes\bra{j}\}_{i,j=1}^{2^{mn}}$ for
$(V\otimes V)^*$ (the dual to $V\otimes V$) by dualising each of the
elements of the basis $\{\ket{i}\otimes\ket{j}\}_{i,j=1}^{2^{mn}}$:
\begin{eqnarray}
  \bra{i}\otimes\bra{j}
  \triangleq
  (-)^{[i][j]}
  {\left( \ket{i}\otimes\ket{j} \right)}^\dagger,
  \label{eq:defnofdualonVoV}
\end{eqnarray}
where we intend
$
  {\left( \ket{i}\otimes\ket{j} \right)}^\dagger
  =
  \ket{i}^\dagger \otimes\ket{j}^\dagger
$.
The conjugate is extended by linearity:
\begin{eqnarray*}
  {\left(
    A \ket{i}\otimes\ket{j}
    +
    B \ket{k}\otimes\ket{l}
  \right)}^\dagger
  \triangleq
  A^* {\left(\ket{i}\otimes\ket{j}\right)}^\dagger
  +
  B^* {\left(\ket{k}\otimes\ket{l}\right)}^\dagger,
\end{eqnarray*}
for scalars $A$ and $B$, where $A^*$ is the complex conjugate of $A$.

The multiplication operations between the dual bases are given by:
\begin{eqnarray}
  (\ket{i}\otimes\ket{j}) \cdot (\bra{k}\otimes\bra{l})
  & \triangleq &
  (-)^{[j][k]}
  (\ketbra{i}{k} \otimes \ketbra{j}{l})
  \label{eq:brabraxketketdefn}
  \\
  (\bra{i}\otimes\bra{j}) \cdot (\ket{k}\otimes\ket{l})
  & \triangleq &
  (-)^{[j][k]}
  \delta^i_k \delta^j_l,
  \label{eq:ketketxbrabradefn}
\end{eqnarray}
analagously to (\ref{eq:XYxketketdefn}).

Thus, we have the natural inner product on the basis
$\{\ket{i}\otimes\ket{j}\}_{i,j=1}^{2^{mn}}$:
\begin{eqnarray}
  \left( \ket{i}\otimes\ket{j}, \ket{k}\otimes\ket{l} \right)
  \triangleq
  (\ket{i}\otimes\ket{j})^\dagger \cdot (\ket{k}\otimes\ket{l}),
  \label{eq:innerproductonVoV}
\end{eqnarray}
which behaves as expected:
\begin{eqnarray*}
  \left( \ket{i}\otimes\ket{j}, \ket{k}\otimes\ket{l} \right)
  & \stackrel{(\ref{eq:defnofdualonVoV})}{=} &
  (-)^{[i][j]}
  (\bra{i}\otimes\bra{j}) \cdot (\ket{k}\otimes\ket{l})
  \\
  & = &
  (-)^{[i][j]}
  (-)^{[j][k]}
  (\braket{i}{k} \otimes \braket{j}{l})
  \\
  & = &
  (-)^{[j]([i]+[k])}
  \delta^i_k \delta^j_l
  =
  \delta^i_k \delta^j_l.
\end{eqnarray*}

Lastly, we will often use the shorthand
$\ketket{i}{j}\triangleq\ket{i}\otimes\ket{j}$ and
$\brabra{i}{j}\triangleq\bra{i}\otimes\bra{j}$.


\subsection{The orthogonal decomposition of $V\otimes V$}
\label{sec:orthogonaldecompositionofVotimesV}

For our modules $V\equiv V_{\Lambda}$, the orthogonal decomposition of
$V\otimes V$ is known, and contains no multiplicities
\cite{DeliusGouldLinksZhang:95b,GouldLinksZhang:96b}. To describe it,
we introduce a little notation.  What follows is strictly true only for
$m\geqslant n$, but the natural isomorphism between $U_q[gl(m|n)]$ and
$U_q[gl(n|m)]$ shows that we need not consider the other case
$n\geqslant m$.

For any $0\leqslant N\leqslant mn$, let $\gamma$
be a nonincreasing sequence of nonnegative integers%
\footnote{%
  We apologise for overloading $\gamma$ -- this one is a sequence,
  not a root.
}
$\gamma_1, \gamma_2, \dots, \gamma_r$, where the $\gamma_k$ satisfy
$\sum_{k=1}^{r}\gamma_k=N$.  We then define a \emph{Young diagram}
$[\gamma]=[\gamma_1,\gamma_2,\dots,\gamma_r]$, to be \emph{allowable}
if it has at most $n$ columns and $m$ rows, viz $r\leqslant m$ and, for
each $k$, $\gamma_k \leqslant n$.  To each such allowable diagram, we
associate a $gl(m|n)$ weight $\lambda_{\gamma}$:
\begin{eqnarray}
  \hspace{-13pt}
  \lambda_{\gamma}
  =
  (
    \dot{0}_{m-r},
    - \gamma_{r},
    \dots,
    - \gamma_{2},
    - \gamma_{1}
    \;|\;
    \dot{r}_{\gamma_{r}},
    \dot{(r-1)}_{\gamma_{r-1} - \gamma_{r}},
    \dots,
    \dot{1}_{\gamma_{1} - \gamma_{2}},
    \dot{0}_{n - \gamma_{1}}
  ).
  \label{eq:Defnoflambdagamma}
\end{eqnarray}
Then, for our specific representations with
$\Lambda=(\dot{0}_m|\dot{\alpha}_n)$, modulo the comments on limiting
$\alpha$ in \S\ref{sec:Uqglmnalphareps} (viz $\alpha$ must be real and
either $\alpha>1-n$ or $\alpha<1-m$), we have the following
irreducible decomposition of $V_{\Lambda}$:
\begin{eqnarray*}
  V_{\Lambda}
  =
  \bigoplus_{[\gamma]}
    V^0_{\Lambda+\lambda_{\gamma}},
\end{eqnarray*}
where the sum is over all possible allowable Young diagrams, and
$V^0_{\Lambda+\lambda_{\gamma}}$ is a
$U_q[gl(m)\oplus gl(n)]$ module of highest weight
$\Lambda+\lambda_{\gamma}$ and $\mathbb{Z}$ graded level $N$.  More
interestingly, we also have the following decomposition:
\begin{eqnarray}
  V_{\Lambda} \otimes V_{\Lambda}
  =
  \bigoplus_{[\gamma]}
    V_{2\Lambda+\lambda_{\gamma}},
  \label{eq:FullUqglmndecomp}
\end{eqnarray}
where again the sum is over all possible allowable Young diagrams, and
$V_{2\Lambda+\lambda_{\gamma}}$ is a $U_q[gl(m|n)]$ module of
highest weight $2\Lambda+\lambda_{\gamma}$ and $\mathbb{Z}$ graded
level $N$.

And now, we introduce an abuse of notation. Instead of the explicit
(\ref{eq:FullUqglmndecomp}), we shall often write:
\begin{eqnarray}
  V \otimes V
  =
  \bigoplus_{k}
    V_k,
  \label{eq:Abuseofnotation}
\end{eqnarray}
where the submodule $V_{k}$ has highest weight $\lambda_k$, and the
generic index $k$ runs over some appropriate index set.
In the \emph{special} case $n=1$, the decomposition of
(\ref{eq:FullUqglmndecomp}) becomes:
\begin{eqnarray}
  V \otimes V
  =
  \bigoplus_{k=0}^{m}
    V_k,
  \label{eq:Uqglm1decomp}
\end{eqnarray}
where $V_k$ has highest weight
$\lambda_k=(\dot{0}_{m-k},\dot{(-1)}_{k}\,|\,2\alpha+k)$, thus the
submodules $V_{0},V_{1},\dots,V_{m}$ are ordered by increasing
$\mathbb{Z}$ graded level (which, in this case, is in fact, $k$).  We
shall be using (\ref{eq:Uqglm1decomp}) to illustrate specific examples
later.

In the following subsections we describe the construction of
orthonormal weight bases
$\mathfrak{B}_k=\{\ket{\Psi^k_l}\}_{l=1}^{\mathrm{dim}(V_k)}$ for each
$V_k$, where $\mathrm{span}(\bigcup_k\mathfrak{B}_k)=V\otimes V$.
Fixing $k$, each $\ket{\Psi^k_l}$ is a linear combination of terms of
the form $\theta_{ij}\ketket{i}{j}$.  As we may safely define the duals
$\bra{\Psi^k_l}\triangleq\ket{\Psi^k_l}^\dagger$ (as expansion of the
conjugate is possible using (\ref{eq:defnofdualonVoV})), it is thus
meaningful to orthonormalise the $\ket{\Psi^k_l}$.

Firstly, in \S\ref{sec:Ahighestweightvector}, we determine a highest
weight vector $\ket{\Psi^k_1}$ for $V_k$; which we use as the starting
point for the KIMC.  Next, in \S\ref{sec:AbasisBk0forVk0} we build a
basis $\mathfrak{B}_k^0$ for the even subalgebra submodule $V_k^0$.  In
this case, as distinct from that of $V^0$ in
\S\ref{sec:Uqglmnalphareps}, $V^k_0$ may not be one-dimensional
(although it contains no weight multiplicities),  so this construction
is nontrivial, although it turns out to be quite straightforward.
Lastly, in \S\ref{sec:overlineBkforVk} and
\S\ref{sec:OrthonormalisingoverlineBkintoBk}, we construct
$\mathfrak{B}_k$ by the actions of the $2^{mn}$ possible combinations
of ordered nonrepeated lowering generators on $\mathfrak{B}_k^0$. A
subtlety in this case is that $V_k$ in general contains weight
multiplicities, so that we must employ a Gram--Schmidt process to
orthonormalise it.


\subsection{A highest weight vector $\ket{\Psi^k_1}$ for $V_k$}
\label{sec:Ahighestweightvector}

We begin the construction of $\mathfrak{B}_k$ with the
deduction of a highest weight vector $\ket{\Psi^k_1}$, of
weight $\lambda_k$.  As we know the weights of the $\ket{i}$, using
(\ref{eq:wtsofketket}), we may immediately write down:
\begin{eqnarray}
  \ket{\Psi^k_1}
  =
  {\textstyle \sum_{ij}}
    \theta_{ij}
    \ketket{i}{j},
  \label{eq:Vkhwv}
\end{eqnarray}
where the sum is only over $i,j$ such that
$
  \mathrm{wt}(\ket{\Psi^k_1})
  =
  \lambda_k
  =
  \mathrm{wt}(\ket{i})+\mathrm{wt}(\ket{j})
$
is satisfied (i.e. we don't know in advance how many terms there are in
the sum).  The coefficients $\theta_{ij}$ are scalar expressions in $q$
and $\alpha$, which we shall determine by the following.

\begin{enumerate}
\item
  We demand that $\ket{\Psi^k_1}$ be annihilated by the actions of (the
  coproducts of) all raising generators. (Actually, it is sufficient
  that it be annihilated by all \emph{simple} raising generators.) As
  $\mathrm{dim}(V_k^0)\neq 1$ in general, $\ket{\Psi^k_1}$ may not be
  annihilated by the actions of \emph{even} lowering generators.

\item
  We demand that $\ket{\Psi^k_1}$ be normalised.
  Examination of (\ref{eq:Vkhwv}) shows that, using
  (\ref{eq:innerproductonVoV}), the magnitude of $\ket{\Psi^k_1}$ is
  $
    \left(
      {\textstyle \sum_{ij}}
        (-)^{[i][j]}
        \theta_{ij}^*
        \theta_{ij}
    \right)^{\frac{1}{2}}
  $,
  where the sum is again only over appropriate $i,j$.  As the
  $\theta_{ij}$ are only determined up to arbitrary phase factors, we
  may replace $\theta_{ij}^* \theta_{ij}$ with $\theta_{ij}^2$, hence
  we demand:
  \begin{eqnarray}
    {\textstyle \sum_{ij}}
      (-)^{[i][j]}
      \theta_{ij}^2
    =
    1.
    \label{eq:Sumofsquares}
  \end{eqnarray}

\end{enumerate}
It turns out that these considerations always yield exactly enough
constraints to uniquely determine $\ket{\Psi^k_1}$ (well, up to an
unimportant phase factor).

Before proceeding, we observe that $\ket{\Psi^0_1}$, the highest weight
vector of the first module $V_0$, is necessarily $\ketket{1}{1}$.  This
provides a check on the methods used to determine the other
$\ket{\Psi^k_1}$.  More generally, we must set up a linear system to
determine the coefficients $\theta_{ij}$, and again,
\textsc{Mathematica} is well-suited to this.


\subsubsection{Illustration: $\ket{\Psi^2_1}$ for $U_q[gl(3|1)]$}

We now use $\ket{\Psi^2_1}$ for $U_q[gl(3|1)]$ to illustrate the entire
process.  Using the notation of
(\ref{eq:Uqglm1decomp}), we have $4$ submodules $V_k$; their highest
weights, dimensions (obtained from
(\ref{eq:KacWeylDimension})) and suitable $\ket{\Psi^k_1}$ are
presented in Table \ref{tab:3hwv}.  Note that for $n=1$, for
$\Lambda=(\dot{0}_m|\alpha)$, in fact (\ref{eq:KacWeylDimension})
degenerates to
$
  \mathrm{dim}(V_k^0)
  =
  \left(
    \begin{array}{@{\hspace{0pt}}c@{\hspace{0pt}}}
      \\[-5mm]
      {\scriptstyle m} \\[-1mm]
      {\scriptstyle k}
    \end{array}
  \right)
$.

\begin{table}[ht]
  \begin{centering}
  \begin{tabular}{
    cr@{\hspace{0pt}}r@{\hspace{2pt}}r@{\hspace{2pt}}r
    @{\hspace{4pt}}c@{\hspace{4pt}}l@{\hspace{0pt}}rcc
    c
  }
    & & & & & & & & \multicolumn{2}{c}{$\mathrm{dim}$} & \\
    $k$ & \multicolumn{7}{c}{$\lambda_k$} &
          $V_{k}^0$ &
          $V_{k}$ &
          $\ket{\Psi^k_1}$
       \\[1mm]
    \hline
    \\[-3mm]
    $0$ & $($ & $ 0,$ & $ 0,$ & $ 0$& $|$ & $ 2\alpha  $ & $)$ & $1$ & $~8$ &
      $\theta_{11} \ketket{1}{1}$
      \\
    \\[-2mm]
    $1$ & $($ & $ 0,$ & $ 0,$ & $-1$& $|$ & $ 2\alpha+1$ & $)$ & $3$ & $24$ &
      $
        \theta_{12} \ketket{1}{2}
        +
        \theta_{21} \ketket{2}{1}
      $
      \\
    \\[-2mm]
    $2$ & $($ & $ 0,$ & $-1,$ & $-1$& $|$ & $ 2\alpha+2$ & $)$ & $3$ & $24$ &
      $
      \begin{array}{r}
        \theta_{15} \ketket{1}{5}
        +
        \theta_{51} \ketket{5}{1}
        \\
        +
        \theta_{23} \ketket{2}{3}
        +
        \theta_{32} \ketket{3}{2}
      \end{array}
      $
      \\
    \\[-2mm]
    $3$ & $($ & $-1,$ & $-1,$ & $-1$& $|$ & $ 2\alpha+3$ & $)$ & $1$ & $~8$ &
      $
      \begin{array}{r}
        \theta_{18} \ketket{1}{8}
        +
        \theta_{81} \ketket{8}{1}
        \\
        +
        \theta_{27} \ketket{2}{7}
        +
        \theta_{72} \ketket{7}{2}
        \\
        +
        \theta_{36} \ketket{3}{6}
        +
        \theta_{63} \ketket{6}{3}
        \\
        +
        \theta_{45} \ketket{4}{5}
        +
        \theta_{54} \ketket{5}{4}
      \end{array}
      $
  \end{tabular}
  \caption{
    The tensor product submodules for the $U_q[gl(3|1)]$ case.
  }
  \label{tab:3hwv}
  \end{centering}
\end{table}

From Table \ref{tab:3hwv}, for $\ket{\Psi^2_1}$, we must determine $4$
coefficients:
\begin{eqnarray}
  \ket{\Psi^2_1}
  =
  \theta_{15} \ketket{1}{5}
  +
  \theta_{51} \ketket{5}{1}
  +
  \theta_{23} \ketket{2}{3}
  +
  \theta_{32} \ketket{3}{2}.
  \label{eq:Psi21defn}
\end{eqnarray}
The $U_q[gl(3|1)]$ simple raising generator set is
$
  SRGS
  =
  \{
    {E^{1}}_{2},
    {E^{2}}_{3},
    {E^{3}}_{4}
  \}
$,
and from (\ref{eq:UqglmnCoproduct}b), the coproducts are:
\begin{eqnarray}
  \left.
  \begin{array}{rcl}
    \Delta({E^{1}}_{2})
    & = &
    ({E^{1}}_{2} \otimes \overline{K}_1^{\frac{1}{2}} K_{2}^{\frac{1}{2}})
    +
    (K_{1}^{\frac{1}{2}} \overline{K}_2^{\frac{1}{2}} \otimes {E^{1}}_{2})
    \\
    \Delta({E^{2}}_{3})
    & = &
    ({E^{2}}_{3} \otimes \overline{K}_2^{\frac{1}{2}} K_{3}^{\frac{1}{2}})
    +
    (K_{2}^{\frac{1}{2}} \overline{K}_3^{\frac{1}{2}} \otimes {E^{2}}_{3})
    \\
    \Delta({E^{3}}_{4})
    & = &
    ({E^{3}}_{4} \otimes \overline{K}_3^{\frac{1}{2}} K_4^{\frac{1}{2}})
    +
    (K_{3}^{\frac{1}{2}} \overline{K}_4^{\frac{1}{2}} \otimes {E^{3}}_{4}),
  \end{array}
  \hspace{68pt}
  \right\}
  \label{eq:3coprodSRGS}
\end{eqnarray}
so the actions that we want are:
\begin{eqnarray*}
  \left(
    \begin{array}{c}
    \Delta({E^{1}}_{2})
    \\
    \Delta({E^{2}}_{3})
    \\
    \Delta({E^{3}}_{4})
    \end{array}
  \right)
  \cdot
  \ket{\Psi^2_1}
  =
  \left(
    \begin{array}{c}
      0
      \\
      0
      \\
      0
    \end{array}
  \right).
\end{eqnarray*}

To evaluate these products, we take the known matrix elements of the
underlying representation (\ref{eq:Uqgl31matrixelements}), and
substitute these into the evaluations of the coproducts
(\ref{eq:3coprodSRGS}). Thus, for the example
$\Delta({E^{3}}_{4}) \cdot \ket{\Psi^2_1}$, we find:
\begin{eqnarray*}
  & &
  \hspace{-40pt}
  \Delta({E^{3}}_{4})
  \cdot
  \ket{\Psi^2_1}
  =
  \left(
    ({E^{3}}_{4}\otimes\overline{K}_3^{\frac{1}{2}} K_4^{\frac{1}{2}})
    +
    (K_{3}^{\frac{1}{2}}\overline{K}_4^{\frac{1}{2}}\otimes {E^{3}}_{4})
  \right)
  \cdot
  \\
  & &
  \qquad \qquad
  \left(
    \theta_{15} \ket{1}\otimes\ket{5}
    +
    \theta_{51} \ket{5}\otimes\ket{1}
    +
    \theta_{23} \ket{2}\otimes\ket{3}
    +
    \theta_{32} \ket{3}\otimes\ket{2}
  \right).
\end{eqnarray*}

To illustrate the multiplication:
\begin{eqnarray*}
  \hspace{-25pt}
  \left(
    {E^{3}}_{4} \otimes \overline{K}_3^{\frac{1}{2}} K_4^{\frac{1}{2}}
  \right)
  \cdot
  \left(
    \ket{5}\otimes\ket{1}
  \right)
  \hspace{-2pt}
  \stackrel{(\ref{eq:XYxketketdefn})}{=}
  \hspace{-2pt}
  {E^{3}}_{4} \cdot \ket{5}
  \otimes
  \overline{K}_3^{\frac{1}{2}} K_4^{\frac{1}{2}} \cdot \ket{1}
  \hspace{-1pt}
  \stackrel{(\ref{eq:Uqgl31matrixelements}i)}{=}
  \hspace{-1pt}
  \overline{q}^{\frac{\alpha}{2}}
  A_{1}^{\frac{1}{2}} \ket{3}
  \otimes
  \ket{1}.
\end{eqnarray*}

At this point, to save space, we introduce a little more notation, which
eliminates the $q$ brackets altogether:
\begin{eqnarray}
  A_{i,j}^z
  & \triangleq &
  ([ i \alpha + j ]_q)^z,
  \qquad
  \mathrm{where~} z\in \{{\textstyle \frac{1}{2}}, 1\}
  \label{eq:DefinitionofAijz}
  \\
  C_{i,j}
  & \triangleq &
  q^{i \alpha+j}
  +
  \overline{q}^{i \alpha+j},
  \label{eq:DefinitionofCij}
\end{eqnarray}
where $i,j\geqslant 0$.
In these expressions, if $i=1$ we shall simply omit it, that is, we
intend:
$
  C_{j}
  \triangleq
  C_{1,j}
  =
  q^{\alpha+j}
  +
  \overline{q}^{\alpha+j}
$
and
$
  A_j^z
  \triangleq
  A_{1,j}^z
  =
  ([\alpha+j]_q)^z
$.
Occasionally, we will write $\overline{A}_{i,j}^z=A_{i,j}^{-z}$ and
$\overline{C}_{i,j}=\left(C_{i,j}\right)^{-1}$.

Using this notation, we have:
\begin{eqnarray*}
  & &
  \hspace{-45pt}
  \Delta({E^{3}}_{4})
  \cdot
  \ket{\Psi^2_1}
  =
  \left(
    ({E^{3}}_{4} \otimes \overline{K}_3^{\frac{1}{2}} K_4^{\frac{1}{2}})
    +
    (K_{3}^{\frac{1}{2}} \overline{K}_4^{\frac{1}{2}} \otimes {E^{3}}_{4})
  \right)
  \cdot
  \\
  & &
  \qquad \quad
  \left(
    \theta_{15} \ket{1}\otimes\ket{5}
    +
    \theta_{51} \ket{5}\otimes\ket{1}
    +
    \theta_{23} \ket{2}\otimes\ket{3}
    +
    \theta_{32} \ket{3}\otimes\ket{2}
  \right)
  \\
  & \stackrel{(\ref{eq:Uqgl31matrixelements})}{=} &
  \left(
    \theta_{15}
    q^{\frac{\alpha}{2}}
    A_1^{\frac{1}{2}}
    +
    \theta_{23}
    \overline{q}^{\frac{\alpha}{2}+\frac{1}{2}}
    A_0^{\frac{1}{2}}
  \right)
  \left(
    \ket{1}\otimes\ket{3}
  \right)
  \\
  & & \qquad
  +
  \left(
    \theta_{51}
    \overline{q}^{\frac{\alpha}{2}}
    A_1^{\frac{1}{2}}
    -
    \theta_{32}
    q^{\frac{\alpha}{2}+\frac{1}{2}}
    A_0^{\frac{1}{2}}
  \right)
  \left(
    \ket{3}\otimes\ket{1}
  \right),
\end{eqnarray*}
and, altogether, for the three $SRGS$ generators, we have:
\begin{eqnarray*}
  \hspace{-0pt}
  \left(
    \begin{array}{@{\hspace{0pt}}c@{\hspace{0pt}}}
    \Delta({E^{1}}_{2})
    \\
    \Delta({E^{2}}_{3})
    \\
    \Delta({E^{3}}_{4})
    \end{array}
  \right)
  \cdot
  \ket{\Psi^2_1}
  =
  \left(
    \begin{array}{@{\hspace{0pt}}c@{\hspace{0pt}}}
    0
    \\
    \left(
      -
      \theta_{23}
      q^{\frac{1}{2}}
      -
      \theta_{32}
      \overline{q}^{\frac{1}{2}}
    \right)
    \left(
      \ket{2}\otimes\ket{2}
    \right)
    \\
    \left(
      \begin{array}{@{\hspace{0pt}}r@{\hspace{0pt}}}
        \left(
          \theta_{15}
          q^{\frac{\alpha}{2}}
          A_1^{\frac{1}{2}}
          +
          \theta_{23}
          \overline{q}^{\frac{\alpha}{2}+\frac{1}{2}}
          A_0^{\frac{1}{2}}
        \right)
        \left(
          \ket{1}\otimes\ket{3}
        \right)
        \\
        +
        \left(
          \theta_{51}
          \overline{q}^{\frac{\alpha}{2}}
          A_1^{\frac{1}{2}}
          -
          \theta_{32}
          q^{\frac{\alpha}{2}+\frac{1}{2}}
          A_0^{\frac{1}{2}}
        \right)
        \left(
          \ket{3}\otimes\ket{1}
        \right)
      \end{array}
    \right)
    \end{array}
  \right).
\end{eqnarray*}
As each component of the RHS must be zero, and the
$\ket{i}\otimes\ket{j}$ are linearly independent, we thus obtain a net
$3$ \emph{linear} constraints on the $\theta_{ij}$ from this set:
\begin{eqnarray*}
  -
  q^{\frac{1}{2}}
  \theta_{23}
  -
  \overline{q}^{\frac{1}{2}}
  \theta_{32}
  & = &
  0
  \\
  q^{\frac{\alpha}{2}}
  A_1^{\frac{1}{2}}
  \theta_{15}
  +
  \overline{q}^{\frac{\alpha}{2}+\frac{1}{2}}
  A_0^{\frac{1}{2}}
  \theta_{23}
  & = &
  0
  \\
  \overline{q}^{\frac{\alpha}{2}}
  A_1^{\frac{1}{2}}
  \theta_{51}
  -
  q^{\frac{\alpha}{2}+\frac{1}{2}}
  A_0^{\frac{1}{2}}
  \theta_{32}
  & = &
  0,
\end{eqnarray*}
better written in matrix form as:
\begin{eqnarray}
  \left(
    \begin{array}{cccc}
      0
      &
      0
      &
      -
      q^{\frac{1}{2}}
      &
      -
      \overline{q}^{\frac{1}{2}}
      \\
      q^{\frac{\alpha}{2}}
      A_1^{\frac{1}{2}}
      &
      0
      &
      \overline{q}^{\frac{\alpha}{2}+\frac{1}{2}}
      A_0^{\frac{1}{2}}
      &
      0
      \\
      0
      &
      \overline{q}^{\frac{\alpha}{2}}
      A_1^{\frac{1}{2}}
      &
      0
      &
      -
      q^{\frac{\alpha}{2}+\frac{1}{2}}
      A_0^{\frac{1}{2}}
    \end{array}
  \right)
  \cdot
  \left(
    \begin{array}{c}
      \theta_{15} \\
      \theta_{51} \\
      \theta_{23} \\
      \theta_{32}
    \end{array}
  \right)
  =
  \left(
    \begin{array}{c}
      0 \\
      0 \\
      0
    \end{array}
  \right).
  \label{eq:3linearsystem}
\end{eqnarray}
Note that we have left the system exactly as supplied by the raw
action of the raising generators. A human calculator might
delete superfluous signs, and perhaps do some factorisation,
but our \textsc{Mathematica} code would require explicit instructions
for this finicky and unnecessary work, so we omit it.

Thus, application of the requirement that $\ket{\Psi^2_1}$ be
annihilated by the simple raising generators yields a linear system of
3 equations in 4 variables.  A final constraint to completely determine
the variables is now obtained by requiring that $\ket{\Psi^2_1}$ be
normalised.

In our example, using (\ref{eq:Sumofsquares}) and (\ref{eq:Psi21defn}),
and recalling the gradings of the basis vectors (in Table
\ref{tab:3wtsgradings}), we thus have the \emph{nonlinear} constraint:
\begin{eqnarray}
  \theta_{15}^2
  +
  \theta_{51}^2
  -
  \theta_{23}^2
  -
  \theta_{32}^2
  =
  1.
  \label{eq:3normofk3}
\end{eqnarray}
Combining the information in (\ref{eq:3linearsystem}) and
(\ref{eq:3normofk3}), we are able to solve for the unknowns
$\theta_{ij}$, uniquely up to the usual overall phase factor.
In practice, we may actually avoid the use of
(\ref{eq:3normofk3}), by first feeding (\ref{eq:3linearsystem}) to the
\textsc{Mathematica} equation solver, which returns us an answer
with a free parameter (the first unknown $\theta_{ij}$).%
\footnote{
  In fact, it returns us \emph{two} answers, differing by a (phase)
  factor of $-1$. We judiciously choose to ignore the second one.
}
Setting that free parameter to unity, we obtain a suitable unnormalised
$\ket{\Psi^2_1}$, which we may immediately normalise.  To illustrate,
we find:
\small
\begin{eqnarray*}
  \hspace{-7pt}
  \ket{\Psi^2_1}
  =
  \frac{
    \overline{q}^{\alpha+1}
    A_0^\frac{1}{2}
  }{
    C_1^\frac{1}{2}
    A_{2,1}^\frac{1}{2}
  }
  \ketket{1}{5}
  +
  \frac{
    q^{\alpha+1}
    A_0^\frac{1}{2}
  }{
    C_1^\frac{1}{2}
    A_{2,1}^\frac{1}{2}
  }
  \ketket{5}{1}
  -
  \frac{
    \overline{q}^\frac{1}{2}
    A_1^\frac{1}{2}
  }{
    C_1^\frac{1}{2}
    A_{2,1}^\frac{1}{2}
  }
  \ketket{2}{3}
  +
  \frac{
    q^\frac{1}{2}
    A_1^\frac{1}{2}
  }{
    C_1^\frac{1}{2}
    A_{2,1}^\frac{1}{2}
  }
  \ketket{3}{2}.
\end{eqnarray*}
\normalsize

To complete the results, we have obtained the highest weight vectors
for each of the $4$ submodules for the $U_q[gl(3|1)]$ case. These are,
after a little factoring:
\small
\begin{eqnarray*}
  \hspace{-30pt}
  \begin{array}{r@{\hspace{3pt}}c@{\hspace{3pt}}l}
    \ket{\Psi^0_1}
    &\hspace{0pt}=\hspace{0pt}&
    \ketket{1}{1}
    \\
    \ket{\Psi^1_1}
    &\hspace{0pt}=\hspace{0pt}&
    \overline{C}_0^{\frac{1}{2}}
    \left(
      q^{\frac{\alpha}{2}}
      \ketket{2}{1}
      -
      \overline{q}^{\frac{\alpha}{2}}
      \ketket{1}{2}
    \right)
    \\
    \ket{\Psi^2_1}
    &\hspace{0pt}=\hspace{0pt}&
    \overline{C}_1^{\frac{1}{2}}
    \overline{A}_{2,1}^{\frac{1}{2}}
    \left(
      \hspace{-2pt}
      \overline{q}^{\alpha+1}
      A_0^\frac{1}{2}
      \ketket{1}{5}
      +
      q^{\alpha+1}
      A_0^\frac{1}{2}
      \ketket{5}{1}
      -
      \overline{q}^{\frac{1}{2}}
      A_1^\frac{1}{2}
      \ketket{2}{3}
      +
      q^\frac{1}{2}
      A_1^\frac{1}{2}
      \ketket{3}{2}
      \hspace{-3pt}
    \right)
    \\
    \ket{\Psi^3_1}
    &\hspace{0pt}=\hspace{0pt}&
    \overline{C}_1^{\frac{1}{2}}
    \overline{C}_2^{\frac{1}{2}}
    \overline{A}_{2,3}^{\frac{1}{2}}
    \times
    \\
    & &
    \hspace{-24pt}
    \left(
      \begin{array}{@{\hspace{0pt}}l@{\hspace{5pt}}l@{\hspace{0pt}}}
        - \;
        \overline{q}^{\frac{3\alpha}{2}+3}
        A_0^\frac{1}{2}
        \ketket{1}{8}
        +
        q^{\frac{3\alpha}{2}+3}
        A_0^\frac{1}{2}
        \ketket{8}{1}
        &
        + \;
        \overline{q}^{\frac{\alpha}{2}+2}
        A_2^\frac{1}{2}
        \ketket{2}{7}
        -
        q^{\frac{\alpha}{2}+2}
        A_2^\frac{1}{2}
        \ketket{7}{2}
        \\
        \hspace{7pt}
        - \;
        \overline{q}^{\frac{\alpha}{2}+1}
        A_2^\frac{1}{2}
        \ketket{3}{6}
        +
        q^{\frac{\alpha}{2}+1}
        A_2^\frac{1}{2}
        \ketket{6}{3}
        &
        + \;
        \overline{q}^{\frac{\alpha}{2}}
        A_2^\frac{1}{2}
        \ketket{4}{5}
        -
        q^{\frac{\alpha}{2}}
        A_2^\frac{1}{2}
        \ketket{5}{4}
      \end{array}
    \right).
  \end{array}
\end{eqnarray*}
\normalsize
Observe the presence of \emph{$q$ graded symmetric combinations} of
$\ketket{i}{j}$ and $\ketket{j}{i}$ in these expressions, viz patterns
of the form
$
  q^{x}
  \ketket{i}{j}
  \pm
  \overline{q}^{x}
  \ketket{j}{i}
$.
This feature
appears repeatedly throughout the bases for the $V_k$.

So, at this stage, we have described how to construct normalised
highest weight vectors $\ket{\Psi^k_1}$ for each $V_k$. An interesting
outstanding point about our process is that the two demands that
$\ket{\Psi^k_1}$ be annihilated by the $m+n-1$ simple raising
generators and that it be normalised, yield \emph{exactly} enough
constraints to determine it uniquely (up to a phase factor). The reason
for this balance lies buried in the combinatorics of the ways the
weights of the underlying module $V$ can be added to yield the weight
$\lambda_k$.

\pagebreak


\subsection{A basis $\mathfrak{B}_k^0$ for $V_k^0$}
\label{sec:AbasisBk0forVk0}

Having determined $\ket{\Psi^k_1}$, we now apply the first stage of the
KIMC to construct the basis $\mathfrak{B}_k^0$.  That is, we construct
basis vectors $\ket{\Psi^k_j}$ by the repeated action of the $m+n-2$
even simple lowering generators (a set which we call $ESLGS$)%
\footnote{%
  We might just as well use $ELGS$, the full set of even lowering
  generators. The tradeoff is that whilst the coproducts of the
  nonsimple generators are more complex, there should be less levels to
  calculate.
}
on the known $\ket{\Psi^k_1}$.  In our case, as $V_k^0$ contains no
weight multiplicities, vectors of distinct weights created in this way
will naturally be orthogonal.  At each stage, we must also check to see
if newly minted vectors are scalar multiples of previously found ones.
To facilitate this, we will normalise each vector as we create it, and
we will also maintain our list of vectors in order of decreasing
weight.  In fact, (\ref{eq:KacWeylDimension}) tells us
$\mathrm{dim}(V_k^0)$, but we will build $\mathfrak{B}_k^0$ as if we
didn't know this.

Recall that in \S\ref{sec:AnorthonrmalbasisB}, we created $\mathbb{Z}$
graded levels of $B$ by repeated applications of the set of odd
lowering generators to the set $\{\ket{1}\}$.  Here, each application
of $ESLGS$ generates an \emph{ungraded} level, which we shall call
$L_i$, where $L_0=\{\ket{\Psi^k_1}\}$. We describe the process in the
following algorithm:

\begin{figure}[h]
  \begin{tabbing}
    \hspace{40pt} \= $i:=0$ \= \\
    \>$L_0:=\{\ket{\Psi^k_1}\}$ \\
    \>while $L_i\neq\emptyset$ \\
          \>\> $L_{i+1}:=(\Delta(ESLGS)\cdot L_i)\;\backslash\;\{0\}$ \\
          \>\> normalise $L_{i+1}$ \\
          \>\> $L_{i+1}:=L_{i+1}\; \backslash \; L$ \\
          \>\> $L := L \cup L_{i+1}$ \\
          \>\> sort $L$ by decreasing weight \\
          \>\> increment $i$ \\
    \>$\mathfrak{B}_k^0:=L$
  \end{tabbing}
\end{figure}

Note that we have taken notational liberties by writing
$\Delta(ESLGS)\cdot L_i$. We do this as it is natural to apply
functions to lists in \textsc{Mathematica}.

Of particular interest here is that the evaluation of the
algebra-module action in the tensor product case utilises the
information encoded in the matrix elements of the underlying
representation to determine when basis vectors are annihilated.  That
is, as different from \S\ref{sec:Uqglmnalphareps}, we do not have to
explicitly implement the annihilation rules dictated by the KIMC.
This observation carries over into the following subsections.

\pagebreak


\subsubsection{Illustration: $\mathfrak{B}_1^0$ for $U_q[gl(3|1)]$}
\label{sec:mathfrakB10}

We illustrate the process by constructing the $3$ dimensional
$\mathfrak{B}_1^0$ for $U_q[gl(3|1)]$. Recall that we
determined:
\begin{eqnarray}
  \ket{\Psi^1_1}
  =
  \overline{C}_0^{\frac{1}{2}}
  \left(
    q^{\frac{\alpha}{2}}
    \ketket{2}{1}
    -
    \overline{q}^{\frac{\alpha}{2}}
    \ketket{1}{2}
  \right),
  \label{eq:Psi11}
\end{eqnarray}
hence
$
  L_0
  =
  \{
    \overline{C}_0^{\frac{1}{2}}
    \left(
      q^{\frac{\alpha}{2}}
      \ketket{2}{1}
      -
      \overline{q}^{\frac{\alpha}{2}}
      \ketket{1}{2}
    \right)
  \}
$.
We also have
$
  ESLGS
  =
  \left\{
    {E^{2}}_{1},
    {E^{3}}_{2}
  \right\}
$,
for which, from (\ref{eq:UqglmnCoproduct}a), we have the coproducts:
\begin{eqnarray}
  \left.
  \begin{array}{rcl}
    \Delta({E^{2}}_{1})
    & = &
    ({E^{2}}_{1}\otimes K_{2}^{\frac{1}{2}}\overline{K}_1^{\frac{1}{2}})
    +
    (\overline{K}_2^{\frac{1}{2}}K_{1}^{\frac{1}{2}}\otimes {E^{2}}_{1})
    \\
    \Delta({E^{3}}_{2})
    & = &
    ({E^{3}}_{2}\otimes K_{3}^{\frac{1}{2}}\overline{K}_2^{\frac{1}{2}})
    +
    (\overline{K}_3^{\frac{1}{2}}K_{2}^{\frac{1}{2}}\otimes{E^{3}}_{2}).
  \end{array}
  \hspace{69pt}
  \right\}
  \label{eq:3coprodESLGS}
\end{eqnarray}
Applying the operators of (\ref{eq:3coprodESLGS}) to
(\ref{eq:Psi11}), we obtain two vectors:
\begin{eqnarray*}
  \left\{
    0,
    \overline{C}_0^{\frac{1}{2}}
    \left(
      \overline{q}^{\frac{\alpha}{2}}
      \ketket{1}{3}
      -
      q^{\frac{\alpha}{2}}
      \ketket{3}{1}
    \right)
  \right\}.
\end{eqnarray*}
We discard the $0$, and find that the second vector
is already normalised and also not found in $L_0$, so we have
$
  L_1
  =
  \{
    \overline{C}_0^{\frac{1}{2}}
    \left(
      \overline{q}^{\frac{\alpha}{2}}
      \ketket{1}{3}
      -
      q^{\frac{\alpha}{2}}
      \ketket{3}{1}
    \right)
  \}
$.
As the components of $L_0$ and $L_1$ are independent, and indeed
ordered, we have, at this stage:
\begin{eqnarray*}
  L
  =
  \{
    \overline{C}_0^{\frac{1}{2}}
    \left(
      q^{\frac{\alpha}{2}}
      \ketket{2}{1}
      -
      \overline{q}^{\frac{\alpha}{2}}
      \ketket{1}{2}
    \right),
    \overline{C}_0^{\frac{1}{2}}
    \left(
      \overline{q}^{\frac{\alpha}{2}}
      \ketket{1}{3}
      -
      q^{\frac{\alpha}{2}}
      \ketket{3}{1}
    \right)
  \}.
\end{eqnarray*}
Repeating this process on $L_1$, we find
$
  L_2
  =
  \left\{
    \overline{C}_0^{\frac{1}{2}}
    \left(
      q^{\frac{\alpha}{2}}
      \ketket{4}{1}
      -
      \overline{q}^{\frac{\alpha}{2}}
      \ketket{1}{4}
    \right)
  \right\}
$.
As before, we normalise (again already OK) and install this vector in
its rightful position in our collection $L$, checking first to see if
we have already met it. We thus have, after ordering:
\begin{eqnarray*}
  L
  =
  \left\{
    \begin{array}{l}
      \overline{C}_0^{\frac{1}{2}}
      \left(
        q^{\frac{\alpha}{2}}
        \ketket{2}{1}
        -
        \overline{q}^{\frac{\alpha}{2}}
        \ketket{1}{2}
      \right)
      \\
      \overline{C}_0^{\frac{1}{2}}
      \left(
        \overline{q}^{\frac{\alpha}{2}}
        \ketket{1}{3}
        -
        q^{\frac{\alpha}{2}}
        \ketket{3}{1}
      \right)
      \\
      \overline{C}_0^{\frac{1}{2}}
      \left(
        q^{\frac{\alpha}{2}}
        \ketket{4}{1}
        -
        \overline{q}^{\frac{\alpha}{2}}
        \ketket{1}{4}
      \right)
    \end{array}
  \right\}.
\end{eqnarray*}
Repeating again, we discover that $L_3=\emptyset$, so the process is
completed. Our orthonormal basis $\mathfrak{B}_1^0$ is thus the above
$L$.

Note that, apart from some factoring, we have left this basis in the
raw form that our \textsc{Mathematica} code yields.  The human
calculator, preferring symmetry, may wish to multiply some of the
vectors by (the phase factor) $-1$, but this is unnecessary for our
purposes. In Appendix \ref{app:Data}, where our results are summarised,
we have make some judicious changes of this nature for readability.

\pagebreak


\subsection{A nonorthogonal basis $\overline{\mathfrak{B}}_k$ for $V_k$}
\label{sec:overlineBkforVk}

At this stage, we have established the orthonormal basis
$\mathfrak{B}_k^0$ for $V_k^0$, and we wish to extend
$\mathfrak{B}_k^0$ to a basis $\mathfrak{B}_k$ for $V_k$.  This process
involves two stages:

\begin{itemize}
\item
  Firstly, we use the KIMC to construct a basis
  $\overline{\mathfrak{B}}_k$ for $V_k$ by the repeated actions of the
  odd lowering generators on the basis $\mathfrak{B}_k^0$, normalising
  and casting out repeated vectors as we go.  This part of this process
  is detailed in this subsection: its appearance is similar to that of
  \S\ref{sec:AbasisBk0forVk0}.

\item
  Unlike $\mathfrak{B}_k^0$ however, the vectors of
  $\overline{\mathfrak{B}}_k$ are not guaranteed to be orthogonal, as
  $V_k$ in general contains some weight multiplicities.  That is,
  distinct vectors of the same weights will generally appear, and these
  are usually nonorthogonal.  To deal with this problem, we apply a
  Gram--Schmidt process to orthonormalise $\overline{\mathfrak{B}}_k$
  into the final $\mathfrak{B}_k$.  To optimise this, we preprocess
  $\overline{\mathfrak{B}}_k$ by ordering its vectors by decreasing
  weight and then partitioning it into weight equivalence classes. We
  then need only apply the Gram--Schmidt process to each equivalence
  class.  The end result is $\mathfrak{B}_k$, the desired orthonormal
  weight basis for $V_k$.  This process is documented in
  \S\ref{sec:OrthonormalisingoverlineBkintoBk}.

\end{itemize}

Thus, we reproduce essentially the same algorithm as that used in
\S\ref{sec:AbasisBk0forVk0}, the essential differences being that the
levels $L_i$ are now $\mathbb{Z}$ graded levels of $V_k$,
and that we act with $OLGS$
(the (full) set of odd lowering generators, see (\ref{eq:OLGS})),
rather than with $ESLGS$.%
\footnote{%
  In this case, if we try to only use $OSLGS$, the set of \emph{simple}
  odd lowering generators, then we miss some of the vectors of each
  level, which are obtained by nonsimple odd lowering generators, i.e.
  products of simple odd lowering generators with even lowering
  generators. The combination of the use $OSLGS$ and $ELGS$ (or, the
  repeated use of $ESLGS$) to build each level would require more
  calculations.
}

\begin{figure}[h]
  \begin{tabbing}
    \hspace{40pt} \= $i:=0$ \= \\
    \>$L_0:=\mathfrak{B}_k^0$ \\
    \>while $L_i\neq\emptyset$ \\
          \>\> $L_{i+1}:=(\Delta(OLGS)\cdot L_i)\;\backslash\;\{0\}$ \\
          \>\> normalise $L_{i+1}$ \\
          \>\> $L_{i+1}:=L_{i+1}\; \backslash \; L$ \\
          \>\> $L := L \cup L_{i+1}$ \\
          \>\> sort $L$ by decreasing weight \\
          \>\> increment $i$ \\
    \>$\overline{\mathfrak{B}}_k:=L$
  \end{tabbing}
\end{figure}

Again, the elements of $\mathfrak{B}_k$ are only unique up to phase
factors, and our code doesn't select these, so the final results
contain various factors of $-1$ that a human calculator would quickly
purge.  Furthermore, the weight ordering covers some elegant symmetries
of the generators. In the results presented in the Appendix, we make
some judicious cosmetic changes for readability.


\subsubsection{Illustration: $\overline{\mathfrak{B}}_1$ for
               $U_q[gl(3|1)]$}
\label{sec:overlinemathfrakB1}

We illustrate the results using $\overline{\mathfrak{B}}_1$ for
$U_q[gl(3|1)]$, continuing the example of \S\ref{sec:mathfrakB10}.
Here, $\overline{\mathfrak{B}}_1$ has $24$ elements, sorted into
equivalence classes of decreasing weight, and judiciously factored.  To
save space, we have written $\nabla\triangleq q+\overline{q}$, and
$\overline{\nabla}\triangleq\nabla^{-1}$.

\renewcommand{\arraystretch}{1.4}
\begin{figure}[htbp]
  \tiny
  \begin{eqnarray*}
    \hspace{-5pt}
    \begin{array}{l}
      \overline{C}_0^{\frac{1}{2}}
      \left(
        -
        \overline{q}^{\frac{\alpha}{2}}
        \ketket{1}{2}
        +
        q^{\frac{\alpha}{2}}
        \ketket{2}{1}
      \right)
      \\[1.0mm]
      \hline
      \\[-2.5mm]
      -\ketket{2}{2}
      \\[1.0mm]
      \hline
      \\[-2.5mm]
      \overline{C}_0^{\frac{1}{2}}
      \left(
        +
        \overline{q}^{\frac{\alpha}{2}}
        \ketket{1}{3}
        -
        q^{\frac{\alpha}{2}}
        \ketket{3}{1}
      \right)
      \\[1.0mm]
      \hline
      \\[-2.5mm]
      \overline{C}_0^{\frac{1}{2}}
      \overline{A}_{2,1}^{\frac{1}{2}}
      \left(
        A_0^{\frac{1}{2}}
        (
          +
          \overline{q}^{\alpha+\frac{1}{2}}
          \ketket{2}{3}
          +
          q^{\alpha+\frac{1}{2}}
          \ketket{3}{2}
        )
        +
        A_1^{\frac{1}{2}}
        (
          +
          \ketket{1}{5}
          -
          \ketket{5}{1}
        )
      \right)
      \\
      \overline{\nabla}^{\frac{1}{2}}
      \left(
        +
        q^{\frac{1}{2}}
        \ketket{2}{3}
        +
        \overline{q}^{\frac{1}{2}}
        \ketket{3}{2}
      \right)
      \\[1.0mm]
      \hline
      \\[-2.5mm]
      \overline{C}_1^{\frac{1}{2}}
      \left(
        -
        q^{\frac{\alpha}{2}+\frac{1}{2}}
        \ketket{2}{5}
        +
        \overline{q}^{\frac{\alpha}{2}+\frac{1}{2}}
        \ketket{5}{2}
      \right)
      \\[1.0mm]
      \hline
      \\[-2.5mm]
      -\ketket{3}{3}
      \\[1.0mm]
      \hline
      \\[-2.5mm]
      \overline{C}_1^{\frac{1}{2}}
      \left(
        +
        q^{\frac{\alpha}{2}+\frac{1}{2}}
        \ketket{3}{5}
        -
        \overline{q}^{\frac{\alpha}{2}+\frac{1}{2}}
        \ketket{5}{3}
      \right)
      \\[1.0mm]
      \hline
      \\[-2.5mm]
      \overline{C}_0^{\frac{1}{2}}
      \left(
        -
        \overline{q}^{\frac{\alpha}{2}}
        \ketket{1}{4}
        +
        q^{\frac{\alpha}{2}}
        \ketket{4}{1}
      \right)
      \\[1.0mm]
      \hline
      \\[-2.5mm]
      \overline{C}_0^{\frac{1}{2}}
      \overline{A}_{2,1}^{\frac{1}{2}}
      \left(
        A_0^{\frac{1}{2}}
        (
          -
          \overline{q}^{\alpha+\frac{1}{2}}
          \ketket{2}{4}
          -
          q^{\alpha+\frac{1}{2}}
          \ketket{4}{2}
        )
        +
        A_1^{\frac{1}{2}}
        (
          -
          \ketket{1}{6}
          +
          \ketket{6}{1}
        )
      \right)
      \\
      \overline{\nabla}^{\frac{1}{2}}
      \left(
        -
        q^{\frac{1}{2}}
        \ketket{2}{4}
        -
        \overline{q}^{\frac{1}{2}}
        \ketket{4}{2}
      \right)
      \\[1.0mm]
      \hline
      \\[-2.5mm]
      \overline{C}_1^{\frac{1}{2}}
      \left(
        +
        q^{\frac{\alpha}{2}+\frac{1}{2}}
        \ketket{2}{6}
        -
        \overline{q}^{\frac{\alpha}{2}+\frac{1}{2}}
        \ketket{6}{2}
      \right)
      \\[1.0mm]
      \hline
      \\[-2.5mm]
      \overline{C}_0^{\frac{1}{2}}
      \overline{A}_{2,1}^{\frac{1}{2}}
      \left(
        A_0^{\frac{1}{2}}
        (
          +
          \overline{q}^{\alpha+\frac{1}{2}}
          \ketket{3}{4}
          +
          q^{\alpha+\frac{1}{2}}
          \ketket{4}{3}
        )
        +
        A_1^{\frac{1}{2}}
       (
          +
          \ketket{1}{7}
          -
          \ketket{7}{1}
       )
      \right)
      \\
      \overline{\nabla}^{\frac{1}{2}}
      \left(
        +
        q^{\frac{1}{2}}
        \ketket{3}{4}
        +
        \overline{q}^{\frac{1}{2}}
        \ketket{4}{3}
      \right)
      \\[1.0mm]
      \hline
      \\[-2.5mm]
      \overline{C}_0^{\frac{1}{2}}
      \overline{C}_1^{\frac{1}{2}}
      \overline{A}_{2,1}^{\frac{1}{2}}
      \left(
        \begin{array}{@{\hspace{0pt}}l@{\hspace{0pt}}l@{\hspace{0pt}}}
          +
          A_2^{\frac{1}{2}}
          \left(
            +
            q^{\frac{\alpha}{2}}
            \ketket{1}{8}
            -
            \overline{q}^{\frac{\alpha}{2}}
            \ketket{8}{1}
          \right)
          &
          +
          A_0^{\frac{1}{2}}
          \left(
            +
            \overline{q}^{\frac{\alpha}{2}+1}
            \ketket{2}{7}
            -
            q^{\frac{\alpha}{2}+1}
            \ketket{7}{2}
          \right)
          \\
          +
          A_0^{\frac{1}{2}}
          \left(
            -
            \overline{q}^{\frac{\alpha}{2}}
            \ketket{3}{6}
            +
            q^{\frac{\alpha}{2}}
            \ketket{6}{3}
          \right)
          &
          +
          A_0^{\frac{1}{2}}
          \left(
            -
            q^{\frac{3\alpha}{2}+1}
            \ketket{4}{5}
            +
            \overline{q}^{\frac{3\alpha}{2}+1}
            \ketket{5}{4}
          \right)
        \end{array}
      \right)
      \\
      \overline{\nabla}^{\frac{1}{2}}
      \overline{C}_1^{\frac{1}{2}}
      \left(
        -
        q^{\frac{\alpha}{2}+1}
        \ketket{3}{6}
        +
        \overline{q}^{\frac{\alpha}{2}+1}
        \ketket{6}{3}
        -
        q^{\frac{\alpha}{2}}
        \ketket{4}{5}
        +
        \overline{q}^{\frac{\alpha}{2}}
        \ketket{5}{4}
      \right)
      \\
      \overline{\nabla}^{\frac{1}{2}}
      \overline{C}_1^{\frac{1}{2}}
      \left(
        -
        q^{\frac{\alpha}{2}+1}
        \ketket{2}{7}
        +
        \overline{q}^{\frac{\alpha}{2}+1}
        \ketket{7}{2}
        -
        q^{\frac{\alpha}{2}}
        \ketket{3}{6}
        +
        \overline{q}^{\frac{\alpha}{2}}
        \ketket{6}{3}
      \right)
      \\[1.0mm]
      \hline
      \\[-2.5mm]
      \overline{C}_1^{\frac{1}{2}}
      \overline{A}_{2,3}^{\frac{1}{2}}
      \left(
        A_1^{\frac{1}{2}}
        (
          -
          \overline{q}^{\frac{1}{2}}
          \ketket{5}{6}
          +
          q^{\frac{1}{2}}
          \ketket{6}{5}
        )
        +
        A_2^{\frac{1}{2}}
        (
          +
          q^{\alpha+1}
          \ketket{2}{8}
          +
          \overline{q}^{\alpha+1}
          \ketket{8}{2}
        )
      \right)
      \\[1.0mm]
      \hline
      \\[-2.5mm]
      \overline{C}_1^{\frac{1}{2}}
      \left(
        +
        q^{\frac{\alpha}{2}+\frac{1}{2}}
        \ketket{3}{7}
        -
        \overline{q}^{\frac{\alpha}{2}+\frac{1}{2}}
        \ketket{7}{3}
      \right)
      \\[1.0mm]
      \hline
      \\[-2.5mm]
      \overline{C}_1^{\frac{1}{2}}
      \overline{A}_{2,3}^{\frac{1}{2}}
      \left(
        A_1^{\frac{1}{2}}
        (
          +
          \overline{q}^{\frac{1}{2}}
          \ketket{5}{7}
          -
          q^{\frac{1}{2}}
          \ketket{7}{5}
        )
        +
        A_2^{\frac{1}{2}}
        (
          -
          q^{\alpha+1}
          \ketket{3}{8}
          -
          \overline{q}^{\alpha+1}
          \ketket{8}{3}
        )
      \right)
      \\[1.0mm]
      \hline
      \\[-2.5mm]
      -\ketket{4}{4}
      \\[1.0mm]
      \hline
      \\[-2.5mm]
      \overline{C}_1^{\frac{1}{2}}
      \left(
        +
        q^{\frac{\alpha}{2}+\frac{1}{2}}
        \ketket{4}{6}
        -
        \overline{q}^{\frac{\alpha}{2}+\frac{1}{2}}
        \ketket{6}{4}
      \right)
      \\[1.0mm]
      \hline
      \\[-2.5mm]
      \overline{C}_1^{\frac{1}{2}}
      \left(
        -
        q^{\frac{\alpha}{2}+\frac{1}{2}}
        \ketket{4}{7}
        +
        \overline{q}^{\frac{\alpha}{2}+\frac{1}{2}}
        \ketket{7}{4}
      \right)
      \\[1.0mm]
      \hline
      \\[-2.5mm]
      \overline{C}_1^{\frac{1}{2}}
      \overline{A}_{2,3}^{\frac{1}{2}}
      \left(
        A_1^{\frac{1}{2}}
        (
          -
          \overline{q}^{\frac{1}{2}}
          \ketket{6}{7}
          +
          q^{\frac{1}{2}}
          \ketket{7}{6}
        )
        +
        A_2^{\frac{1}{2}}
        (
          +
          q^{\alpha+1}
          \ketket{4}{8}
          +
          \overline{q}^{\alpha+1}
          \ketket{8}{4}
        )
      \right)
    \end{array}
  \end{eqnarray*}
  \normalsize
\end{figure}
\renewcommand{\arraystretch}{1}

\pagebreak


\subsection{Orthonormalising $\overline{\mathfrak{B}}_k$ into
            $\mathfrak{B}_k$}
\label{sec:OrthonormalisingoverlineBkintoBk}

At this stage, the KIMC has yielded a nonorthogonal (although
normalised!) $q$ graded symmetric basis $\overline{\mathfrak{B}}_k$ for
$V_k$.  With a view to constructing the projector $P_k$ onto $V_k$
(see \S\ref{sec:Projectors}), we require dual bases for $V_k$ and
its dual $V_k^*$.  From knowledge of $\overline{\mathfrak{B}}_k$ there
are two obvious ways to construct these dual bases:

\begin{itemize}
\item
  We might continue to regard $\overline{\mathfrak{B}}_k$ as
  our basis for $V_k$, and construct a (nonorthogonal) dual tensor
  product basis by the inversion of an overlap (i.e. metric) matrix.
  This process is described for $\overline{\mathfrak{B}}_1$ for
  $U_q[gl(3|1)]$ in \cite{GeGouldZhangZhou:98a}, although those authors
  don't actually implement it.

\item
  Alternatively, we can orthonormalise $\overline{\mathfrak{B}}_k$ into
  $\mathfrak{B}_k$ using a Gram--Schmidt process. The dual basis
  $\mathfrak{B}_k^*$ is then naturally orthonormal, and indeed trivial
  to write down.

\end{itemize}
Implementation of both of these methods has shown that, apart from
being substantially more efficient, the latter method yields more
tractable and symmetric results, so we choose it.  Not only is it more
elegant, but happily, the Gram--Schmidt process also allows us to
maintain the $q$ graded symmetry of the basis vectors.

A basic principle in \emph{numerical} computation is to never invert a
matrix unless absolutely necessary, as the process is both
computationally inefficient and (numerically) unstable. The same
comment about computational inefficiency certainly holds for the
inversion of symbolic matrices. More seriously, for our current
purposes, the inversion can bog down altogether due to difficulties in
the simplification of algebraic expressions; a feature we might call
`symbolic instability'.


\subsubsection{Illustration: $\mathfrak{B}_1$ for $U_q[gl(3|1)]$}

We illustrate the Gram--Schmidt process by continuing the example from
\S\ref{sec:overlinemathfrakB1}, that of $\overline{\mathfrak{B}}_1$ for
$U_q[gl(3|1)]$.  Observe that the partitions of
$\overline{\mathfrak{B}}_1$ include three of size $2$ and one of size
$3$.  To convert $\overline{\mathfrak{B}}_1$ to $\mathfrak{B}_1$, we
must orthogonalise each of those partitions. To illustrate, for the
largest partition, we obtain the $3$ vectors:

\tiny
\begin{eqnarray*}
  \hspace{-10pt}
  \left\{
  \begin{array}{@{\hspace{0pt}}c@{\hspace{0pt}}}
    \overline{C}_0^{\frac{1}{2}}
    \overline{C}_1^{\frac{1}{2}}
    \overline{A}_{2,1}^{\frac{1}{2}}
    \left(
    \begin{array}{@{\hspace{0pt}}c@{\hspace{0pt}}}
      A_2^{\frac{1}{2}}
      \left(
        +
        q^{\frac{\alpha}{2}}
        \ketket{1}{8}
        -
        \overline{q}^{\frac{\alpha}{2}}
        \ketket{8}{1}
      \right)
      +
      A_0^{\frac{1}{2}}
      \left(
        +
        \overline{q}^{\frac{\alpha}{2}+1}
        \ketket{2}{7}
        -
        q^{\frac{\alpha}{2}+1}
        \ketket{7}{2}
      \right)
      \\
      +
      A_0^{\frac{1}{2}}
      \left(
        -
        \overline{q}^{\frac{\alpha}{2}}
        \ketket{3}{6}
        +
        q^{\frac{\alpha}{2}}
        \ketket{6}{3}
      \right)
      +
      A_0^{\frac{1}{2}}
      \left(
        -
        q^{\frac{3\alpha}{2}+1}
        \ketket{4}{5}
        +
        \overline{q}^{\frac{3\alpha}{2}+1}
        \ketket{5}{4}
      \right)
      \end{array}
    \right)
    \\
    \overline{C}_1
    \overline{A}_{2,1}^{\frac{1}{2}}
    \left(
      \begin{array}{@{\hspace{0pt}}c@{\hspace{0pt}}}
        A_0^{\frac{1}{2}}
        A_2^{\frac{1}{2}}
        \overline{A}_1^{\frac{1}{2}}
        \left(
          -
          q^{\frac{\alpha}{2}}
          \ketket{1}{8}
          +
          \overline{q}^{\frac{\alpha}{2}}
          \ketket{8}{1}
        \right)
        +
        A_0
        \overline{A}_1^{\frac{1}{2}}
        \left(
          -
          \overline{q}^{\frac{\alpha}{2}+1}
          \ketket{2}{7}
          +
          q^{\frac{\alpha}{2}+1}
          \ketket{7}{2}
        \right)
        \\
        +
        A_1^{\frac{1}{2}}
        \left(
          -
          q^{\frac{3\alpha}{2}+1}
          \ketket{3}{6}
          +
          \overline{q}^{\frac{3\alpha}{2}+1}
          \ketket{6}{3}
        \right)
        +
        A_1^{\frac{1}{2}}
        \left(
          -
          \overline{q}^{\frac{\alpha}{2}}
          \ketket{4}{5}
          +
          q^{\frac{\alpha}{2}}
          \ketket{5}{4}
        \right)
      \end{array}
    \right)
    \\
    \overline{C}_1
    \overline{A}_{2,3}^{\frac{1}{2}}
    \left(
      \begin{array}{@{\hspace{0pt}}c@{\hspace{0pt}}}
        A_0^{\frac{1}{2}}
        A_2^{\frac{1}{2}}
        \overline{A}_1^{\frac{1}{2}}
        \left(
          +
          q^{\frac{\alpha}{2}}
          \ketket{1}{8}
          -
          \overline{q}^{\frac{\alpha}{2}}
          \ketket{8}{1}
        \right)
        +
        A_2
        \overline{A}_1^{\frac{1}{2}}
        \left(
          -
          q^{\frac{3\alpha}{2}+1}
          \ketket{2}{7}
          +
          \overline{q}^{\frac{3\alpha}{2}+1}
          \ketket{7}{2}
        \right)
        \\
        +
        A_1^{\frac{1}{2}}
        \left(
          -
          \overline{q}^{\frac{\alpha}{2}+1}
          \ketket{3}{6}
          +
          q^{\frac{\alpha}{2}+1}
          \ketket{6}{3}
        \right)
        +
        A_1^{\frac{1}{2}}
        \left(
          +
          \overline{q}^{\frac{\alpha}{2}}
          \ketket{4}{5}
          -
          q^{\frac{\alpha}{2}}
          \ketket{5}{4}
        \right)
      \end{array}
    \right)
  \end{array}
  \right\}.
\end{eqnarray*}
\normalsize

Observe that, again, our \textsc{Mathematica} code has done some
nontrivial work in simplifying the algebraic expressions in $q$ and
$\alpha$. This work would present a significant barrier for a human
calculator.

\pagebreak


\section{Projectors and R matrices for $V\otimes V$}
\label{sec:ProjectorsRMatrices}

\subsection{Projectors $P_k$ onto the $V_k$}
\label{sec:Projectors}

At this stage, for each of the submodules $V_k\subset V\otimes V$, we
have an orthonormal basis
$
  \mathfrak{B}_k
  =
  \left\{
    \ket{\Psi^k_l}\right
  \}_{l=1}^{\mathrm{dim}(V_k)}
$
and, for each of their duals $V_k^*$, a corresponding dual basis
$
  \mathfrak{B}_k^*
  =
  \left\{
    \bra{\Psi^k_l}\right
  \}_{l=1}^{\mathrm{dim}(V_k)}
$,
where $\bra{\Psi^k_l}\triangleq\ket{\Psi^k_l}^\dagger$.
Using these dual bases, it is a simple matter to construct the
projectors $P_k:V\otimes V\to V_k$:
\begin{eqnarray*}
  P_k
  =
  \sum_{j=1}^{\mathrm{dim}(V_k)}
    \ketbra{\Psi^k_l}{\Psi^k_l};
\end{eqnarray*}
note that we must use (\ref{eq:ketketxbrabradefn}) for the
multiplication of tensor products.

We now make a change of notation.  As we did for the matrix elements in
\S\ref{sec:Matrixelements}, we replace $\ketbra{i}{j}$ with the
elementary matrix $e^i_j\in M_{2^{mn}}$. We then use the
notation $e^{ik}_{jl}
\in M_{2^{2mn}}
$ to indicate the two dimensional matrix form of
the usual elementary rank $4$ tensor $e^i_j\otimes e^k_l$, obtained by
inserting a copy of $e^k_l$ at each location of $e^i_j$.  We find that
the $P_k$ are in general quite sparse, that is, only a small fraction
of their $2^{4mn}$ components are nonzero.


\subsubsection{Illustration: $P_0$ for $U_q[gl(3|1)]$}

We illustrate using $P_0$ for the case $U_q[gl(3|1)]$,
which has $125$ (out of $2^{12}=4096$) nonzero components.  We present
these components below, using horizontal lines to separate equivalence
classes of symmetry.
\small
\begin{eqnarray*}
  \hspace{-10pt}
  \begin{array}{l}
    \left\{
      e^{1 1}_{1 1}
    \right\}
    \\[1mm]
    \hline
    \\[-3mm]
    \frac{
      1
    }{
      C_0
    }
    \hspace{-1pt}
    \left\{
    \begin{array}{@{\hspace{-1pt}}l@{\hspace{-1pt}}}
      q^{\alpha}
      \hspace{-2pt}
      \left\{
        e^{1 2}_{1 2},
        e^{1 3}_{1 3},
        e^{1 4}_{1 4}
      \right\}
      \\
      \overline{q}^{\alpha}
      \hspace{-2pt}
      \left\{
        e^{2 1}_{2 1},
        e^{3 1}_{3 1},
        e^{4 1}_{4 1}
      \right\}
      \end{array}
    \right\},
    \;
    \frac{
      A_1
    }{
      C_0
      A_{2,1}
    }
    \hspace{-1pt}
    \left\{
    \begin{array}{@{\hspace{-1pt}}l@{\hspace{-1pt}}}
      q^{2\alpha}
      \hspace{-2pt}
      \left\{
        e^{1 5}_{1 5},
        e^{1 6}_{1 6},
        e^{1 7}_{1 7}
      \right\}
      \\
      \overline{q}^{2\alpha}
      \hspace{-2pt}
      \left\{
        e^{5 1}_{5 1},
        e^{6 1}_{6 1},
        e^{7 1}_{7 1}
      \right\}
    \end{array}
    \right\},
    \;
    \frac{
      A_2
    }{
      C_0
      C_1
      A_{2,1}
    }
    \left\{
      \begin{array}{@{\hspace{-1pt}}l@{\hspace{-1pt}}}
        q^{3\alpha}
        \hspace{-2pt}
        \left\{
          e^{1 8}_{1 8}
        \right\}
        \\
        \overline{q}^{3\alpha}
        \hspace{-2pt}
        \left\{
          e^{8 1}_{8 1}
        \right\}
      \end{array}
    \right\}
    \\
    \frac{
      A_0
    }{
      C_0
      A_{2,1}
    }
    \left\{
      \begin{array}{@{\hspace{0pt}}l@{\hspace{0pt}}}
      \overline{q}
      \left\{
        e^{2 3}_{2 3},
        e^{2 4}_{2 4},
        e^{3 4}_{3 4}
      \right\}
      \\
      q
      \left\{
        e^{3 2}_{3 2},
        e^{4 2}_{4 2},
        e^{4 3}_{4 3}
      \right\}
      \end{array}
    \right\},
    \quad
    \frac{
      A_0
    }{
      C_0
      C_1
      A_{2,1}
    }
    \left\{
      \begin{array}{@{\hspace{0pt}}l@{\hspace{0pt}}}
        q^{\alpha}
        \left\{
        \begin{array}{@{\hspace{0pt}}l@{\hspace{0pt}}}
          \overline{q}^{2}
          e^{2 7}_{2 7},
          e^{3 6}_{3 6},
          q^{2}
          e^{4 5}_{4 5}
        \end{array}
        \right\}
        \\
        \overline{q}^{\alpha}
        \left\{
        \begin{array}{@{\hspace{0pt}}l@{\hspace{0pt}}}
          q^{2}
          e^{7 2}_{7 2},
          e^{6 3}_{6 3},
          \overline{q}^{2}
          e^{5 4}_{5 4}
        \end{array}
      \right\}
      \end{array}
    \right\}
    \\[3mm]
    \hline
    \\[-3mm]
    \frac{
      1
    }{
      C_0
    }
    \left\{
      \begin{array}{@{\hspace{0pt}}l@{\hspace{0pt}}}
        -
        \left\{
          \begin{array}{@{\hspace{0pt}}l@{\hspace{0pt}}}
            e^{1 2}_{2 1},
            e^{1 3}_{3 1},
            e^{1 4}_{4 1}
          \end{array}
        \right\}
        \\
        +
        \left\{
          \begin{array}{@{\hspace{0pt}}l@{\hspace{0pt}}}
            e^{2 1}_{1 2},
            e^{3 1}_{1 3},
            e^{4 1}_{1 4}
          \end{array}
        \right\}
      \end{array}
    \right\},
    \quad
    \frac{
      A_1
    }{
      C_0
      A_{2,1}
    }
    \left\{
      \begin{array}{@{\hspace{0pt}}l@{\hspace{0pt}}}
        e^{1 5}_{5 1},
        e^{1 6}_{6 1},
        e^{1 7}_{7 1}
        \\
        e^{5 1}_{1 5},
        e^{6 1}_{1 6},
        e^{7 1}_{1 7}
      \end{array}
    \right\},
    \quad
    \frac{
      A_2
    }{
      C_0
      C_1
      A_{2,1}
    }
    \left\{
      \begin{array}{@{\hspace{0pt}}l@{\hspace{0pt}}}
        - e^{1 8}_{8 1}
        \\
        + e^{8 1}_{1 8}
      \end{array}
    \right\}
    \\
    -
    \frac{
      A_0
    }{
      C_0
      A_{2,1}
    }
    \left\{
      \begin{array}{@{\hspace{0pt}}l@{\hspace{0pt}}}
        e^{2 3}_{3 2},
        e^{2 4}_{4 2},
        e^{3 4}_{4 3}
        \\
        e^{4 2}_{2 4},
        e^{3 2}_{2 3},
        e^{4 3}_{3 4}
      \end{array}
    \right\},
    \quad
    \frac{
      A_0
    }{
      C_0
      C_1
      A_{2,1}
    }
    \left\{
      \begin{array}{@{\hspace{0pt}}l@{\hspace{0pt}}}
        +
        \left\{
          \begin{array}{@{\hspace{0pt}}l@{\hspace{0pt}}}
            e^{2 7}_{7 2},
            e^{3 6}_{6 3},
            e^{4 5}_{5 4}
          \end{array}
        \right\}
        \\
        -
        \left\{
          \begin{array}{@{\hspace{0pt}}l@{\hspace{0pt}}}
            e^{7 2}_{2 7},
            e^{6 3}_{3 6},
            e^{5 4}_{4 5}
          \end{array}
        \right\}
      \end{array}
    \right\}
    \\[3mm]
    \hline
    \\[-3mm]
    \frac{
      A_0^{\frac{1}{2}}
      A_1^{\frac{1}{2}}
    }{
      C_0
      A_{2,1}
    }
    \left\{
      \begin{array}{@{\hspace{0pt}}rr@{\hspace{0pt}}}
        q^{\alpha-\frac{1}{2}}
        \left\{
          \begin{array}{@{\hspace{0pt}}l@{\hspace{0pt}}}
          -
          \left\{
            \begin{array}{@{\hspace{0pt}}l@{\hspace{0pt}}}
              e^{1 5}_{2 3},
              e^{1 6}_{2 4},
              e^{1 7}_{3 4}
            \end{array}
          \right\}
          \\
          +
          \left\{
            \begin{array}{@{\hspace{0pt}}l@{\hspace{0pt}}}
              e^{2 3}_{1 5},
              e^{2 4}_{1 6},
              e^{3 4}_{1 7}
            \end{array}
          \right\}
          \end{array}
        \right\},
        &
        \overline{q}^{\alpha-\frac{1}{2}}
        \left\{
          \begin{array}{@{\hspace{0pt}}l@{\hspace{0pt}}}
          +
          \left\{
            \begin{array}{@{\hspace{0pt}}l@{\hspace{0pt}}}
              e^{5 1}_{3 2},
              e^{6 1}_{4 2},
              e^{7 1}_{4 3}
            \end{array}
          \right\}
          \\
          -
          \left\{
            \begin{array}{@{\hspace{0pt}}l@{\hspace{0pt}}}
              e^{3 2}_{5 1},
              e^{4 2}_{6 1},
              e^{4 3}_{7 1}
            \end{array}
          \right\}
          \end{array}
        \right\}
        \\
        q^{\alpha+\frac{1}{2}}
        \left\{
          \begin{array}{@{\hspace{0pt}}l@{\hspace{0pt}}}
          +
          \left\{
            \begin{array}{@{\hspace{0pt}}l@{\hspace{0pt}}}
              e^{1 5}_{3 2},
              e^{1 6}_{4 2},
              e^{1 7}_{4 3}
            \end{array}
          \right\}
          \\
          -
          \left\{
            \begin{array}{@{\hspace{0pt}}l@{\hspace{0pt}}}
              e^{3 2}_{1 5},
              e^{4 2}_{1 6},
              e^{4 3}_{1 7}
            \end{array}
          \right\}
          \end{array}
        \right\},
        &
        \overline{q}^{\alpha+\frac{1}{2}}
        \left\{
          \begin{array}{@{\hspace{0pt}}l@{\hspace{0pt}}}
          -
          \left\{
            \begin{array}{@{\hspace{0pt}}l@{\hspace{0pt}}}
              e^{5 1}_{2 3},
              e^{6 1}_{2 4},
              e^{7 1}_{3 4}
            \end{array}
          \right\}
          \\
          +
          \left\{
            \begin{array}{@{\hspace{0pt}}l@{\hspace{0pt}}}
              e^{2 3}_{5 1},
              e^{2 4}_{6 1},
              e^{3 4}_{7 1}
            \end{array}
          \right\}
          \end{array}
        \right\}
      \end{array}
    \right\}
    \\
    \frac{
      A_0^{\frac{1}{2}}
      A_2^{\frac{1}{2}}
    }{
      C_0
      C_1
      A_{2,1}
    }
    \left\{
      \begin{array}{@{\hspace{0pt}}rr@{\hspace{0pt}}}
        q^{2\alpha}
        \left\{
          \begin{array}{@{\hspace{0pt}}l@{\hspace{0pt}}}
            - \overline{q}
            e^{1 8}_{2 7},
            +
            e^{1 8}_{3 6},
            - q
            e^{1 8}_{4 5}
            \\
            + \overline{q}
            e^{2 7}_{1 8},
            -
            e^{3 6}_{1 8},
            + q
            e^{4 5}_{1 8}
          \end{array}
        \right\},
        &
        \overline{q}^{2\alpha}
        \left\{
          \begin{array}{@{\hspace{0pt}}l@{\hspace{0pt}}}
            + q
            e^{8 1}_{7 2},
            -
            e^{8 1}_{6 3},
            + \overline{q}
            e^{8 1}_{5 4}
            \\
            - q
            e^{7 2}_{8 1},
            +
            e^{6 3}_{8 1},
            - \overline{q}
            e^{5 4}_{8 1}
          \end{array}
        \right\}
        \\
        q^{\alpha}
        \left\{
          \begin{array}{@{\hspace{0pt}}l@{\hspace{0pt}}}
            + q
            e^{1 8}_{7 2},
            -
            e^{1 8}_{6 3},
            + \overline{q}
            e^{1 8}_{5 4}
            \\
            + q
            e^{7 2}_{1 8},
            -
            e^{6 3}_{1 8},
            + \overline{q}
            e^{5 4}_{1 8}
          \end{array}
        \right\},
        &
        \overline{q}^{\alpha}
        \left\{
          \begin{array}{@{\hspace{0pt}}l@{\hspace{0pt}}}
            + \overline{q}
            e^{8 1}_{2 7},
            -
            e^{8 1}_{3 6},
            + q
            e^{8 1}_{4 5}
            \\
            + \overline{q}
            e^{2 7}_{8 1},
            -
            e^{3 6}_{8 1},
            + q
            e^{4 5}_{8 1}
          \end{array}
        \right\}
      \end{array}
    \right\}
    \\
    \frac{
      A_0
    }{
      C_0
      C_1
      A_{2,1}
    }
    \left\{
      \begin{array}{@{\hspace{0pt}}rr@{\hspace{0pt}}}
        q^{\alpha}
        \left\{
          \begin{array}{@{\hspace{0pt}}l@{\hspace{0pt}}}
            - \overline{q}
            e^{2 7}_{3 6},
            +
            e^{2 7}_{4 5},
            - q
            e^{3 6}_{4 5}
            \\
            - \overline{q}
            e^{3 6}_{2 7},
            +
            e^{4 5}_{2 7},
            - q
            e^{4 5}_{3 6}
          \end{array}
        \right\},
        &
        \overline{q}^{\alpha}
        \left\{
          \begin{array}{@{\hspace{0pt}}l@{\hspace{0pt}}}
            - q
            e^{7 2}_{6 3},
            +
            e^{7 2}_{5 4},
            - \overline{q}
            e^{6 3}_{5 4}
            \\
            - q
            e^{6 3}_{7 2},
            +
            e^{5 4}_{7 2},
            - \overline{q}
            e^{5 4}_{6 3}
          \end{array}
        \right\}
        \\
        q
        \left\{
          \begin{array}{@{\hspace{0pt}}l@{\hspace{0pt}}}
            +
            e^{7 2}_{3 6},
            - q
            e^{7 2}_{4 5},
            +
            e^{6 3}_{4 5}
            \\
            -
            e^{3 6}_{7 2},
            + q
            e^{4 5}_{7 2},
            -
            e^{4 5}_{6 3}
          \end{array}
        \right\},
        &
        \overline{q}
        \left\{
          \begin{array}{@{\hspace{0pt}}l@{\hspace{0pt}}}
            -
            e^{2 7}_{6 3},
            + \overline{q}
            e^{2 7}_{5 4},
            -
            e^{3 6}_{5 4}
            \\
            +
            e^{6 3}_{2 7},
            - \overline{q}
            e^{5 4}_{2 7},
            +
            e^{5 4}_{3 6}
          \end{array}
        \right\}
      \end{array}
    \right\}.
  \end{array}
\end{eqnarray*}
\normalsize


\subsection{R matrices $\check{R}^{m,n}(u)$ and $\check{R}^{m,n}$}

We may now form the trigonometric R matrix
$\check{R}^{m,n}(u)$ as a weighted sum of the projectors, where the
weights are the eigenvalues of $\check{R}^{m,n}(u)$ on the submodules.
For the \emph{special} case of our representations labeled
$\Lambda=(\dot{0}_m|\dot{\alpha}_n)$, these eigenvalues are actually
\emph{known} \cite{DeliusGouldLinksZhang:95b}.  The quantum R matrix is
then the spectral limit
$\check{R}^{m,n}=\lim_{u\to\infty}\check{R}^{m,n}(u)$.
Again, what follows here strictly applies only to the $m\geqslant n$
case, but given the natural isomorphism between $U_q[gl(m|n)]$ and
$U_q[gl(n|m)]$, this is unimportant.

Again using the notation introduced in \S\ref{sec:ThesubmodulesVk},
especially noting the definition of $\lambda_{\gamma}$ in
(\ref{eq:Defnoflambdagamma}), we have the following expression for
$\check{R}^{m,n} (u)$, normalised such that its `first' component (i.e.
the coefficient of $e^{1 1}_{1 1}$) is unity (for applications, other
normalisations may be applicable, e.g. see \cite{DeWit:99e}):
\begin{eqnarray}
  \check{R}^{m,n} (u)
  =
  \sum_{[\gamma]}
    \Xi_{2\Lambda+\lambda_\gamma} (u)
    P_{2\Lambda+\lambda_\gamma},
  \label{eq:ExpansionofRmnu}
\end{eqnarray}
where, again, as in (\ref{eq:FullUqglmndecomp}),
the sum is over all allowable Young diagrams $[\gamma]$ and
$P_{2\Lambda+\lambda_\gamma}$ is the projector onto the submodule
$V_{2\Lambda+\lambda_\gamma}$.
Recalling that we intend $\gamma=\gamma_1,\gamma_2,\dots,\gamma_r$, the
eigenvalue $\Xi_{2\Lambda+\lambda_\gamma} (u)$ is:
\begin{eqnarray}
  \Xi_{2\Lambda+\lambda_\gamma} (u)
  =
  \prod_{j=1}^{r}
    \prod_{i=1}^{\gamma_j}
      \frac{
        [\alpha+j-i+u]_q
      }{
        [\alpha+j-i-u]_q
      },
  \label{eq:ExplicitXi}
\end{eqnarray}
where, for the empty Young diagram case, we intend
$\Xi_{2\Lambda}(u)=1$. Note that we have substituted $x=q^{-2u}$ from
the original $x$ used in the \emph{multiplicative} Yang--Baxter
equations of \cite{DeliusGouldLinksZhang:95b}: \emph{our} Yang--Baxter
equations are \emph{additive} in variable $u$ (see
\S\ref{sec:YangBaxter}).

Then, as $\check{R}^{m,n}$ is the spectral limit
$\lim_{u\to\infty}\check{R}^{m,n}(u)$, we have:
\begin{eqnarray}
  \check{R}^{m,n}
  =
  \sum_{[\gamma]}
    \xi_{2\Lambda+\lambda_\gamma}
    P_{2\Lambda+\lambda_\gamma},
  \label{eq:ExpansionofRmn}
\end{eqnarray}
where
$
  \xi_{2\Lambda+\lambda_\gamma} (u)
  \triangleq
  \lim_{u\to\infty}
    \Xi_{2\Lambda+\lambda_\gamma} (u)
$.
Again, the coefficient of $e^{1 1}_{1 1}$ is unity.  Evaluating this
limit:
\begin{eqnarray*}
  \lim_{u\to\infty}
    \Xi_{2\Lambda+\lambda_\gamma} (u)
  =
  \prod_{j=1}^{r}
    \prod_{i=1}^{\gamma_j}
      (-)
      q^{2(\alpha+j-i)}
  =
  \prod_{j=1}^{r}
    {(-)}^{\gamma_j}
      q^{\gamma_j(2\alpha+2j-\gamma_j-1)},
\end{eqnarray*}
where we have applied the observation:
$
  \lim_{u\to\infty}
    \frac{
      [X+u]_q
    }{
      [X-u]_q
    }
    =
    -
    q^{+2X}
$.
Strictly, this requires that $|q|>1$, and this is perhaps not so
sensible, as $q$ is in some sense, a \emph{deformation} parameter, so
should be \emph{small}.  If, instead, we assume that $|q|<1$, then the
limit becomes $-q^{-2X}$, and the expressions for
$\xi_{2\Lambda+\lambda_\gamma}$ both above and below (in
(\ref{eq:Explicitxi})), remain valid under the mapping
$q\mapsto\overline{q}$. In the final analysis, this simply means that
we obtain quantum R matrices related by $q\mapsto\overline{q}$.  As
above, we are considering \emph{generic} $q$ (i.e. $q$ is not a root of
unity), and to ensure that, we can demand $|q|\neq 1$.

These considerations aside, we obtain:
\begin{eqnarray}
  \xi_{2\Lambda+\lambda_\gamma}
  =
  {(-)}^{N}
  q^{\sum_{j=1}^{r} \gamma_j(2\alpha+2j-\gamma_j-1)},
  \label{eq:Explicitxi}
\end{eqnarray}
where $N=\sum_{j=1}^{r} \gamma_j$ is the $\mathbb{Z}$ graded level
of $V_{2\Lambda+\lambda_\gamma}$ (cf.
\S\ref{sec:orthogonaldecompositionofVotimesV}).  Note that we intend
$\xi_{2\Lambda}=1$, in agreement with $\Xi_{2\Lambda}(u)=1$.

Thus, we have two methods to compute $\check{R}^{m,n}$.  Firstly, we
may explicitly evaluate it as the spectral limit of
$\check{R}^{m,n}(u)$, itself computed by substituting
(\ref{eq:ExplicitXi}) into (\ref{eq:ExpansionofRmnu}).  Secondly, we
may directly substitute (\ref{eq:Explicitxi}) into
(\ref{eq:ExpansionofRmn}), bypassing the construction of
$\check{R}^{m,n}(u)$ altogether. This method is of course less
computationally expensive, so in practice, we use it, but we have also
implemented the former method, which is useful for checking
consistency.


\subsubsection{Illustration:
               The R matrix decompositions for $U_q[gl(m|1)]$}

Again resurrecting the abusive notion of (\ref{eq:Abuseofnotation}), we
may replace the explicit (\ref{eq:ExpansionofRmnu}) and
(\ref{eq:ExpansionofRmn}) with:
\begin{eqnarray*}
  \check{R}^{m,n} (u)
  =
  \sum_{k}
    \Xi_k(u)
    P_k,
  \qquad \qquad
  \check{R}^{m,n}
  =
  \sum_{k}
    \xi_k
    P_k.
\end{eqnarray*}
For the special case $n=1$, which we are interested in for physical
applications, these simplify to:
\begin{eqnarray*}
  \check{R}^{m,1} (u)
  =
  \sum_{k=0}^{m}
    \Xi_k(u)
    P_k,
  \qquad \qquad
  \check{R}^{m,1}
  =
  \sum_{k=0}^{m}
    \xi_k
    P_k,
\end{eqnarray*}
where $\Xi_k(u)$ and $\xi_k$ are the following eigenvalues
(\cite{GouldLinksZhang:96b}, and cf. \cite{GeGouldZhangZhou:98a}):
\begin{eqnarray*}
  \Xi_k(u)
  =
  \prod_{j=0}^{k-1}
    \frac{
      [\alpha + j + u]_q
    }{
      [\alpha + j - u]_q
    },
  \qquad
  \qquad
  \xi_k
  =
  {(-)}^{k}
  q^{k(2\alpha+k-1)};
\end{eqnarray*}
here we of course again intend $\Xi_0(u)=\xi_0=1$.

\pagebreak


\subsection{Yang--Baxter equations}
\label{sec:YangBaxter}

To be certain, $\check{R}^{m,n} (u)$ satisfies the following
graded version of the (additive) trigonometric Yang--Baxter equation
(TYBE):
\begin{eqnarray}
   \begin{array}{l}
     {(-)}^{[b'][c'] + [a'][c] + [a][b] + [b'][b'']}
     {\check{R} (u)    }^{c'' b''}_{b' c'}
     {\check{R} (u + v)}^{c'  a''}_{a' c }
     {\check{R} (v)    }^{b'  a' }_{a  b }
     \\[1mm]
     \qquad
     \quad
     =
     {(-)}^{[a'][b'] + [a][c'] + [b][c] + [b][b']}
     {\check{R} (v)    }^{b'' a''}_{a' b'}
     {\check{R} (u + v)}^{c'' a' }_{a  c'}
     {\check{R} (u)    }^{c'  b' }_{b  c },
   \end{array}
   \label{eq:ComponentGTYBE}
\end{eqnarray}
where we have written $[a]$ for $[\ket{a}]$.
The parity factors in
(\ref{eq:ComponentGTYBE}) may be removed by the following
transformation (e.g. see \cite{DeWit:98}):
\begin{eqnarray*}
  \begin{array}{ll}
    \check{R}^{a' b'}_{a b} (u)
    \mapsto
    {(-)}^{[a]([b]+[b'])}
    \check{R}^{a' b'}_{a b} (u),
  \end{array}
\end{eqnarray*}
after which $\check{R}(u)$ which satisfied (\ref{eq:ComponentGTYBE})
now satisfies the usual ungraded TYBE:
\begin{eqnarray}
  {\check{R} (u)    }^{c'' b''}_{b' c'}
  {\check{R} (u + v)}^{c'  a''}_{a' c }
  {\check{R} (v)    }^{b'  a' }_{a  b }
  =
  {\check{R} (v)    }^{b'' a''}_{a' b'}
  {\check{R} (u + v)}^{c'' a' }_{a  c'}
  {\check{R} (u)    }^{c'  b' }_{b  c },
  \label{eq:ComponentTYBE}
\end{eqnarray}
written in noncomponent form as:
\begin{eqnarray}
  \check{R}_{12} (u)
  \check{R}_{23} (u + v)
  \check{R}_{12} (v)
  =
  \check{R}_{23} (v)
  \check{R}_{12} (u + v)
  \check{R}_{23} (u).
  \label{eq:NoncomponentTYBE}
\end{eqnarray}
In the spectral limit $\check{R}=\lim_{u\to\infty}\check{R}(u)$, this
of course becomes a quantum Yang--Baxter equation (QYBE):
\begin{eqnarray}
  \check{R}_{12}
  \check{R}_{23}
  \check{R}_{12}
  =
  \check{R}_{23}
  \check{R}_{12}
  \check{R}_{23},
  \label{eq:NoncomponentQYBE}
\end{eqnarray}
viz
$
  (\check{R} \otimes I)
  (I \otimes \check{R})
  (\check{R} \otimes I)
  =
  (I \otimes \check{R})
  (\check{R} \otimes I)
  (I \otimes \check{R})
$,
familiar as the braid relation
$\sigma_1 \sigma_2 \sigma_1 = \sigma_2 \sigma_1 \sigma_2$.

Defining $R(u) \triangleq P \check{R}(u)$, where
$P$ is a permutation operator, yields a trigonometric R matrix
$R(u)$ satisfying the following version of
(\ref{eq:NoncomponentTYBE}):
\begin{eqnarray}
  R_{12} (u)
  R_{13} (u+v)
  R_{23} (v)
  =
  R_{23} (v)
  R_{13} (u+v)
  R_{12} (u).
  \label{eq:NoncomponentAlternativeTYBE}
\end{eqnarray}
This transformation amounts to the mapping:
$
  {R (u)}^{a' b'}_{a b}
  =
  {\check{R} (u)}^{b' a'}_{a b}
$.
In component form, (\ref{eq:NoncomponentAlternativeTYBE}) is
more symmetric than (\ref{eq:ComponentTYBE}):
\begin{eqnarray*}
  {R (u)  }^{b'' c''}_{b' c'}
  {R (u+v)}^{a'' c' }_{a' c }
  {R (v)  }^{a'  b' }_{a  b }
  =
  {R (v)  }^{a'' b''}_{a' b'}
  {R (u+v)}^{a'  c''}_{a  c'}
  {R (u)  }^{b'  c' }_{b  c }.
\end{eqnarray*}

\pagebreak


\subsection{An alternative construction of the quantum R matrix}

In \S\ref{sec:ThesubmodulesVk} and \S\ref{sec:ProjectorsRMatrices}, we
described the construction of trigonometric and quantum R matrices
corresponding to the tensor products of representations of highest
weight $\Lambda$, from explicit knowledge of the decomposition of the
tensor product, and the eigenvalues of the R matrices on the subspaces
of the decomposition.  This method is limited to such situations where
this data is \emph{known}.

An alternative approach to the construction of \emph{quantum} R
matrices sidesteps the construction of bases $\mathfrak{B}_k$ and
projectors $P_k$ altogether, instead using only the knowledge of the
matrix elements. That is, say that we know the `universal' (i.e.
algebraic) form of the R matrix, i.e.:

\begin{eqnarray}
  \mathcal{R}
  =
  \sum_i
    a_i \otimes b_i,
  \label{eq:UniversalRDecomposition}
\end{eqnarray}
where $a_i, b_i\in U_q[gl(m|n)]$, for some hopefully finite sum over
indices $i$. Then, for any particular representation of highest weight
$\Lambda$, we may obtain a quantum R matrix
$\pi_{\Lambda\otimes\Lambda}(\mathcal{R})$ satisfying a parameter-free
version of (\ref{eq:NoncomponentAlternativeTYBE}) from $\mathcal{R}$,
by simply replacing the $a_i$ and $b_i$ with their matrix
representations.  In fact, for $U_q[gl(m|n)]$, in
\cite{KhoroshkinTolstoy:91} we find formulae for $\mathcal{R}$ of the
form (\ref{eq:UniversalRDecomposition}), so this method is feasible.

Implementation of this approach has significant advantages over the
current method; apart from being simpler, and greatly reducing
computational effort, it is considerably more general in that it does
not require knowledge of the tensor product decomposition.

A substantial loss is that is does not yield explicit
\emph{trigonometric} R matrices, so we cannot use this method for
\emph{physical} applications -- recall that our primary intended
application is topological, being the construction of link invariants,
for which we only require \emph{quantum} R matrices.  However, there is
a further loss. An interesting alternative approach to constructing
quantum R matrices involves constructing families of distinct but
`gauge equivalent' quantum R matrices, starting from a single
trigonometric R matrix.  (Details of investigations into this for the
$U_q[gl(2|1)]$ case appear in \cite{LinksDeWit:2000}.) Clearly, we also
lose this alternative approach if we can't construct trigonometric R
matrices.  Furthermore, the method described in the present work
incorporates the foreknowledge of the eigenvalues of the quantum and
trigonometric R matrices. This knowledge is special to our
representations, and it may be be used to assist analysis of the
associated link invariants (again, see \cite{LinksDeWit:2000}). In the
alternative construction, we have no such knowledge in general.

The only barrier to implementing this method is that the results of
\cite{KhoroshkinTolstoy:91} are presented in a somewhat abstract form
which would require considerable modification before being directly
useful in the framework described herein. This approach is outside the
scope of this paper; it is left as a future project.

\pagebreak


\section{Implementation and results}
\label{sec:ImplementationResults}

The entire process has been implemented as a suite of functions in the
interpreted environment of \textsc{Mathematica}.  The procedure to
construct the R matrices requires several stages (viz the algorithms of
\S\ref{sec:UqglmnPBW} to \S\ref{sec:ProjectorsRMatrices}), and the
several thousand lines of \textsc{Mathematica} code are broken down
into functional units to achieve this. The code is available on request
from the author.


\subsection{Data structure for the $U_q[gl(m|n)]$ generators}

A challenging issue in implementation is to find a consistent data
structure to represent the algebra generators.  This problem arises as
the Cartan generators $K_a$ often appear exponentiated as $K_a^N$,
where $N$ is not necessarily a positive integer.  Whilst the
unexponentiated $K_{a}$ (equivalently the ${E^a}_a$) \emph{are} a basis
for the Cartan subalgebra, if we use them as our data structures, then
we must deal with the problem of how to express and manipulate their
exponentials.  Thus it proves pragmatic to regard the more general
`generators' $K_a^N$ as the logical units for computation.  Beyond
this, our data structure must uniformly incorporate the non-Cartan
${E^{a}}_{b}$, which are never exponentiated.  The following
\textsc{Mathematica} pattern integrates these two disparate expressions
into a coherent form:
\begin{center}
  \texttt{Generator[Uqgl[m, n], a, b, N]},
\end{center}
where $\mathtt{m}$ and $\mathtt{n}$ are the fixed $m$ and $n$ defining
the algebra, and, if $\mathtt{a}\neq\mathtt{b}$ (and \texttt{N} is
fixed as $1$), then we intend ${E^{a}}_{b}$, and if
$\mathtt{a}=\mathtt{b}$, then we intend $K_a^N$.

Using this pattern, we are able to implement the collection of PBW
commutators of \S\ref{sec:Commutations} in only a few hundred lines of
code.


\subsection{General comments}

Beyond the data structure, the following aspects of the code are
specifically interesting as expositions of the use of
\textsc{Mathematica}.

\begin{enumerate}
\item
  Implementation of the PBW commutators to normal order strings of
  generators uses the repeated application of \texttt{Rule}s to find a
  fixed point.

\item
  The code establishes and solves a system of linear equations to
  determine the parameters defining the highest weight vectors
  $\ket{\Psi^k_1}$ of the $V_k$, where even the size of this system is
  not specified in advance (see \S\ref{sec:Ahighestweightvector}).
  Doing this manually is particularly finicky due to the semantic
  complexity of the expressions involved. (In the literature, this
  process is specifically \emph{avoided} whenever possible; e.g.
  \cite{GeGouldZhangZhou:98a} contains two different kludges.)

\item
  Throughout, it has proven possible to maintain tensor product vectors
  and rank $4$ tensor components in $q$ graded symmetric combinations,
  which facilitates both the simplification of expressions and the
  presentation of the output. This has been achieved by applying
  rewrite rules to nonsymmetric expressions in a carefully controlled
  manner.

\end{enumerate}


\subsection{Limitations}

The computations are computationally inefficient, and this is due
mostly to the fact that the algorithms used are direct, and not
refined.  Although there are no theoretical limits to $m$ and $n$,
computer storage and human patience mean that a current reasonable
practical limit is $mn\leqslant 4$. Whilst the technically difficult
translation of the interpreted \textsc{Mathematica} code into a
compiled language would increase the speed of the computations
enormously, storage requirements would still limit $mn$ to perhaps
$7$ in the general case.


\subsection{Summary of results}
\label{sec:Results}

Both $\check{R}^{m,1} (u)$ and $\check{R}^{m,1}$ have been obtained for
$m=1,2,3,4$.  Fixing $m$, these are rank $4$ tensors, where the tensor
indices run from $1$ to $2^{m}$ (i.e. the dimension of the underlying
representation).  Thus, they contain $D_m \triangleq 2^{4m}$ (albeit
mostly zero) components.  Where $N'_m$ and $N_m$ denote the number of
nonzero components of $\check{R}^{m,1} (u)$ and $\check{R}^{m,1}$
respectively; as $\check{R}^{m,1}$ is the spectral limit of
$\check{R}^{m,1} (u)$, we find, as expected, that $N_m \leqslant N'_m$.
We also find that $N'_m = 6^m$ (why?).  Where $s_m \triangleq N_m/D_m$
denotes the sparsity of $\check{R}^{m,1}$, we find that $s_m$ rapidly
decreases with increasing $m$. Table \ref{tab:data} presents this data
for $m=1,2,3,4$.

\begin{table}[ht]
  \begin{centering}
  \begin{tabular}{crrrr}
    $m$ & $  D_m$ & $N'_m$ & $N_m$ & $s_m (\%)$ \\[1mm]
    \hline
        &         &        &       &            \\[-3mm]
    $1$ & $   16$ & $   6$ & $  5$ & $31.3$     \\
    $2$ & $  256$ & $  36$ & $ 26$ & $10.2$     \\
    $3$ & $ 4096$ & $ 216$ & $139$ & $ 3.4$     \\
    $4$ & $65536$ & $1296$ & $758$ & $ 1.1$
  \end{tabular}
  \caption{
    Data for R matrix construction.
  }
  \label{tab:data}
  \end{centering}
\end{table}

For good measure, we also record the numbers of components of each of
the $m+1$ projectors for the above cases, finding them similar to
$N_m'$ and $N_m$. This data is included in Table
\ref{tab:dimsubprojsize}, which also records submodule dimensions.

\begin{table}[ht]
  \begin{centering}
  \begin{tabular}{ccc}
    $m$ & submodule dimensions & projector sizes \\[1mm]
    \hline
        &                  &                          \\[-3mm]
    $1$ & $     2,2      $ & $          5,5         $ \\
    $2$ & $    4,8,4     $ & $       25,34,25       $ \\
    $3$ & $  8,24,24,8   $ & $   125,199,199,125    $ \\
    $4$ & $16,64,96,64,16$ & $625,1124,1254,1124,625$
  \end{tabular}
  \caption{
    Dimensions of the submodules in the tensor product, and the number
    of nonzero components of the corresponding projectors.
  }
  \label{tab:dimsubprojsize}
  \end{centering}
\end{table}

\pagebreak

Computer run times involved are listed in Table \ref{tab:timingdata}.
These show a rapid increase in cost with increasing $m$. This is
accompanied by an exorbitant increase in storage required.

\begin{table}[ht]
  \begin{centering}
  \begin{tabular}{lrrrr}
                         & \multicolumn{4}{c}{$m$} \\
                         & $   1$ & $   2$ & $    3$ & $   4$ \\[1mm]
    \hline
                         &        &        &         &        \\[-3mm]
    Representations      & $0.10$ & $ 0.5$ & $  2.8$ & $  20$ \\
    TP submodule bases   & $0.74$ & $ 7.4$ & $ 72.6$ & $1559$ \\
    TP projectors        & $0.12$ & $ 4.0$ & $ 55.4$ & $ 877$ \\
    $\check{R}^{m,1} (u)$& $0.27$ & $ 6.5$ & $ 89.5$ & $1736$ \\
    $\check{R}^{m,1}$    & $0.04$ & $ 1.4$ & $ 18.3$ & $ 570$
  \end{tabular}
  \caption{%
    Run times for various calculations (in CPU seconds), using
    \textsc{Mathematica 3} on a \textsc{Sun ULTRA 60} computer. The
    timing for construction of the $5$ submodule bases for
    $U_q[gl(4|1)]$ is the sum of the timings for the individual
    modules, viz:  $1559 = 27 + 201 + 1061 + 218 + 53$.
  }
  \label{tab:timingdata}
  \end{centering}
\end{table}

Of these, the $U_q[gl(1|1)]$ case can be done by hand in
a few hours \cite{BrackenGouldZhangDelius:94b}; the complete
$U_q[gl(2|1)]$ case appears in my PhD thesis \cite{DeWit:98}, and took
several weeks to do by hand; partial details (i.e. up to calculation of
$3$ out of $4$ of the $\overline{\mathfrak{B}}_k$) for the
$U_q[gl(3|1)]$ case appear in \cite{GeGouldZhangZhou:98a}; whilst the
$U_q[gl(4|1)]$ case is entirely new.

Of interest is that $\check{R}^{m,1}$ generally contains non-binomial
irreducible polynomial factors as well as various $q$ brackets.

By direct substitution, we have been able to verify that each
$\check{R}^{m,1} (u)$ satisfies (\ref{eq:NoncomponentTYBE})%
\footnote{%
  We didn't check $\check{R}^{4,1} (u)$, as the computations
  would have been excessively expensive.
}
and that each $\check{R}^{m,1}$ satisfies (\ref{eq:NoncomponentQYBE}).


\section*{Acknowledgements}

My research at Kyoto University in 1999 and 2000 was funded by a
Postdoctoral Fellowship for Foreign Researchers (\# P99703), provided
by the Japan Society for the Promotion of Science.
D\={o}mo arigat\={o} gozaimashita!

Some of this work was completed in March 1999, under the direction of
Mark Gould as part of ongoing research at The University of Queensland,
Australia.  I further thank Jon Links of the same institute for
continuing helpful discussions and general bonhomie.  I also wish to
thank Hiroshi Yamada of Kitami Institute of Technology, Hokkaido, for
hospitality during August 2000, during which some of the writing of
this work was completed.


\pagebreak

\appendix

\section{Proofs of various lemmas}
\label{app:Lemmas}

\subsection{A commutation lemma -- the proof of
            (\protect \ref{eq:FullCNCCommutator})}
\label{app:CommutatorLemma}

\begin{lemma}
  \begin{eqnarray*}
    K_{a}^N
    {E^b}_{c}
    =
    q_a^{N (\delta^a_b-\delta^a_{c})}
    {E^b}_{c}
    K_{a}^N,
  \end{eqnarray*}
  for any meaningful indices $b, c$, and any power $N$.
  \label{lem:CommutatorLemma}
\end{lemma}

\vspace{\baselineskip}
\noindent\textbf{Proof:}

Firstly, we show the following result:
\begin{eqnarray}
  K_{a}
  {E^b}_{c}
  =
  q_a^{(\delta^a_b-\delta^a_{c})}
  {E^b}_{c}
  K_{a},
  \label{eq:HalfwayCommutator}
\end{eqnarray}
for \emph{any} meaningful indices $b\neq c$, not just for the simple
generators where we have $|b-c|=1$.

To see this, we first consider the case $b>c$, so that ${E^b}_c$ is
a lowering generator.  We use induction on $b-c$ to show the result,
assuming (\ref{eq:HalfwayCommutator})
for some $b>c$, where we know that it is true for $c=b-1$ by
(\ref{eq:Donkey}). Then:
\begin{eqnarray*}
  \hspace{-15pt}
  \begin{array}{rcl}
    K_{a}
    {E^b}_{c-1}
    & \stackrel{(\ref{eq:UqglmnNonSimpleGenerators}a)}{=} &
    K_{a}
    \left(
      {E^b}_{c}
      {E^c}_{c-1}
      -
      q_c
      {E^c}_{c-1}
      {E^b}_{c}
    \right)
    \\
    & = &
    \left(
      K_{a}
      {E^b}_{c}
    \right)
    {E^c}_{c-1}
    -
    q_c
    \left(
      K_{a}
      {E^c}_{c-1}
    \right)
    {E^b}_{c}
    \\
    & \stackrel{(\ref{eq:HalfwayCommutator})}{=} &
    q_a^{(\delta^a_b-\delta^a_{c})}
    {E^b}_{c}
    \left(
      K_{a}
      {E^c}_{c-1}
    \right)
    -
    q_c
    q_a^{(\delta^a_c-\delta^a_{c-1})}
    \left(
      {E^c}_{c-1}
      K_{a}
    \right)
    {E^b}_{c}
    \\
    & \stackrel{(\ref{eq:Donkey})}{=} &
    q_a^{(\delta^a_b-\delta^a_{c})}
    {E^b}_{c}
    q_a^{(\delta^a_c-\delta^a_{c-1})}
    {E^c}_{c-1}
    K_{a}
    -
    \\
    & &
    \qquad
    q_c
    q_a^{(\delta^a_c-\delta^a_{c-1})}
    {E^c}_{c-1}
    q_a^{(\delta^a_b-\delta^a_{c})}
    {E^b}_{c}
    K_{a}
    \\
    & = &
    q_a^{(\delta^a_b-\delta^a_{c}+\delta^a_c-\delta^a_{c-1})}
    \left(
      {E^b}_{c}
      {E^c}_{c-1}
      -
      q_c
      {E^c}_{c-1}
      {E^b}_{c}
    \right)
    K_{a}
    \\
    & \stackrel{(\ref{eq:UqglmnNonSimpleGenerators}a)}{=} &
    q_a^{(\delta^a_b-\delta^a_{c-1})}
    {E^b}_{c-1}
    K_{a}.
  \end{array}
\end{eqnarray*}
Thus, (\ref{eq:HalfwayCommutator}) is true for ${E^b}_{c-1}$ if it is
true for ${E^b}_c$, and as it is true for $c=b-1$, thus it is true for
all $b>c$.  The proof for raising generators follows by a trivial
analogy.

Next, observe that setting $b=c$ in (\ref{eq:HalfwayCommutator}) also
yields a true statement. That is, (\ref{eq:HalfwayCommutator}) then
states $K_a {E^b}_b={E^b}_b K_a$, which is equivalent to
(\ref{eq:CartanGensCommute}) when we replace $K_{a}$ with
$q^{{(-)}^{[a]} {E^a}_a}$, and expand the exponential as a power
series.

Finally, replacement of $q$ with $q^N$ in (\ref{eq:HalfwayCommutator})
and again using $K_{a}= q^{{(-)}^{[a]}{E^{a}}_{a}}$ allows us
to write the more general statement in (\ref{eq:FullCNCCommutator}).
\hfill$\Box$

\vspace{3mm}

Having proven this, in retrospect, it appears that some use of the
`generalised Lusztig automorphisms' \cite{DeWit:2000,Zhang:93} might
facilitate a more elegant proof of the recursion in
(\ref{eq:HalfwayCommutator}).


\subsection{A coproduct lemma -- the proof of
            (\protect \ref{eq:FullCoproduct})}
\label{app:CoproductLemma}


\begin{lemma}
  \begin{eqnarray*}
    \hspace{-20pt}
    \Delta({E^a}_b)
    & = &
    {E^a}_b
    \otimes
    K_{a}^{\frac{1}{2}S^a_b}
    \overline{K}_b^{\frac{1}{2}S^a_b}
    +
    \overline{K}_a^{\frac{1}{2}S^a_b}
    K_{b}^{\frac{1}{2}S^a_b}
    \otimes
    {E^a}_b
    -
    \\
    & &
    \qquad
    S^a_b
    \sum_{c}
      \Delta_c
      \left(
        \overline{K}_c^{\frac{1}{2}S^a_b}
        K_{b}^{\frac{1}{2}S^a_b}
        {E^a}_c
        \otimes
        {E^c}_b
        K_{a}^{\frac{1}{2}S^a_b}
        \overline{K}_c^{\frac{1}{2}S^a_b}
      \right),
  \end{eqnarray*}
  for \emph{any} valid indices $a,b$ (including $a=b$), where the sum
  ranges over all $c$ strictly between $a$ and $b$, and is simply
  ignored if $|a-b|\leqslant 1$.
  \label{lem:FullCoproduct}
\end{lemma}

\vspace{\baselineskip}
\noindent\textbf{Proof:}

Firstly, note that, for simple generators i.e.  $|a-b|\leqslant 1$,
the result is just the coproduct of (\ref{eq:SimpleFullCoproduct}):
\begin{eqnarray*}
  \Delta({E^a}_b)
  =
  {E^a}_b
  \otimes
  K_{a}^{\frac{1}{2}S^a_b}
  \overline{K}_b^{\frac{1}{2}S^a_b}
  +
  \overline{K}_{a}^{\frac{1}{2}S^a_b}
  K_{b}^{\frac{1}{2}S^a_b}
  \otimes
  {E^a}_b,
\end{eqnarray*}
as the sum is ignored.  Specifically, it applies to the
Cartan generators $a=b$, being
$
  \Delta({E^{a}}_{a})
  =
  {E^{a}}_{a} \otimes \mathrm{Id}
  +
  \mathrm{Id} \otimes {E^{a}}_{a}
$,
which is equivalent to (\ref{eq:UqglmnCoproduct}c) once we make
the identification $K_a = q^{{(-)}^{[a]} {E^{a}}_{a}}$.

More generally, given (\ref{eq:SimpleFullCoproduct}), and the expansion
of the nonsimple generators (\ref{eq:UqglmnNonSimpleGenerators}), the
following straightforward induction on $|a-b|$ shows that the result
follows for arbitrary \emph{non}simple generators ${E^a}_b$.

We first deal with the lowering case $a>b$, viz $S^a_b=+1$.
For the inductive step, we assume the
truth of our statement for $a-b=p$, for some $p>1$, and we show that
this implies its truth for $a-b=p+1$. Thus, beginning with:
\begin{eqnarray*}
  \Delta
  \left(
    {E^a}_{b-1}
  \right)
  & \stackrel{(\ref{eq:UqglmnNonSimpleGenerators}a)}{=} &
  \Delta
  \left(
    {E^a}_{b}
    {E^b}_{b-1}
    -
    q_b
    {E^b}_{b-1}
    {E^a}_{b}
  \right)
  \\
  & \stackrel{(\ref{eq:Deltaisahomo})}{=} &
  \Delta \left( {E^a}_{b} \right)
  \Delta \left( {E^b}_{b-1} \right)
  -
  q_b
  \Delta \left( {E^b}_{b-1} \right)
  \Delta \left( {E^a}_{b} \right),
\end{eqnarray*}
and expanding $\Delta\left({E^a}_{b}\right)$ using our
inductive hypothesis and
$\Delta\left({E^b}_{b-1}\right)$ from the definition
(\ref{eq:UqglmnCoproduct}a), we obtain:
\begin{eqnarray*}
  & &
  \hspace{-40pt}
  \Delta
  \left(
    {E^a}_{b-1}
  \right)
  =
  \\
  & &
  \hspace{-35pt}
  \left\{
  \begin{array}{ll}
    +
    {E^{a}}_{b}
    {E^{b}}_{b-1}
    \otimes
    K_{a}^{\frac{1}{2}}
    \overline{K}_{b}^{\frac{1}{2}}
    K_{b}^{\frac{1}{2}}
    \overline{K}_{b-1}^{\frac{1}{2}}
    &
    -
    q_b
    {E^{b}}_{b-1}
    {E^{a}}_{b}
    \otimes
    K_{b}^{\frac{1}{2}}
    \overline{K}_{b-1}^{\frac{1}{2}}
    K_{a}^{\frac{1}{2}}
    \overline{K}_{b-1}^{\frac{1}{2}}
    \\
    +
    \overline{K}_{a}^{\frac{1}{2}}
    K_{b}^{\frac{1}{2}}
    \overline{K}_{b}^{\frac{1}{2}}
    K_{b-1}^{\frac{1}{2}}
    \otimes
    {E^{a}}_{b}
    {E^{b}}_{b-1}
    &
    -
    q_b
    \overline{K}_{b}^{\frac{1}{2}}
    K_{b-1}^{\frac{1}{2}}
    \overline{K}_{a}^{\frac{1}{2}}
    K_{b}^{\frac{1}{2}}
    \otimes
    {E^{b}}_{b-1}
    {E^{a}}_{b}
    \\
    +
    \overline{K}_{a}^{\frac{1}{2}}
    K_{b}^{\frac{1}{2}}
    {E^{b}}_{b-1}
    \otimes
    {E^{a}}_{b}
    K_{b}^{\frac{1}{2}}
    \overline{K}_{b-1}^{\frac{1}{2}}
    &
    -
    q_b
    {E^{b}}_{b-1}
    \overline{K}_{a}^{\frac{1}{2}}
    K_{b}^{\frac{1}{2}}
    \otimes
    K_{b}^{\frac{1}{2}}
    \overline{K}_{b-1}^{\frac{1}{2}}
    {E^{a}}_{b}
    \\
    +
    {E^{a}}_{b}
    \overline{K}_{b}^{\frac{1}{2}}
    K_{b-1}^{\frac{1}{2}}
    \otimes
    K_{a}^{\frac{1}{2}}
    \overline{K}_{b}^{\frac{1}{2}}
    {E^{b}}_{b-1}
    &
    -
    q_b
    \overline{K}_{b}^{\frac{1}{2}}
    K_{b-1}^{\frac{1}{2}}
    {E^{a}}_{b}
    \otimes
    {E^{b}}_{b-1}
    K_{a}^{\frac{1}{2}}
    \overline{K}_{b}^{\frac{1}{2}}
  \end{array}
  \right\}
  -
  \\
  & &
  \hspace{-30pt}
  \sum_{c=b+1}^{a-1}
    (q_c - \overline{q}_c)
    \left\{
      \begin{array}{ll}
        +
        \overline{K}_{c}^{\frac{1}{2}}
        K_{b}^{\frac{1}{2}}
        {E^{a}}_{c}
        {E^{b}}_{b-1}
        \otimes
        {E^{c}}_{b}
        K_{a}^{\frac{1}{2}}
        \overline{K}_{c}^{\frac{1}{2}}
        K_{b}^{\frac{1}{2}}
        \overline{K}_{b-1}^{\frac{1}{2}}
        \\
        \quad
        -
        q_b
        {E^{b}}_{b-1}
        \overline{K}_{c}^{\frac{1}{2}}
        K_{b}^{\frac{1}{2}}
        {E^{a}}_{c}
        \otimes
        K_{b}^{\frac{1}{2}}
        \overline{K}_{b-1}^{\frac{1}{2}}
        {E^{c}}_{b}
        K_{a}^{\frac{1}{2}}
        \overline{K}_{c}^{\frac{1}{2}}
        \\
        +
        \overline{K}_{c}^{\frac{1}{2}}
        K_{b}^{\frac{1}{2}}
        {E^{a}}_{c}
        \overline{K}_{b}^{\frac{1}{2}}
        K_{b-1}^{\frac{1}{2}}
        \otimes
        {E^{c}}_{b}
        K_{a}^{\frac{1}{2}}
        \overline{K}_{c}^{\frac{1}{2}}
        {E^{b}}_{b-1}
        \\
        \quad
        -
        q_b
        \overline{K}_{b}^{\frac{1}{2}}
        K_{b-1}^{\frac{1}{2}}
        \overline{K}_{c}^{\frac{1}{2}}
        K_{b}^{\frac{1}{2}}
        {E^{a}}_{c}
        \otimes
        {E^{b}}_{b-1}
        {E^{c}}_{b}
        K_{a}^{\frac{1}{2}}
        \overline{K}_{c}^{\frac{1}{2}}
      \end{array}
    \right\}
    \\
    & & \triangleq
    \left( t_1 + t_2 + t_3 + t_4 \right)
    -
    \sum_{c=b+1}^{a-1}
      \Delta_c
      (t_5 + t_6).
\end{eqnarray*}
Note that in the above, no parity factors appear, as ${E^{a}}_{b}$ and
${E^{b}}_{b-1}$ cannot both simultaneously be odd.
Next, we examine the terms $t_i$:
\begin{eqnarray*}
  t_1
  & = &
  \left(
    {E^{a}}_{b}
    {E^{b}}_{b-1}
    -
    q_b
    {E^{b}}_{b-1}
    {E^{a}}_{b}
  \right)
  \otimes
  K_{a}^{\frac{1}{2}}
  \overline{K}_{b-1}^{\frac{1}{2}}
  \stackrel{(\ref{eq:UqglmnNonSimpleGenerators}a)}{=}
  {E^{a}}_{b-1}
  \otimes
  K_{a}^{\frac{1}{2}}
  \overline{K}_{b-1}^{\frac{1}{2}}
  \\
  t_2
  & = &
  \overline{K}_{a}^{\frac{1}{2}}
  K_{b-1}^{\frac{1}{2}}
  \otimes
  \left(
    {E^{a}}_{b}
    {E^{b}}_{b-1}
    -
    q_b
    {E^{b}}_{b-1}
    {E^{a}}_{b}
  \right)
  \stackrel{(\ref{eq:UqglmnNonSimpleGenerators}a)}{=}
  \overline{K}_{a}^{\frac{1}{2}}
  K_{b-1}^{\frac{1}{2}}
  \otimes
  {E^{a}}_{b-1}
  \\
  t_3
  & \stackrel{(\ref{eq:FullCNCCommutator})}{=} &
  \left(
    \overline{K}_{a}^{\frac{1}{2}}
    K_{b}^{\frac{1}{2}}
    {E^{b}}_{b-1}
    \otimes
    {E^{a}}_{b}
    K_{b}^{\frac{1}{2}}
    \overline{K}_{b-1}^{\frac{1}{2}}
  \right)
  \left(
    1
    -
    q_b
    \overline{q}_b^{\frac{1}{2}}
    \overline{q}_b^{\frac{1}{2}}
  \right)
  =
  0
  \\
  t_4
  & \stackrel{(\ref{eq:FullCNCCommutator})}{=} &
  \left(
    \overline{K}_{b}^{\frac{1}{2}}
    K_{b-1}^{\frac{1}{2}}
    {E^{a}}_{b}
    \otimes
    {E^{b}}_{b-1}
    K_{a}^{\frac{1}{2}}
    \overline{K}_{b}^{\frac{1}{2}}
  \right)
  \left(
    \overline{q}_b^{\frac{1}{2}}
    \overline{q}_b^{\frac{1}{2}}
    -
    q_b
  \right)
  \\
  & = &
  -
  \Delta_b
  \left(
    \overline{K}_{b}^{\frac{1}{2}}
    K_{b-1}^{\frac{1}{2}}
    {E^{a}}_{b}
    \otimes
    {E^{b}}_{b-1}
    K_{a}^{\frac{1}{2}}
    \overline{K}_{b}^{\frac{1}{2}}
  \right)
  \\
  t_5
  & \stackrel{(\ref{eq:FullCNCCommutator})}{=} &
  \left(
    \overline{K}_{c}^{\frac{1}{2}}
    K_{b}^{\frac{1}{2}}
    {E^{a}}_{c}
    {E^{b}}_{b-1}
    \otimes
    {E^{c}}_{b}
    K_{a}^{\frac{1}{2}}
    \overline{K}_{c}^{\frac{1}{2}}
    K_{b}^{\frac{1}{2}}
    \overline{K}_{b-1}^{\frac{1}{2}}
  \right)
  \left(
    1
    -
    q_b
    \overline{q}_b^{\frac{1}{2}}
    \overline{q}_b^{\frac{1}{2}}
  \right)
  =
  0
  \\
  t_6
  & \stackrel{(\ref{eq:FullCNCCommutator})}{=} &
  \overline{K}_{c}^{\frac{1}{2}}
  K_{b-1}^{\frac{1}{2}}
  {E^{a}}_{c}
  \otimes
  \left(
    {E^{c}}_{b}
    {E^{b}}_{b-1}
    -
    q_b
    {E^{b}}_{b-1}
    {E^{c}}_{b}
  \right)
  K_{a}^{\frac{1}{2}}
  \overline{K}_{c}^{\frac{1}{2}}
  \\
  & \stackrel{(\ref{eq:UqglmnNonSimpleGenerators}a)}{=} &
  \overline{K}_{c}^{\frac{1}{2}}
  K_{b-1}^{\frac{1}{2}}
  {E^{a}}_{c}
  \otimes
  {E^{c}}_{b-1}
  K_{a}^{\frac{1}{2}}
  \overline{K}_{c}^{\frac{1}{2}}.
\end{eqnarray*}
Thus, we have:
\begin{eqnarray*}
  \Delta({E^a}_{b-1})
  & = &
  {E^a}_{b-1}
  \otimes
  K_{a}^{\frac{1}{2}} \overline{K}_{b-1}^{\frac{1}{2}}
  +
  \overline{K}_a^{\frac{1}{2}} K_{b-1}^{\frac{1}{2}}
  \otimes
  {E^a}_{b-1}
  -
  \\
  & &
  \quad
    \Delta_b
    \left(
      \overline{K}_b^{\frac{1}{2}} K_{b-1}^{\frac{1}{2}} {E^a}_b
      \otimes
      {E^b}_{b-1} K_{a}^{\frac{1}{2}} \overline{K}_b^{\frac{1}{2}}
    \right)
  -
  \\
  & &
  \quad \quad
  {\displaystyle \sum_{c=b+1}^{a-1}}
    \Delta_c
    \left(
      \overline{K}_c^{\frac{1}{2}} K_{b-1}^{\frac{1}{2}} {E^a}_c
      \otimes
      {E^c}_{b-1} K_{a}^{\frac{1}{2}} \overline{K}_c^{\frac{1}{2}}
    \right)
  \\
  & = &
  {E^a}_{b-1}
  \otimes
  K_{a}^{\frac{1}{2}} \overline{K}_{b-1}^{\frac{1}{2}}
  +
  \overline{K}_a^{\frac{1}{2}} K_{b-1}^{\frac{1}{2}}
  \otimes
  {E^a}_{b-1}
  -
  \\
  & &
  \quad
  {\displaystyle \sum_{c=(b-1)+1}^{a-1}}
    \Delta_c
    \left(
      \overline{K}_c^{\frac{1}{2}}
      K_{b-1}^{\frac{1}{2}}
      {E^a}_c
      \otimes
      {E^c}_{b-1}
      K_{a}^{\frac{1}{2}}
      \overline{K}_c^{\frac{1}{2}}
    \right)
\end{eqnarray*}
Thus, the result is true for $a-b=p+1$ if it is true for $a-b=p$, hence
it is true for all lowering generators ${E^a}_b$, where $a>b$.

Now observe that the definition of $\Delta$ for simple raising
generators is obtained from that for simple lowering generators by the
mapping $q\mapsto\overline{q}$. Also, the definition of the nonsimple
raising generators is obtained from that of the nonsimple lowering
generators under the same mapping. Together, these definitions imply
that the expression for $\Delta$ for the nonsimple raising generators
may be obtained from that for $\Delta$ of the nonsimple lowering
generators under that mapping $q\mapsto\overline{q}$.

\hfill$\Box$

It is also possible to deduce the coproducts of nonsimple generators
via an entirely different approach, using the operator $L(x)$ defined
in \cite{Zhang:92}, but we have not followed this up.


\pagebreak


\section{Explicit results for the $U_q[gl(3|1)]$ case}
\label{app:Data}

\noindent
Here, $[i]=0$ for $i\in\{1;5,6,7\}$ and $[i]=1$ for $i\in\{2,3,4;8\}$.
Apart from all the notational conventions mentioned in the main text
(recall $\Delta\triangleq q - \overline{q}$ of
\S\ref{sec:MatrixelementsE34} and the $A_{i,j}^z$ and $C_{i,j}$ of
(\ref{eq:DefinitionofAijz}) and (\ref{eq:DefinitionofCij}) in
\S\ref{sec:Ahighestweightvector}), we add a couple more to condense the
results.

\begin{itemize}
\item
  We use the following notation as a shorthand for the $q$ graded
  symmetric combination of tensor product vectors:
  \begin{eqnarray*}
    q^{x}_{\pm} \ketket{i}{j}
    & \triangleq &
    q^{x}
    \ketket{i}{j}
    \pm
    \overline{q}^{x}
    \ketket{j}{i}.
  \end{eqnarray*}

  This notation leads to such eyesores as
  ``$q^{0}_{\pm}\ketket{i}{j}$''.
  In these expressions, we shall choose $i \leqslant j$:
  observe that if $j>i$, we may
  replace $q^{x}_{\pm}\ketket{i}{j}$ with
  $
    \pm q^{-x}_{\pm}\ketket{j}{i}
    \equiv
    \pm\overline{q}^{x}_{\pm}\ketket{j}{i}
  $.

\item
  To convert the \emph{graded} R matrices into the equivalent ungraded
  objects, simply multiply all terms in \textbf{boldface} by $-1$.

  Having done this conversion the following notation is a convenient
  shorthand for the $q$ graded symmetric combination of rank $4$
  tensors in the resultant ungraded R matrices:
  \begin{eqnarray*}
    q^{x}_{\pm} e^{i k}_{j l}
    & \triangleq &
    q^{x} e^{i k}_{j l} \pm \overline{q}^{x} e^{k i}_{l j}.
  \end{eqnarray*}

\end{itemize}

\pagebreak

\input{data/rep3}

\input{data/basis30}

\input{data/basis31}

\input{data/basis32}

\input{data/basis33}

\input{data/trigrmx3}

\pagebreak

\input{data/qrmx3}

\pagebreak


\bibliographystyle{plain}
\bibliography{DeWit99c}

\end{document}

%% file: data/rep3.tex
\subsection{Matrix elements of the $U_q[gl(3|1)]$ generators}
\label{app:MatrixElements}

For completeness, here we present matrix elements
$\pi_{(0,0,0\,|\,\alpha)}(X)$ for \emph{all} the $U_q[gl(3|1)]$
generators $X$ (including even the ${E^a}_a$, for comparison with the
$K_a$).

\begin{eqnarray*}
  \hspace{-30pt}
  \begin{array}{r@{\hspace{3pt}}c@{\hspace{3pt}}l} 
    \pi ( {E^1}_1 )
    & = &
    - e^4_4
    - e^6_6
    - e^7_7
    - e^8_8
    \\
    \pi ( {E^2}_2 )
    & = &
    - e^3_3
    - e^5_5
    - e^7_7
    - e^8_8
    \\
    \pi ( {E^3}_3 )
    & = &
    - e^2_2
    - e^5_5
    - e^6_6
    - e^8_8
    \\
    \pi ( K_1 )
    & = &
                    e^1_1
    +               e^2_2
    +               e^3_3
    + \overline{q}  e^4_4
    +               e^5_5
    + \overline{q}  e^6_6
    + \overline{q}  e^7_7
    + \overline{q}  e^8_8
    \\
    \pi ( K_2 )
    & = &
                    e^1_1
    +               e^2_2
    + \overline{q}  e^3_3
    +               e^4_4
    + \overline{q}  e^5_5
    +               e^6_6
    + \overline{q}  e^7_7
    + \overline{q}  e^8_8
    \\
    \pi ( K_3 )
    & = &
                    e^1_1
    + \overline{q}  e^2_2
    +               e^3_3
    +               e^4_4
    + \overline{q}  e^5_5
    + \overline{q}  e^6_6
    +               e^7_7
    + \overline{q}  e^8_8
    \\
    \pi ( {E^4}_4 )
    & = &
       \alpha        e^1_1
    + (\alpha + 1)  (e^2_2 + e^3_3 + e^4_4)
    + (\alpha + 2)  (e^5_5 + e^6_6 + e^7_7)
    + (\alpha + 3)   e^8_8
    \\
    \pi ( K_4 )
    & = &
      \overline{q}^{\alpha}      e^1_1
    + \overline{q}^{\alpha + 1}  (e^2_2 + e^3_3 + e^4_4)
    + \overline{q}^{\alpha + 2}  (e^5_5 + e^6_6 + e^7_7)
    + \overline{q}^{\alpha + 3}  e^8_8
    \\
    \pi ( {E^1}_2 )
    & = &
    - e^3_4
    - e^5_6
    \\ 
    \pi ( {E^2}_1 )
    & = &
    - e^4_3
    - e^6_5
    \\
    \pi ( {E^2}_3 )
    & = &
    - e^2_3
    - e^6_7
    \\
    \pi ( {E^3}_2 )
    & = &
    - e^3_2
    - e^7_6
    \\ 
    \pi ( {E^3}_4 )
    & = &
      A_{0}^{\frac{1}{2}} e^1_2
    + A_{1}^{\frac{1}{2}} e^3_5
    + A_{1}^{\frac{1}{2}} e^4_6
    + A_{2}^{\frac{1}{2}} e^7_8
    \\ 
    \pi ( {E^4}_3 )
    & = &
      A_{0}^{\frac{1}{2}} e^2_1
    + A_{1}^{\frac{1}{2}} e^5_3
    + A_{1}^{\frac{1}{2}} e^6_4
    + A_{2}^{\frac{1}{2}} e^8_7
    \\
    \pi ( {E^1}_3 )
    & = &
    e^5_7
    -
    \overline{q}
    e^2_4
    \\ 
    \pi ( {E^3}_1 )
    & = &
    e^7_5
    -
    q
    e^4_2
    \\
    \pi ( {E^2}_4 )
    & = &
    \overline{q}
    A_{0}^{\frac{1}{2}}
    e^1_3
    -
    A_{1}^{\frac{1}{2}}
    e^2_5
    +
    \overline{q}
    A_{1}^{\frac{1}{2}}
    e^4_7
    -
    A_{2}^{\frac{1}{2}}
    e^6_8
    \\
    \pi ( {E^4}_2 )
    & = &
    q
    A_{0}^{\frac{1}{2}}
    e^3_1
    -
    A_{1}^{\frac{1}{2}}
    e^5_2
    +
    q
    A_{1}^{\frac{1}{2}}
    e^7_4
    -
    A_{2}^{\frac{1}{2}}
    e^8_6
    \\
    \pi ( {E^1}_4 )
    & = &
    -
    \overline{q}^2
    A_{0}^{\frac{1}{2}}
    e^1_4
    -
    \overline{q}
    A_{1}^{\frac{1}{2}}
    e^2_6
    -
    \overline{q}
    A_{1}^{\frac{1}{2}}
    e^3_7
    +
    A_{2}^{\frac{1}{2}}
    e^5_8
    \\
    \pi ( {E^4}_1 )
    & = &
    -
    q^2
    A_{0}^{\frac{1}{2}}
    e^4_1
    -
    q
    A_{1}^{\frac{1}{2}}
    e^6_2
    -
    q
    A_{1}^{\frac{1}{2}}
    e^7_3
    +
    A_{2}^{\frac{1}{2}}
    e^8_5
  \end{array}
\end{eqnarray*}

%% file: data/basis30.tex
\subsection{The basis $\mathfrak{B}_0$ for
              $V_{0}\equiv V_{({0, 0, 0}\,|\,2 \alpha)}$}

The $8$ vectors in this basis are:
\small
\begin{eqnarray*}
  \hspace{-10pt}
  \begin{array}{r}
    \ketket{1}{1}
    \\
    \overline{C}_0^{\frac{1}{2}}
    q^{\frac{\alpha}{2}}_{+}
    \left\{
      \ketket{1}{2},
      \ketket{1}{3},
      \ketket{1}{4}
    \right\}
    \\
    \overline{C}_0^{\frac{1}{2}}
    \overline{A}_{2,1}^{\frac{1}{2}}
    \left\{
    \begin{array}{@{\hspace{0pt}}c@{\hspace{0pt}}}
      A_{1}^{\frac{1}{2}}
      q^{\alpha}_{+}
      \ketket{1}{5}  
      +
      A_{0}^{\frac{1}{2}} \overline{q}^{{\frac{1}{2}}}_{-}
      \ketket{2}{3},
      \\
      A_{1}^{\frac{1}{2}} q^{\alpha}_{+}
      \ketket{1}{6}
      +
      A_{0}^{\frac{1}{2}} \overline{q}^{{\frac{1}{2}}}_{-}
      \ketket{2}{4},
      \\
      A_{1}^{\frac{1}{2}} q^{\alpha}_{+}
      \ketket{1}{7}  
      +
      A_{0}^{\frac{1}{2}} \overline{q}^{{\frac{1}{2}}}_{-}
      \ketket{3}{4}
    \end{array}
    \right\}
    \\
    \overline{C}_{0}^{\frac{1}{2}}
    \overline{C}_{1}^{\frac{1}{2}}
    \overline{A}_{2,1}^{\frac{1}{2}}
    \left[
      \begin{array}{@{\hspace{0pt}}c@{\hspace{0pt}}}
        A_{0}^{\frac{1}{2}}
        \left(
          q^{\frac{\alpha}{2}}_{+}
          \ketket{3}{6}
          -
          q^{\frac{\alpha}{2}-1}_{+}
          \ketket{2}{7}
          -
          q^{\frac{\alpha}{2}+1}_{+}
          \ketket{4}{5}
        \right)
        -
        A_{2}^{\frac{1}{2}}
        q^{\frac{3\alpha}{2}}_{+}
        \ketket{1}{8}
      \end{array}
    \right]
  \end{array}
\end{eqnarray*}
\normalsize

%% file: data/basis31.tex
\subsection{The basis $\mathfrak{B}_{1}$ for
              $V_{1}\equiv V_{({0,0,-1}\,|\,2\alpha+1)}$}

The $24$ vectors in this basis are:
\small
\begin{eqnarray*}
  \hspace{-27pt}
  \begin{array}{r}
    \left\{
    \begin{array}{@{\hspace{0pt}}c@{\hspace{0pt}}}
      \ketket{2}{2},
      \ketket{3}{3},
      \ketket{4}{4}
    \end{array}
    \right\}
    \\
    \overline{C}_{0}^{\frac{1}{2}}
    \overline{q}^{{\frac{\alpha}{2}}}_{-}
    \left\{
    \begin{array}{@{\hspace{0pt}}c@{\hspace{0pt}}}
      \ketket{1}{2},
      \ketket{1}{3},
      \ketket{1}{4}
    \end{array}
    \right\}
    \\
    \overline{C}_{1}^{\frac{1}{2}}
    q^{\frac{\alpha}{2}+\frac{1}{2}}_{-}
    \left\{
    \begin{array}{@{\hspace{0pt}}c@{\hspace{0pt}}}
      \ketket{2}{5},
      \ketket{3}{5},
      \ketket{2}{6},
      \ketket{3}{7},
      \ketket{4}{7},
      \ketket{4}{6}
    \end{array}
    \right\}
    \\
    \overline{C}_{0}^{\frac{1}{2}}
    \overline{A}_{2,1}^{\frac{1}{2}}
    \left\{
    \begin{array}{@{\hspace{0pt}}c@{\hspace{0pt}}}
      A_{1}^{\frac{1}{2}} q^{0}_{-}
      \ketket{1}{5}
      +
      A_{0}^{\frac{1}{2}} \overline{q}^{\alpha+\frac{1}{2}}_{+}
      \ketket{2}{3},
      \\
      A_{1}^{\frac{1}{2}} q^{0}_{-}
      \ketket{1}{6}  
      +
      A_{0}^{\frac{1}{2}} \overline{q}^{\alpha+\frac{1}{2}}_{+}
      \ketket{2}{4},
      \\
      A_{1}^{\frac{1}{2}} q^{0}_{-}
      \ketket{1}{7}
      -
      A_{0}^{\frac{1}{2}} \overline{q}^{\alpha+\frac{1}{2}}_{+}
      \ketket{3}{4}
    \end{array}
    \right\}
    \\
    \overline{C}_{1}^{\frac{1}{2}}
    \overline{A}_{2,1}^{\frac{1}{2}}
    \left\{
    \begin{array}{@{\hspace{0pt}}c@{\hspace{0pt}}}
      A_{0}^{\frac{1}{2}} q^{0}_{-}
      \ketket{1}{6}
      -
      A_{1}^{\frac{1}{2}}
      q^{\alpha+\frac{1}{2}}_{+}
      \ketket{2}{4},
      \\
      A_{0}^{\frac{1}{2}} q^{0}_{-}
      \ketket{1}{5}  
      -
      A_{1}^{\frac{1}{2}}
      q^{\alpha+\frac{1}{2}}_{+}
      \ketket{2}{3},
      \\
      A_{0}^{\frac{1}{2}} q^{0}_{-}
      \ketket{1}{7}  
      -
      A_{1}^{\frac{1}{2}}
      q^{\alpha+\frac{1}{2}}_{+}
      \ketket{3}{4}
    \end{array}
    \right\}
    \\
    \overline{C}_{1}^{\frac{1}{2}}
    \overline{A}_{2,3}^{\frac{1}{2}}
    \left\{
    \begin{array}{@{\hspace{0pt}}c@{\hspace{0pt}}}
      A_{2}^{\frac{1}{2}}
      q^{\alpha+1}_{+}
      \ketket{2}{8}
      -
      A_{1}^{\frac{1}{2}}
      \overline{q}^{\frac{1}{2}}_{-}
      \ketket{5}{6},
      \\
      A_{2}^{\frac{1}{2}}
      q^{\alpha+1}_{+}
      \ketket{3}{8}  
      -
      A_{1}^{\frac{1}{2}}
      \overline{q}^{\frac{1}{2}}_{-}
      \ketket{5}{7},
      \\
      A_{2}^{\frac{1}{2}}
      q^{\alpha+1}_{+}
      \ketket{4}{8}
      -
      A_{1}^{\frac{1}{2}}
      \overline{q}^{\frac{1}{2}}_{-}
      \ketket{6}{7}
    \end{array}
    \right\}
    \\
    \overline{C}_{0}^{\frac{1}{2}}
    \overline{C}_{1}^{\frac{1}{2}}
    \overline{A}_{2,1}^{\frac{1}{2}}
    \left[
      \begin{array}{@{\hspace{0pt}}c@{\hspace{0pt}}}
        A_{2}^{\frac{1}{2}}
        q^{\frac{\alpha}{2}}_{-}
        \ketket{1}{8}
        +
        A_{0}^{\frac{1}{2}}
        \overline{q}^{\frac{\alpha}{2}+1}_{-}
        \ketket{2}{7}
        -
        A_{0}^{\frac{1}{2}}
        \overline{q}^{\frac{\alpha}{2}}_{-}
        \ketket{3}{6}
        -
        A_{0}^{\frac{1}{2}}
        q^{\frac{3\alpha}{2}+1}_{-}
        \ketket{4}{5}
      \end{array}
    \right]
    \\
    \overline{C}_1
    \overline{A}_{1}^{\frac{1}{2}}
    \overline{A}_{2,1}^{\frac{1}{2}}
    \left[
      \begin{array}{@{\hspace{0pt}}c@{\hspace{0pt}}}
        A_{0}^{\frac{1}{2}}
        A_{2}^{\frac{1}{2}}
        q^{\frac{\alpha}{2}}_{-}
        \ketket{1}{8}
        +
        A_{0}
        \overline{q}^{\frac{\alpha}{2}+1}_{-}
        \ketket{2}{7}
        +
        A_{1}
        q^{\frac{3\alpha}{2}+1}_{-}
        \ketket{3}{6}
        +
        A_{1}
        \overline{q}^{\frac{\alpha}{2}}_{-}
        \ketket{4}{5}
      \end{array}
    \right]
    \\
    \overline{C}_1
    \overline{A}_{1}^{\frac{1}{2}}
    \overline{A}_{2,3}^{\frac{1}{2}}
    \left[
      \begin{array}{@{\hspace{0pt}}c@{\hspace{0pt}}}
        A_{0}^{\frac{1}{2}}
        A_{2}^{\frac{1}{2}}
        q^{\frac{\alpha}{2}}_{-}
        \ketket{1}{8}
        -
        A_{2}
        q^{\frac{3\alpha}{2}+1}_{-}
        \ketket{2}{7}
        -
        A_{1}
        \overline{q}^{\frac{\alpha}{2}+1}_{-}
        \ketket{3}{6}
        +
        A_{1}
        \overline{q}^{\frac{\alpha}{2}}_{-}
        \ketket{4}{5}
      \end{array}
    \right]
  \end{array}
\end{eqnarray*}

\normalsize

%% file: data/basis32.tex
\subsection{The basis $\mathfrak{B}_2$ for
              $V_{2}\equiv V_{({0,-1,-1}\,|\,2\alpha+2)}$}

The $24$ vectors in this basis are:
\small
\begin{eqnarray*}
  \hspace{-36pt}
  \begin{array}{r}
    \overline{C}_1
    \overline{A}_{1}^{\frac{1}{2}}
    \overline{A}_{2,1}^{\frac{1}{2}}
    \left[
    \begin{array}{@{\hspace{0pt}}c@{\hspace{0pt}}}
      A_{0}^{\frac{1}{2}}
      A_{2}^{\frac{1}{2}}
      \overline{q}^{\frac{\alpha}{2}+1}_{+}\ketket{1}{8}
      +
      A_{0}
      \overline{q}^{\frac{3\alpha}{2}+2}_{+}\ketket{2}{7}
      +
      A_{1} 
      q^{\frac{\alpha}{2}}_{+}\ketket{3}{6}
      -
      A_{1} 
      q^{\frac{\alpha}{2}+1}_{+}\ketket{4}{5}
    \end{array}
    \right]
    \\
    \overline{C}_1
    \overline{A}_{1}^{\frac{1}{2}}
    \overline{A}_{2,3}^{\frac{1}{2}}
    \left[
    \begin{array}{@{\hspace{0pt}}c@{\hspace{0pt}}}
      A_{0}^{\frac{1}{2}}
      A_{2}^{\frac{1}{2}}
      \overline{q}^{\frac{\alpha}{2}+1}_{+}
      \ketket{1}{8}
      -
      A_{2} 
      q^{\frac{\alpha}{2}}_{+}
      \ketket{2}{7}
      -
      A_{1}
      \overline{q}^{\frac{3\alpha}{2}+2}_{+}
      \ketket{3}{6}
      -
      A_{1} 
      q^{\frac{\alpha}{2}+1}_{+} \ketket{4}{5}
    \end{array}
    \right]
    \\
    \overline{C}_1^{\frac{1}{2}}
    \overline{C}_2^{\frac{1}{2}}
    \overline{A}_{2,3}^{\frac{1}{2}}
    \left[
    \begin{array}{@{\hspace{0pt}}c@{\hspace{0pt}}}
      A_{0}^{\frac{1}{2}}
      \overline{q}^{\frac{\alpha}{2}+1}_{+}\ketket{1}{8}
      -
      A_{2}^{\frac{1}{2}}
      q^{\frac{\alpha}{2}}_{+}\ketket{2}{7}
      +
      A_{2}^{\frac{1}{2}} 
      q^{\frac{\alpha}{2}+1}_{+}\ketket{3}{6}
      +
      A_{2}^{\frac{1}{2}}
      \overline{q}^{\frac{3\alpha}{2}+2}_{+}\ketket{4}{5}
    \end{array}
    \right]
    \\
    \overline{C}_1^{\frac{1}{2}}
    \overline{A}_{2,1}^{\frac{1}{2}}
    \left\{
    \begin{array}{@{\hspace{0pt}}c@{\hspace{0pt}}}
      A_{0}^{\frac{1}{2}}
      \overline{q}^{\alpha+1}_{+}
      \ketket{1}{5}
      -
      A_{1}^{\frac{1}{2}}
      \overline{q}^{\frac{1}{2}}_{-}
      \ketket{2}{3},
      \\
      A_{0}^{\frac{1}{2}}
      \overline{q}^{\alpha+1}_{+}
      \ketket{1}{6} 
      -
      A_{1}^{\frac{1}{2}}
      \overline{q}^{\frac{1}{2}}_{-}
      \ketket{2}{4},
      \\
      A_{0}^{\frac{1}{2}}
      \overline{q}^{\alpha+1}_{+}
      \ketket{1}{7}
      -
      A_{1}^{\frac{1}{2}}
      \overline{q}^{\frac{1}{2}}_{-}
      \ketket{3}{4}
    \end{array}
    \right\}
    \\
    \overline{C}_1^{\frac{1}{2}}
    \overline{A}_{2,3}^{\frac{1}{2}}
    \left\{
    \begin{array}{@{\hspace{0pt}}c@{\hspace{0pt}}}
      A_{2}^{\frac{1}{2}}
      q^{0}_{-}
      \ketket{2}{8} 
      -
      A_{1}^{\frac{1}{2}}
      \overline{q}^{\alpha+\frac{3}{2}}_{+}
      \ketket{5}{6},
      \\
      A_{2}^{\frac{1}{2}}
      q^{0}_{-}
      \ketket{3}{8}
      -
      A_{1}^{\frac{1}{2}}
      \overline{q}^{\alpha+\frac{3}{2}}_{+}
      \ketket{5}{7},
      \\
      A_{2}^{\frac{1}{2}}
      q^{0}_{-}
      \ketket{4}{8} 
      -
      A_{1}^{\frac{1}{2}}
      \overline{q}^{\alpha+\frac{3}{2}}_{+}
      \ketket{6}{7}
    \end{array}
    \right\}
    \\
    \overline{C}_2^{\frac{1}{2}}
    \overline{A}_{2,3}^{\frac{1}{2}}
    \left\{
    \begin{array}{@{\hspace{0pt}}c@{\hspace{0pt}}}
      A_{1}^{\frac{1}{2}}
      q^{0}_{-}\ketket{3}{8} 
      +
      A_{2}^{\frac{1}{2}} 
      q^{\alpha+\frac{3}{2}}_{+}
      \ketket{5}{7},
      \\
      A_{1}^{\frac{1}{2}}
      q^{0}_{-}
      \ketket{4}{8}
      +
      A_{2}^{\frac{1}{2}} 
      q^{\alpha+\frac{3}{2}}_{+}
      \ketket{6}{7},
      \\
      A_{1}^{\frac{1}{2}}
      q^{0}_{-}
      \ketket{2}{8}
      +
      A_{2}^{\frac{1}{2}} 
      q^{\alpha+\frac{3}{2}}_{+}
      \ketket{5}{6}
    \end{array}
    \right\}
    \\
    \overline{C}_1^{\frac{1}{2}}
    \overline{q}^{\frac{\alpha}{2}+\frac{1}{2}}_{+}
    \left\{
    \begin{array}{@{\hspace{0pt}}c@{\hspace{0pt}}}
      \ketket{2}{6},
      \ketket{2}{5},
      \ketket{3}{5},
      \ketket{3}{7},
      \ketket{4}{6},
      \ketket{4}{7}
    \end{array}
    \right\}
    \\
    \overline{C}_2^{\frac{1}{2}}
    q^{\frac{\alpha}{2}+1}_{+}
    \left\{
    \begin{array}{@{\hspace{0pt}}c@{\hspace{0pt}}}
      \ketket{5}{8},
      \ketket{7}{8},
      \ketket{6}{8}
    \end{array}
    \right\}
    \\
    \left\{
    \begin{array}{@{\hspace{0pt}}c@{\hspace{0pt}}}
      \ketket{7}{7},
      \ketket{6}{6},
      \ketket{5}{5}
    \end{array}
    \right\}
  \end{array}
\end{eqnarray*}

\normalsize

%% file: data/basis33.tex
\subsection{The basis $\mathfrak{B}_{3}$ for
              $V_{3}\equiv V_{({-1,-1,-1}\,|\,2\alpha+3)}$}

The $8$ vectors in this basis are:
\small
\begin{eqnarray*}
  \hspace{-20pt}
  \begin{array}{r}
    \overline{C}_1^{\frac{1}{2}}
    \overline{C}_2^{\frac{1}{2}}
    \overline{A}_{2,3}^{\frac{1}{2}}
    \hspace{-1pt}
    \left[
    \begin{array}{@{\hspace{0pt}}c@{\hspace{0pt}}}
      A_{0}^{\frac{1}{2}}
      \overline{q}^{\frac{3\alpha}{2}+3}_{-}
      \ketket{1}{8}
      -
      A_{2}^{\frac{1}{2}}
      \left(
        \overline{q}^{\frac{\alpha}{2}+2}_{-}
        \ketket{2}{7}
        -
        \overline{q}^{\frac{\alpha}{2}+1}_{-}
        \ketket{3}{6}
        +
        \overline{q}^{\frac{\alpha}{2}}_{-}
        \ketket{4}{5}
      \right)
    \end{array}
    \right]
    \\
    \overline{C}_2^{\frac{1}{2}}
    \overline{A}_{2,3}^{\frac{1}{2}}
    \left\{
    \begin{array}{@{\hspace{0pt}}c@{\hspace{0pt}}}
      A_{1}^{\frac{1}{2}}
      \overline{q}^{\alpha+2}_{+}
      \ketket{2}{8}  
      +
      A_{2}^{\frac{1}{2}}
      \overline{q}^{\frac{1}{2}}_{-}
      \ketket{5}{6},
      \\
      A_{1}^{\frac{1}{2}}
      \overline{q}^{\alpha+2}_{+}
      \ketket{3}{8}
      +
      A_{2}^{\frac{1}{2}}
      \overline{q}^{\frac{1}{2}}_{-}
      \ketket{5}{7},
      \\
      A_{1}^{\frac{1}{2}}
      \overline{q}^{\alpha+2}_{+}
      \ketket{4}{8}  
      +
      A_{2}^{\frac{1}{2}}
      \overline{q}^{\frac{1}{2}}_{-}
      \ketket{6}{7}
    \end{array}
    \right\}
    \\
    \overline{C}_2^{\frac{1}{2}}
    \overline{q}^{\frac{\alpha}{2}+1}_{-}
    \left\{
    \begin{array}{@{\hspace{0pt}}c@{\hspace{0pt}}}
      \ketket{5}{8},
      \ketket{6}{8},
      \ketket{7}{8}
    \end{array}
    \right\}
    \\
    \ketket{8}{8}
  \end{array}
\end{eqnarray*}

\normalsize

%% file: data/trigrmx3.tex
\subsection{The trigonometric R matrix {$\check{R}^{3,1}(u)$}}

For the listing of the components of $\check{R}^{3,1}(u)$, we invoke a
little more notation:
\begin{eqnarray*}
  S^\pm_i
  & \triangleq &
  [ \alpha + i \pm u ]_q,
  \\
  U_i^z
  & \triangleq &
  {[ u - i ]_q}^z,
  \qquad
  \mathrm{where~} z\in \{1, 2\},
\end{eqnarray*}
and $i \in \{ 0,1,2\}$. With this, $\check{R}^{3,1} (u)$ has $216$ nonzero
components:

\small

\begin{eqnarray*}
  & &
  \hspace{-45pt}
  1
  \left\{
    \begin{array}{@{\hspace{0mm}}c@{\hspace{0mm}}}
      e^{1 1}_{1 1}
    \end{array}
  \right\},
  \qquad
  \frac{
    S_{0}^{+}
  }{
    S_{0}^{-}
  }
  \left\{
    \begin{array}{@{\hspace{0mm}}c@{\hspace{0mm}}}
      e^{2 2}_{2 2},
      e^{3 3}_{3 3},
      e^{4 4}_{4 4}
    \end{array}
  \right\},
  \qquad
  \frac{
    S_{0}^{+}
    S_{1}^{+}
  }{
    S_{0}^{-}
    S_{1}^{-}
  }
  \left\{
    \begin{array}{@{\hspace{0mm}}c@{\hspace{0mm}}}
      e^{5 5}_{5 5},
      e^{6 6}_{6 6},
      e^{7 7}_{7 7}
    \end{array}
  \right\},
  \qquad
  \frac{
    S_{0}^{+}
    S_{1}^{+}
    S_{2}^{+}
  }{
    S_{0}^{-}
    S_{1}^{-}
    S_{2}^{-}
  }
  \left\{
    \begin{array}{@{\hspace{0mm}}c@{\hspace{0mm}}}
      e^{8 8}_{8 8}
    \end{array}
  \right\}
  \\[1mm]
  \hline
  \\[1mm]
  & &
  \hspace{-45pt}
  \frac{
    A_{0}
  }{
    S_{0}^{-}
  }
  \left\{
    \begin{array}{@{\hspace{0mm}}c@{\hspace{0mm}}}
      \overline{q}^{u}
      \left\{
      \begin{array}{@{\hspace{0mm}}c@{\hspace{0mm}}}
        e^{1 2}_{1 2},
        e^{1 3}_{1 3},
        e^{1 4}_{1 4}
      \end{array}
      \right\}
      \\
      q^{u}
      \left\{
      \begin{array}{@{\hspace{0mm}}c@{\hspace{0mm}}}
        e^{2 1}_{2 1},
        e^{3 1}_{3 1},
        e^{4 1}_{4 1}
      \end{array}
      \right\}
    \end{array}
  \right\},
  \qquad
  \frac{
    A_{2}
    S_{0}^{+}
    S_{1}^{+}
  }{
    S_{0}^{-}
    S_{1}^{-}
    S_{2}^{-}
  }
  \left\{
    \begin{array}{@{\hspace{0mm}}c@{\hspace{0mm}}}
      q^{u}
      \left\{
      \begin{array}{@{\hspace{0mm}}c@{\hspace{0mm}}}
        e^{8 7}_{8 7},
        e^{8 6}_{8 6},
        e^{8 5}_{8 5}
      \end{array}
      \right\}
      \\
      \overline{q}^{ u}
      \left\{
      \begin{array}{@{\hspace{0mm}}c@{\hspace{0mm}}}
        e^{7 8}_{7 8},
        e^{6 8}_{6 8},
        e^{5 8}_{5 8}
      \end{array}
      \right\}
    \end{array}
  \right\}
  \\
  & &
  \hspace{-45pt}
  \frac{
    1
  }{
    \Delta^2
    S_{0}^{-}
    S_{1}^{-}
  }
  \left\{
    \begin{array}{@{\hspace{0mm}}c@{\hspace{0mm}}}
      f_1 (q)
      \left\{
      \begin{array}{@{\hspace{0mm}}c@{\hspace{0mm}}}
        e^{2 3}_{2 3},
        e^{2 4}_{2 4},
        e^{3 4}_{3 4}
      \end{array}
      \right\}
      \\
      f_1 (\overline{q})
      \left\{
      \begin{array}{@{\hspace{0mm}}c@{\hspace{0mm}}}
        e^{3 2}_{3 2},
        e^{4 2}_{4 2},
        e^{4 3}_{4 3}
      \end{array}
      \right\}
    \end{array}
  \right\},
  \qquad
  \frac{
    1
  }{
  \Delta^2
    S_{0}^{-}
    S_{1}^{-}
    S_{2}^{-}
  }
  \left\{
    \begin{array}{@{\hspace{0mm}}c@{\hspace{0mm}}}
      f_2 (q)
      \left\{
      \begin{array}{@{\hspace{0mm}}c@{\hspace{0mm}}}
        e^{7 6}_{7 6},
        e^{7 5}_{7 5},
        e^{6 5}_{6 5}
      \end{array}
      \right\}
      \\
      f_2 (\overline{q})
      \left\{
      \begin{array}{@{\hspace{0mm}}c@{\hspace{0mm}}}
        e^{6 7}_{6 7},
        e^{5 7}_{5 7},
        e^{5 6}_{5 6}
      \end{array}
      \right\}
    \end{array}
  \right\}
  \\
  & &
  \hspace{-45pt}
  \frac{
    A_{0}
    A_{1}
  }{
    S_{0}^{-}
    S_{1}^{-}
  }
  \left\{
    \begin{array}{@{\hspace{0mm}}c@{\hspace{0mm}}}
      q^{2 u}
      \left\{
      \begin{array}{@{\hspace{0mm}}c@{\hspace{0mm}}}
        e^{5 1}_{5 1},
        e^{6 1}_{6 1},
        e^{7 1}_{7 1}
      \end{array}
      \right\}
      \\
      \overline{q}^{2 u}
      \left\{
      \begin{array}{@{\hspace{0mm}}c@{\hspace{0mm}}}
        e^{1 5}_{1 5},
        e^{1 6}_{1 6},
        e^{1 7}_{1 7}
      \end{array}
      \right\}
    \end{array}
  \right\},
  \qquad
  \frac{
    A_{1}
    A_{2}
    S_{0}^{+}
  }{
    S_{0}^{-}
    S_{1}^{-}
    S_{2}^{-}
  }
  \left\{
    \begin{array}{@{\hspace{0mm}}c@{\hspace{0mm}}}
      q^{2 u}
      \left\{
      \begin{array}{@{\hspace{0mm}}c@{\hspace{0mm}}}
        e^{8 4}_{8 4},
        e^{8 3}_{8 3},
        e^{8 2}_{8 2}
      \end{array}
      \right\}
      \\
      \overline{q}^{2 u}
      \left\{
      \begin{array}{@{\hspace{0mm}}c@{\hspace{0mm}}}
        e^{4 8}_{4 8},
        e^{3 8}_{3 8},
        e^{2 8}_{2 8}
      \end{array}
      \right\}
    \end{array}
  \right\}
  \\
  & &
  \hspace{-45pt}
  \frac{
    A_{1}
  }{
  \Delta^2
    S_{0}^{-}
    S_{1}^{-}
    S_{2}^{-}
  }
  \left\{
    \begin{array}{@{\hspace{0mm}}c@{\hspace{0mm}}}
      f_3 (q)
      \left\{
      \begin{array}{@{\hspace{0mm}}c@{\hspace{0mm}}}
        e^{2 7}_{2 7}
      \end{array}
      \right\},
      f_4 (q)
      \left\{
      \begin{array}{@{\hspace{0mm}}c@{\hspace{0mm}}}
        e^{3 6}_{3 6}
      \end{array}
      \right\},
      f_5 (q)
      \left\{
      \begin{array}{@{\hspace{0mm}}c@{\hspace{0mm}}}
        e^{4 5}_{4 5}
      \end{array}
      \right\}
      \\
      f_3 (\overline{q})
      \left\{
      \begin{array}{@{\hspace{0mm}}c@{\hspace{0mm}}}
        e^{7 2}_{7 2}
      \end{array}
      \right\},
      f_4 (\overline{q})
      \left\{
      \begin{array}{@{\hspace{0mm}}c@{\hspace{0mm}}}
        e^{6 3}_{6 3}
      \end{array}
      \right\},
      f_5 (\overline{q})
      \left\{
      \begin{array}{@{\hspace{0mm}}c@{\hspace{0mm}}}
        e^{5 4}_{5 4}
      \end{array}
      \right\}
    \end{array}
  \right\}
  \\
  & &
  \hspace{-45pt}
  \frac{
    A_{1}
    S_{0}^{+}
  }{
    S_{0}^{-}
    S_{1}^{-}
  }
  \left\{
    \begin{array}{@{\hspace{0mm}}c@{\hspace{0mm}}}
    q^{u}
    \left\{
    \begin{array}{@{\hspace{0mm}}c@{\hspace{0mm}}}
      e^{5 2}_{5 2},
      e^{5 3}_{5 3},
      e^{6 2}_{6 2},
      e^{6 4}_{6 4},
      e^{7 3}_{7 3},
      e^{7 4}_{7 4}
    \end{array}
    \right\}
    \\
    \overline{q}^{u}
    \left\{
    \begin{array}{@{\hspace{0mm}}c@{\hspace{0mm}}}
      e^{2 5}_{2 5},
      e^{3 5}_{3 5},
      e^{2 6}_{2 6},
      e^{4 6}_{4 6},
      e^{3 7}_{3 7},
      e^{4 7}_{4 7}
    \end{array}
    \right\}
  \end{array}
  \right\},
  \qquad
  \frac{
    A_{0}
    A_{1}
    A_{2}
  }{
    S_{0}^{-}
    S_{1}^{-}
    S_{2}^{-}
  }
  \left\{
  \begin{array}{@{\hspace{0mm}}c@{\hspace{0mm}}}
    q^{3 u}
    \left\{
    \begin{array}{@{\hspace{0mm}}c@{\hspace{0mm}}}
      e^{8 1}_{8 1}
    \end{array}
    \right\}
      \\
    \overline{q}^{3 u}
    \left\{
    \begin{array}{@{\hspace{0mm}}c@{\hspace{0mm}}}
      e^{1 8}_{1 8}
    \end{array}
    \right\}
  \end{array}
  \right\}
  \\[1mm]
  \hline
  \\[1mm]
  & &
  \hspace{-45pt}
  \frac{
    U_{0}
  }{
    S_{0}^{-}
  }
  \left\{
    \begin{array}{@{\hspace{0mm}}c@{\hspace{0mm}}}
      \mathbf{+ 1}
      \left\{
      \begin{array}{@{\hspace{0mm}}c@{\hspace{0mm}}}
        e^{1 2}_{2 1},
        e^{1 3}_{3 1},
        e^{1 4}_{4 1}
      \end{array}
      \right\}
      \\
      - 1
      \left\{
      \begin{array}{@{\hspace{0mm}}c@{\hspace{0mm}}}
        e^{2 1}_{1 2},
        e^{3 1}_{1 3},
        e^{4 1}_{1 4}
      \end{array}
      \right\}
    \end{array}
  \right\},
  \;\;
  \frac{
    U_{0}
    S_{0}^{+}
    S_{1}^{+}
  }{
    S_{0}^{-}
    S_{1}^{-}
    S_{2}^{-}
  }
  \left\{
    \begin{array}{@{\hspace{0mm}}c@{\hspace{0mm}}}
      \mathbf{+ 1}
      \left\{
      \begin{array}{@{\hspace{0mm}}c@{\hspace{0mm}}}
        e^{7 8}_{8 7},
        e^{6 8}_{8 6},
        e^{5 8}_{8 5}
      \end{array}
      \right\}
      \\
      - 1
      \left\{
      \begin{array}{@{\hspace{0mm}}c@{\hspace{0mm}}}
        e^{8 7}_{7 8},
        e^{8 6}_{6 8},
        e^{8 5}_{5 8}
      \end{array}
      \right\}
    \end{array}
  \right\},
  \;\;
  \frac{
    U_{0}
    U_{1}
    U_{2}
  }{
    S_{0}^{-}
    S_{1}^{-}
    S_{2}^{-}
  }
  \left\{
    \begin{array}{@{\hspace{0mm}}c@{\hspace{0mm}}}
      \mathbf{+ 1}
      \left\{
      \begin{array}{@{\hspace{0mm}}c@{\hspace{0mm}}}
        e^{1 8}_{8 1}
      \end{array}
      \right\}
      \\
      - 1
      \left\{
      \begin{array}{@{\hspace{0mm}}c@{\hspace{0mm}}}
        e^{8 1}_{1 8}
      \end{array}
      \right\}
    \end{array}
  \right\}
  \\
  & &
  \hspace{-45pt}
  \frac{
    U_{0}
    U_{1}
  }{
    S_{0}^{-}
    S_{1}^{-}
  }
  \left\{
    \begin{array}{@{\hspace{0mm}}c@{\hspace{0mm}}}
        e^{1 5}_{5 1},
        e^{1 6}_{6 1},
        e^{1 7}_{7 1}
      \\
        e^{5 1}_{1 5},
        e^{6 1}_{1 6},
        e^{7 1}_{1 7}
    \end{array}
  \right\},
  \qquad
  -
  \frac{
    U_{0}
    U_{1}
    S_{0}^{+}
  }{
    S_{0}^{-}
    S_{1}^{-}
    S_{2}^{-}
  }
  \left\{
    \begin{array}{@{\hspace{0mm}}c@{\hspace{0mm}}}
      e^{8 4}_{4 8},
      e^{8 3}_{3 8},
      e^{8 2}_{2 8}
      \\
      e^{4 8}_{8 4},
      e^{3 8}_{8 3},
      e^{2 8}_{8 2}
    \end{array}
  \right\}
  \\
  & &
  \hspace{-45pt}
  -
  \frac{
    U_{0}^2
  }{
    S_{0}^{-}
    S_{1}^{-}
  }
  \left\{
    \begin{array}{@{\hspace{0mm}}c@{\hspace{0mm}}}
        e^{2 3}_{3 2},
        e^{2 4}_{4 2},
        e^{3 4}_{4 3}
      \\
        e^{3 2}_{2 3},
        e^{4 2}_{2 4},
        e^{4 3}_{3 4}
    \end{array}
  \right\},
  \qquad
  -
  \frac{
    U_{0}^2
    S_{0}^{+}
  }{
    S_{0}^{-}
    S_{1}^{-}
    S_{2}^{-}
  }
  \left\{
    \begin{array}{@{\hspace{0mm}}c@{\hspace{0mm}}}
        e^{7 6}_{6 7},
        e^{7 5}_{5 7},
        e^{6 5}_{5 6}
      \\
        e^{6 7}_{7 6},
        e^{5 7}_{7 5},
        e^{5 6}_{6 5}
    \end{array}
  \right\}
  \\
  & &
  \hspace{-45pt}
  \frac{
    U_{0}^2
    U_{1}
  }{
    S_{0}^{-}
    S_{1}^{-}
    S_{2}^{-}
  }
  \left\{
    \begin{array}{@{\hspace{0mm}}c@{\hspace{0mm}}}
      \mathbf{+ 1}
      \left\{
      \begin{array}{@{\hspace{0mm}}c@{\hspace{0mm}}}
        e^{5 4}_{4 5},
        e^{6 3}_{3 6},
        e^{7 2}_{2 7}
      \end{array}
      \right\}
      \\
      - 1
      \left\{
      \begin{array}{@{\hspace{0mm}}c@{\hspace{0mm}}}
        e^{4 5}_{5 4},
        e^{3 6}_{6 3},
        e^{2 7}_{7 2}
      \end{array}
      \right\}
    \end{array}
  \right\},
  \qquad
  \frac{
    U_{0}
    S_{0}^{+}
  }{
    S_{0}^{-}
    S_{1}^{-}
  }
  \left\{
    \begin{array}{@{\hspace{0mm}}c@{\hspace{0mm}}}
      \mathbf{- 1}
      \left\{
      \begin{array}{@{\hspace{0mm}}c@{\hspace{0mm}}}
        e^{5 2}_{2 5},
        e^{5 3}_{3 5},
        e^{6 2}_{2 6},
        e^{6 4}_{4 6},
        e^{7 3}_{3 7},
        e^{7 4}_{4 7}
      \end{array}
      \right\}
      \\
      + 1
      \left\{
      \begin{array}{@{\hspace{0mm}}c@{\hspace{0mm}}}
        e^{2 5}_{5 2},
        e^{3 5}_{5 3},
        e^{2 6}_{6 2},
        e^{4 6}_{6 4},
        e^{3 7}_{7 3},
        e^{4 7}_{7 4}
      \end{array}
      \right\}
    \end{array}
  \right\}
  \\[1mm]
  \hline
  \\[1mm]
  & &
  \hspace{-45pt}
  \frac{
    A_{0}^{\frac{1}{2}}
    A_{1}^{\frac{1}{2}}
    U_{0}
  }{
    S_{0}^{-}
    S_{1}^{-}
  }
  \left\{
    \begin{array}{@{\hspace{0mm}}c@{\hspace{0mm}}}
    q^{u + \frac{1}{2}}
    \left\{
    \begin{array}{@{\hspace{0mm}}c@{\hspace{0mm}}}
      \mathbf{- 1}
      \left\{
      \begin{array}{@{\hspace{0mm}}c@{\hspace{0mm}}}
        e^{5 1}_{3 2},
        e^{6 1}_{4 2},
        e^{7 1}_{4 3}
      \end{array}
      \right\}
      \\
      + 1
      \left\{
      \begin{array}{@{\hspace{0mm}}c@{\hspace{0mm}}}
        e^{3 2}_{5 1},
        e^{4 2}_{6 1},
        e^{4 3}_{7 1}
      \end{array}
      \right\}
    \end{array}
    \right\}
    \\
    q^{u - \frac{1}{2}}
    \left\{
    \begin{array}{@{\hspace{0mm}}c@{\hspace{0mm}}}
      \mathbf{+ 1}
      \left\{
      \begin{array}{@{\hspace{0mm}}c@{\hspace{0mm}}}
        e^{5 1}_{2 3},
        e^{6 1}_{2 4},
        e^{7 1}_{3 4}
      \end{array}
      \right\}
      \\
      - 1
      \left\{
      \begin{array}{@{\hspace{0mm}}c@{\hspace{0mm}}}
        e^{2 3}_{5 1},
        e^{2 4}_{6 1},
        e^{3 4}_{7 1}
      \end{array}
      \right\}
    \end{array}
    \right\}
    \\
    \overline{q}^{u + \frac{1}{2}}
    \left\{
    \begin{array}{@{\hspace{0mm}}c@{\hspace{0mm}}}
      \mathbf{+ 1}
      \left\{
      \begin{array}{@{\hspace{0mm}}c@{\hspace{0mm}}}
        e^{1 5}_{2 3},
        e^{1 6}_{2 4},
        e^{1 7}_{3 4}
      \end{array}
      \right\}
      \\
      - 1
      \left\{
      \begin{array}{@{\hspace{0mm}}c@{\hspace{0mm}}}
        e^{2 3}_{1 5},
        e^{2 4}_{1 6},
        e^{3 4}_{1 7}
      \end{array}
      \right\}
    \end{array}
    \right\}
    \\
    \overline{q}^{u - \frac{1}{2}}
    \left\{
    \begin{array}{@{\hspace{0mm}}c@{\hspace{0mm}}}
      \mathbf{- 1}
      \left\{
      \begin{array}{@{\hspace{0mm}}c@{\hspace{0mm}}}
        e^{1 5}_{3 2},
        e^{1 6}_{4 2},
        e^{1 7}_{4 3}
      \end{array}
      \right\}
      \\
      + 1
      \left\{
      \begin{array}{@{\hspace{0mm}}c@{\hspace{0mm}}}
        e^{3 2}_{1 5},
        e^{4 2}_{1 6},
        e^{4 3}_{1 7}
      \end{array}
      \right\}
    \end{array}
    \right\}
  \end{array}
  \right\},
  \hspace{-1pt}
  \frac{
    A_{1}^{\frac{1}{2}}
    A_{2}^{\frac{1}{2}}
    U_{0}
    S_{0}^{+}
  }{
    S_{0}^{-}
    S_{1}^{-}
    S_{2}^{-}
  }
  \left\{
    \begin{array}{@{\hspace{0mm}}c@{\hspace{0mm}}}
    q^{u + \frac{1}{2}}
    \left\{
    \begin{array}{@{\hspace{0mm}}c@{\hspace{0mm}}}
      \mathbf{+ 1}
      \left\{
      \begin{array}{@{\hspace{0mm}}c@{\hspace{0mm}}}
        e^{6 5}_{8 2},
        e^{7 5}_{8 3},
        e^{7 6}_{8 4}
      \end{array}
      \right\}
      \\
      - 1
      \left\{
      \begin{array}{@{\hspace{0mm}}c@{\hspace{0mm}}}
        e^{8 2}_{6 5},
        e^{8 3}_{7 5},
        e^{8 4}_{7 6}
      \end{array}
      \right\}
    \end{array}
    \right\}
    \\
    q^{u - \frac{1}{2}}
    \left\{
    \begin{array}{@{\hspace{0mm}}c@{\hspace{0mm}}}
      \mathbf{- 1}
      \left\{
      \begin{array}{@{\hspace{0mm}}c@{\hspace{0mm}}}
        e^{5 6}_{8 2},
        e^{5 7}_{8 3},
        e^{6 7}_{8 4}
      \end{array}
      \right\}
      \\
      + 1
      \left\{
      \begin{array}{@{\hspace{0mm}}c@{\hspace{0mm}}}
        e^{8 2}_{5 6},
        e^{8 3}_{5 7},
        e^{8 4}_{6 7}
      \end{array}
      \right\}
    \end{array}
    \right\}
    \\
    \overline{q}^{u + \frac{1}{2}}
    \left\{
    \begin{array}{@{\hspace{0mm}}c@{\hspace{0mm}}}
      \mathbf{- 1}
      \left\{
      \begin{array}{@{\hspace{0mm}}c@{\hspace{0mm}}}
        e^{5 6}_{2 8},
        e^{5 7}_{3 8},
        e^{6 7}_{4 8}
      \end{array}
      \right\}
      \\
      + 1
      \left\{
      \begin{array}{@{\hspace{0mm}}c@{\hspace{0mm}}}
        e^{2 8}_{5 6},
        e^{3 8}_{5 7},
        e^{4 8}_{6 7}
      \end{array}
      \right\}
    \end{array}
    \right\}
    \\
    \overline{q}^{u - \frac{1}{2}}
    \left\{
    \begin{array}{@{\hspace{0mm}}c@{\hspace{0mm}}}
      \mathbf{+ 1}
      \left\{
      \begin{array}{@{\hspace{0mm}}c@{\hspace{0mm}}}
        e^{6 5}_{2 8},
        e^{7 5}_{3 8},
        e^{7 6}_{4 8}
      \end{array}
      \right\}
      \\
      - 1
      \left\{
      \begin{array}{@{\hspace{0mm}}c@{\hspace{0mm}}}
        e^{2 8}_{6 5},
        e^{3 8}_{7 5},
        e^{4 8}_{7 6}
      \end{array}
      \right\}
    \end{array}
    \right\}
  \end{array}
  \right\}
  \\
  & &
  \hspace{-45pt}
  \frac{
    U_{0}
  }{
    \Delta^{2}
    S_{0}^{-}
    S_{1}^{-}
    S_{2}^{-}
  }
  \left\{
  \begin{array}{@{\hspace{0mm}}c@{\hspace{0mm}}}
    f_6(q)
    \left\{
    \begin{array}{@{\hspace{0mm}}c@{\hspace{0mm}}}
      \mathbf{+ 1}
      \left\{
      \begin{array}{@{\hspace{0mm}}c@{\hspace{0mm}}}
            e^{6 3}_{4 5},
        - q e^{7 2}_{4 5},
            e^{7 2}_{3 6}
      \end{array}
      \right\}
      \\
      - 1
      \left\{
      \begin{array}{@{\hspace{0mm}}c@{\hspace{0mm}}}
            e^{4 5}_{6 3},
        - q e^{4 5}_{7 2},
            e^{3 6}_{7 2}
      \end{array}
      \right\}
    \end{array}
    \right\},
    f_6 (\overline{q})
    \left\{
    \begin{array}{@{\hspace{0mm}}c@{\hspace{0mm}}}
      \mathbf{+ 1}
      \left\{
      \begin{array}{@{\hspace{0mm}}c@{\hspace{0mm}}}
                       e^{5 4}_{3 6},
        - \overline{q} e^{5 4}_{2 7},
                       e^{6 3}_{2 7}
      \end{array}
      \right\}
      \\
      - 1
      \left\{
      \begin{array}{@{\hspace{0mm}}c@{\hspace{0mm}}}
                       e^{3 6}_{5 4},
        - \overline{q} e^{2 7}_{5 4},
                       e^{2 7}_{6 3}
      \end{array}
      \right\}
    \end{array}
    \right\}
  \end{array}
  \right\}
  \\
  & &
  \hspace{-45pt}
  \frac{
    A_{1}
    U_{0}^2
  }{
    S_{0}^{-}
    S_{1}^{-}
    S_{2}^{-}
  }
  \left\{
  \begin{array}{@{\hspace{0mm}}c@{\hspace{0mm}}}
    \overline{q}^{u}    
    \left\{
    \begin{array}{@{\hspace{0mm}}c@{\hspace{0mm}}}
      -
      q
      \left\{
      \begin{array}{@{\hspace{0mm}}c@{\hspace{0mm}}}
        e^{3 6}_{4 5}
        \\
        e^{4 5}_{3 6}
      \end{array}
      \right\},
      \left\{
      \begin{array}{@{\hspace{0mm}}c@{\hspace{0mm}}}
        e^{2 7}_{4 5}
        \\
        e^{4 5}_{2 7}
      \end{array}
      \right\},
      - \overline{q}
      \left\{
      \begin{array}{@{\hspace{0mm}}c@{\hspace{0mm}}}
        e^{3 6}_{2 7}
        \\
        e^{2 7}_{3 6}
      \end{array}
      \right\}
    \end{array}
    \right\}
    \\
    q^{u}    
    \left\{
    \begin{array}{@{\hspace{0mm}}c@{\hspace{0mm}}}
      - \overline{q}
      \left\{
      \begin{array}{@{\hspace{0mm}}c@{\hspace{0mm}}}
        e^{6 3}_{5 4}
        \\
        e^{5 4}_{6 3}
      \end{array}
      \right\},
      \left\{
      \begin{array}{@{\hspace{0mm}}c@{\hspace{0mm}}}
        e^{7 2}_{5 4}
        \\
        e^{5 4}_{7 2}
      \end{array}
      \right\},
      - q
      \left\{
      \begin{array}{@{\hspace{0mm}}c@{\hspace{0mm}}}
        e^{6 3}_{7 2}
        \\
        e^{7 2}_{6 3}
      \end{array}
      \right\}
    \end{array}
    \right\}
  \end{array}
  \right\}
  \\
  & &
  \hspace{-45pt}
  \frac{
    A_{0}^{\frac{1}{2}}
    A_{2}^{\frac{1}{2}}
    U_{0}
    U_{1}
  }{
    S_{0}^{-}
    S_{1}^{-}
    S_{2}^{-}
  }
  \left\{
  \begin{array}{@{\hspace{0mm}}c@{\hspace{0mm}}}
    \overline{q}^{u}    
    \left\{
    \begin{array}{@{\hspace{0mm}}c@{\hspace{0mm}}}
      q
      \left\{
      \begin{array}{@{\hspace{0mm}}c@{\hspace{0mm}}}
        e^{1 8}_{7 2}
        \\
        e^{7 2}_{1 8}
      \end{array}
      \right\},
      -
      \left\{
      \begin{array}{@{\hspace{0mm}}c@{\hspace{0mm}}}
        e^{1 8}_{6 3}
        \\
        e^{6 3}_{1 8}
      \end{array}
      \right\},
      \overline{q}
      \left\{
      \begin{array}{@{\hspace{0mm}}c@{\hspace{0mm}}}
        e^{1 8}_{5 4}
        \\
        e^{5 4}_{1 8}
      \end{array}
      \right\}
    \end{array}
    \right\}
    \\
    q^{u}
    \left\{
    \begin{array}{@{\hspace{0mm}}c@{\hspace{0mm}}}
      \overline{q}
      \left\{
      \begin{array}{@{\hspace{0mm}}c@{\hspace{0mm}}}
        e^{8 1}_{2 7}
        \\
        e^{2 7}_{8 1}
      \end{array}
      \right\},
      -
      \left\{
      \begin{array}{@{\hspace{0mm}}c@{\hspace{0mm}}}
        e^{8 1}_{3 6}
        \\
        e^{3 6}_{8 1}
      \end{array}
      \right\},
      q
      \left\{
      \begin{array}{@{\hspace{0mm}}c@{\hspace{0mm}}}
        e^{8 1}_{4 5}
        \\
        e^{4 5}_{8 1}
      \end{array}
      \right\}
    \end{array}
    \right\}
  \end{array}
  \right\}
  \\
  & &
  \hspace{-45pt}
  \frac{
    A_{0}^{\frac{1}{2}}
    A_{1}
    A_{2}^{\frac{1}{2}}
    U_{0}
  }{
    S_{0}^{-}
    S_{1}^{-}
    S_{2}^{-}
  }
  \left\{
    \begin{array}{@{\hspace{0mm}}c@{\hspace{0mm}}}
    \overline{q}^{2 u}    
    \left\{
    \begin{array}{@{\hspace{0mm}}c@{\hspace{0mm}}}
      \overline{q}
      \left\{
      \begin{array}{@{\hspace{0mm}}c@{\hspace{0mm}}}
        \mathbf{+ e^{1 8}_{2 7}}
        \\
                - e^{2 7}_{1 8}
      \end{array}
      \right\},
      -
      \left\{
      \begin{array}{@{\hspace{0mm}}c@{\hspace{0mm}}}
        \mathbf{+ e^{1 8}_{3 6}}
        \\
                - e^{3 6}_{1 8}
      \end{array}
      \right\},
      q 
      \left\{
      \begin{array}{@{\hspace{0mm}}c@{\hspace{0mm}}}
        \mathbf{+ e^{1 8}_{4 5}}
        \\
                - e^{4 5}_{1 8}
      \end{array}
      \right\}
    \end{array}
    \right\}
    \\
    q^{2 u}    
    \left\{
      \begin{array}{@{\hspace{0mm}}c@{\hspace{0mm}}}
      q
      \left\{
      \begin{array}{@{\hspace{0mm}}c@{\hspace{0mm}}}
        \mathbf{+ e^{7 2}_{8 1}}
        \\
                - e^{8 1}_{7 2}
      \end{array}
      \right\},
      -
      \left\{
      \begin{array}{@{\hspace{0mm}}c@{\hspace{0mm}}}
        \mathbf{+ e^{6 3}_{8 1}}
        \\
                - e^{8 1}_{6 3}
      \end{array}
      \right\},
      \overline{q}
      \left\{
      \begin{array}{@{\hspace{0mm}}c@{\hspace{0mm}}}
        \mathbf{+ e^{5 4}_{8 1}}
        \\
                - e^{8 1}_{5 4}
      \end{array}
      \right\}
    \end{array}
    \right\}
    \end{array}
  \right\},
\end{eqnarray*}
\normalsize
where:
\small
\begin{eqnarray*}
  f_1 (q)
  & = &
  - 2 \overline{q}
  + (q^{1 + 2 \alpha} + \overline{q}^{1 + 2 \alpha})
  - \overline{q}^{2 u} (q - \overline{q})
  \\
  f_2 (q)
  & = &
  - 2 q
  + (\overline{q}^{3 + 2 \alpha} + q^{3 + 2 \alpha})
  + q^{2 u} (q - \overline{q})
  \\
  f_3 (q)
  & = &
  - \overline{q}^{u}
  (
    2 \overline{q}^{2}
    - (q^{2 + 2 \alpha} + \overline{q}^{2 + 2 \alpha})
    + \overline{q}^{2 u} (q^{2} - \overline{q}^{2})
  )
  \\
  f_4 (q)
  & = &
  \overline{q}^{u}
  (
    - 2 
    + (q^{2 + 2 \alpha} + \overline{q}^{2 + 2 \alpha})
    - (q^{2 - 2 u} + \overline{q}^{2 - 2 u})
    + (q^{2u} + \overline{q}^{2u})
  )
  \\
  f_5 (q)
  & = &
  \overline{q}^{u}
  (
    - 2 q^{2}
    + (q^{2 + 2 \alpha} + \overline{q}^{2 + 2 \alpha})
    + q^{2 u} (q^{2} - \overline{q}^{2})
  )
  \\
  f_6 (q)
  & = &
  q (q + \overline{q})
  - (q^{2 + 2 \alpha} + \overline{q}^{2 + 2 \alpha})
  - q^{2 u - 1} (q - \overline{q}).
\end{eqnarray*}

\normalsize

%% file: data/qrmx3.tex
\subsection{The quantum R matrix {$\check{R}^{3,1}$}}

$\check{R}^{3,1}$ has $139$ nonzero components:

\small

\begin{eqnarray*}
  & &
  \hspace{-45pt}
  1
  \left\{
    \begin{array}{@{\hspace{0mm}}c@{\hspace{0mm}}}
       e^{1 1}_{1 1}
    \end{array}
  \right\},
  \qquad
  -
  q^{2 \alpha}
  \left\{
    \begin{array}{@{\hspace{0mm}}c@{\hspace{0mm}}}
      e^{2 2}_{2 2},
      e^{3 3}_{3 3},
      e^{4 4}_{4 4},
    \end{array}
  \right\},
  \qquad
  q^{4 \alpha + 2}
  \left\{
    \begin{array}{@{\hspace{0mm}}c@{\hspace{0mm}}}
       e^{5 5}_{5 5},
       e^{6 6}_{6 6},
       e^{7 7}_{7 7}
    \end{array}
  \right\},
  \qquad
  -
  q^{6 \alpha + 6}
  \left\{
    \begin{array}{@{\hspace{0mm}}c@{\hspace{0mm}}}
      e^{8 8}_{8 8}
    \end{array}
  \right\}
  \\[1mm]
  \hline
  \\[1mm]
  & &
  \hspace{-45pt}
  -
  \Delta
  q^{\alpha}
  A_{0}
  \left\{
    \begin{array}{@{\hspace{0mm}}c@{\hspace{0mm}}}
      e^{2 1}_{2 1},
      e^{3 1}_{3 1},
      e^{4 1}_{4 1}
    \end{array}
  \right\},
  \qquad
  -
  \Delta
  q^{5 \alpha + 4}
  A_{2}
  \left\{
    \begin{array}{@{\hspace{0mm}}c@{\hspace{0mm}}}
      e^{8 7}_{8 7},
      e^{8 6}_{8 6},
      e^{8 5}_{8 5}
    \end{array}
  \right\}
  \\
  & &
  \hspace{-45pt}
  \Delta
  q^{2 \alpha + 1}
  \left\{
    \begin{array}{@{\hspace{0mm}}c@{\hspace{0mm}}}
       e^{3 2}_{3 2},
       e^{4 2}_{4 2},
       e^{4 3}_{4 3}
    \end{array}
  \right\},
  \qquad
  -
  \Delta
  q^{4 \alpha + 3}
  \left\{
    \begin{array}{@{\hspace{0mm}}c@{\hspace{0mm}}}
      e^{7 6}_{7 6},
      e^{7 5}_{7 5},
      e^{6 5}_{6 5}
    \end{array}
  \right\}
  \\
  & &
  \hspace{-45pt}
  \Delta^2
  q^{2 \alpha + 1}
  A_{0}
  A_{1}
  \left\{
    \begin{array}{@{\hspace{0mm}}c@{\hspace{0mm}}}
       e^{5 1}_{5 1},
       e^{6 1}_{6 1},
       e^{7 1}_{7 1}
    \end{array}
  \right\},
  \qquad
  -
  \Delta^2
  q^{4 \alpha + 3}
  A_{1}
  A_{2}
  \left\{
    \begin{array}{@{\hspace{0mm}}c@{\hspace{0mm}}}
      e^{8 4}_{8 4},
      e^{8 3}_{8 3},
      e^{8 2}_{8 2}
    \end{array}
  \right\}
  \\
  & &
  \hspace{-45pt}
  -
  \Delta^2
  q^{3 \alpha + 2}
  A_{1}
  \left\{
    \begin{array}{@{\hspace{0mm}}c@{\hspace{0mm}}}
       e^{6 3}_{6 3}
    \end{array}
  \right\},
  \qquad
  -
  \Delta
  q^{3 \alpha + 3}
  A_{1}
  (q^2 - \overline{q}^2)
  \left\{
  \begin{array}{@{\hspace{0mm}}c@{\hspace{0mm}}}
    e^{7 2}_{7 2}
  \end{array}
  \right\}
  \\
  & &
  \hspace{-45pt}
  \Delta
  q^{3 \alpha + 1}
  A_{1}
  \left\{
    \begin{array}{@{\hspace{0mm}}c@{\hspace{0mm}}}
       e^{5 2}_{5 2},
       e^{5 3}_{5 3},
       e^{6 2}_{6 2},
       e^{6 4}_{6 4},
       e^{7 3}_{7 3},
       e^{7 4}_{7 4}
    \end{array}
  \right\},
  \qquad
  -
  \Delta^3
  q^{3 \alpha + 3}
  A_{0}
  A_{1}
  A_{2}
  \left\{
    \begin{array}{@{\hspace{0mm}}c@{\hspace{0mm}}}
      e^{8 1}_{8 1}
    \end{array}
  \right\}
  \\[1mm]
  \hline
  \\[1mm]
  & &
  \hspace{-45pt}
  q^{\alpha}
  \left\{
    \begin{array}{@{\hspace{0mm}}c@{\hspace{0mm}}}
      \mathbf{- 1}
      \left\{
      \begin{array}{@{\hspace{0mm}}c@{\hspace{0mm}}}
       e^{1 2}_{2 1},
       e^{1 3}_{3 1},
       e^{1 4}_{4 1}
      \end{array}
      \right\}
      \\
      + 1
      \left\{
      \begin{array}{@{\hspace{0mm}}c@{\hspace{0mm}}}
        e^{2 1}_{1 2},
        e^{3 1}_{1 3},
        e^{4 1}_{1 4}
      \end{array}
      \right\}
    \end{array}
  \right\},
  \qquad
  q^{5 \alpha + 4}
  \left\{
    \begin{array}{@{\hspace{0mm}}c@{\hspace{0mm}}}
      \mathbf{- 1}
      \left\{
      \begin{array}{@{\hspace{0mm}}c@{\hspace{0mm}}}
        e^{7 8}_{8 7},
        e^{6 8}_{8 6},
        e^{5 8}_{8 5}
      \end{array}
      \right\}
      \\
      + 1
      \left\{
      \begin{array}{@{\hspace{0mm}}c@{\hspace{0mm}}}
        e^{8 7}_{7 8},
        e^{8 6}_{6 8},
        e^{8 5}_{5 8}
      \end{array}
      \right\}
    \end{array}
  \right\}
  \\
  & &
  \hspace{-45pt}
  q^{2 \alpha}
  \left\{
    \begin{array}{@{\hspace{0mm}}c@{\hspace{0mm}}}
      e^{1 5}_{5 1},
      e^{1 6}_{6 1},
      e^{1 7}_{7 1}
      \\
      e^{5 1}_{1 5},
      e^{6 1}_{1 6},
      e^{7 1}_{1 7}
    \end{array}
  \right\},
  \qquad
  -
  q^{4 \alpha + 2}
  \left\{
    \begin{array}{@{\hspace{0mm}}c@{\hspace{0mm}}}
      e^{8 4}_{4 8},
      e^{8 3}_{3 8},
      e^{8 2}_{2 8}
      \\
      e^{4 8}_{8 4},
      e^{3 8}_{8 3},
      e^{2 8}_{8 2}
    \end{array}
  \right\}
  \\
  & &
  \hspace{-45pt}
  -
  q^{2 \alpha + 1}
  \left\{
    \begin{array}{@{\hspace{0mm}}c@{\hspace{0mm}}}
      e^{2 3}_{3 2},
      e^{2 4}_{4 2},
      e^{3 4}_{4 3}
      \\
      e^{3 2}_{2 3},
      e^{4 2}_{2 4},
      e^{4 3}_{3 4}
    \end{array}
  \right\},
  \qquad
  q^{4 \alpha + 3}
  \left\{
    \begin{array}{@{\hspace{0mm}}c@{\hspace{0mm}}}
       e^{7 6}_{6 7},
       e^{7 5}_{5 7},
       e^{6 5}_{5 6}
       \\
       e^{6 7}_{7 6},
       e^{5 7}_{7 5},
       e^{5 6}_{6 5}
    \end{array}
  \right\}
  \\
  & &
  \hspace{-45pt}
  q^{3 \alpha + 1}
  \left\{
    \begin{array}{@{\hspace{0mm}}c@{\hspace{0mm}}}
      \mathbf{- 1}
      \left\{
      \begin{array}{@{\hspace{0mm}}c@{\hspace{0mm}}}
        e^{5 2}_{2 5},
        e^{6 2}_{2 6},
        e^{5 3}_{3 5},
        e^{7 3}_{3 7},
        e^{6 4}_{4 6},
        e^{7 4}_{4 7}
      \end{array}
      \right\}
      \\
      + 1
      \left\{
      \begin{array}{@{\hspace{0mm}}c@{\hspace{0mm}}}
        e^{2 5}_{5 2},
        e^{2 6}_{6 2},
        e^{3 5}_{5 3},
        e^{3 7}_{7 3},
        e^{4 6}_{6 4},
        e^{4 7}_{7 4}
      \end{array}
      \right\}
    \end{array}
  \right\},
  q^{3 \alpha + 2}
  \left\{
    \begin{array}{@{\hspace{0mm}}c@{\hspace{0mm}}}
      \mathbf{- 1}
      \left\{
      \begin{array}{@{\hspace{0mm}}c@{\hspace{0mm}}}
        e^{7 2}_{2 7},
        e^{6 3}_{3 6},
        e^{5 4}_{4 5}
      \end{array}
      \right\}
      \\
      + 1
      \left\{
      \begin{array}{@{\hspace{0mm}}c@{\hspace{0mm}}}
        e^{2 7}_{7 2},
        e^{3 6}_{6 3},
        e^{4 5}_{5 4}
      \end{array}
      \right\}
    \end{array}
  \right\},
  q^{3 \alpha}
  \left\{
    \begin{array}{@{\hspace{0mm}}c@{\hspace{0mm}}}
      \mathbf{- 1}
      \left\{
      \begin{array}{@{\hspace{0mm}}c@{\hspace{0mm}}}
        e^{1 8}_{8 1}
      \end{array}
      \right\}
      \\
      + 1
      \left\{
      \begin{array}{@{\hspace{0mm}}c@{\hspace{0mm}}}
        e^{8 1}_{1 8}
      \end{array}
      \right\}
    \end{array}
  \right\}
  \\[1mm]
  \hline
  \\[1mm]
  & &
  \hspace{-45pt}
  \Delta
  q^{2 \alpha + 1}
  A_{0}
  A_{1}
  \left\{
    \overline{q}^{\frac{1}{2}}
    \left\{
    \begin{array}{@{\hspace{0mm}}c@{\hspace{0mm}}}
      \mathbf{+ 1}
      \left\{
      \begin{array}{@{\hspace{0mm}}c@{\hspace{0mm}}}
        e^{5 1}_{2 3},
        e^{6 1}_{2 4},
        e^{7 1}_{3 4}
      \end{array}
      \right\}
      \\
      - 1
      \left\{
      \begin{array}{@{\hspace{0mm}}c@{\hspace{0mm}}}
         e^{2 3}_{5 1},
         e^{2 4}_{6 1},
         e^{3 4}_{7 1}
      \end{array}
      \right\}
    \end{array}
    \right\},
    q^{\frac{1}{2}}
    \left\{
    \begin{array}{@{\hspace{0mm}}c@{\hspace{0mm}}}
      \mathbf{- 1}
      \left\{
      \begin{array}{@{\hspace{0mm}}c@{\hspace{0mm}}}
        e^{5 1}_{3 2},
        e^{6 1}_{4 2},
        e^{7 1}_{4 3}
      \end{array}
      \right\}
      \\
      + 1
      \left\{
      \begin{array}{@{\hspace{0mm}}c@{\hspace{0mm}}}
        e^{3 2}_{5 1},
        e^{4 2}_{6 1},
        e^{4 3}_{7 1}
      \end{array}
      \right\}
    \end{array}
    \right\}
  \right\}
  \\
  & &
  \hspace{-45pt}
  \Delta
  q^{4 \alpha + 3}
  A_{1}
  A_{2}
  \left\{
    \overline{q}^{\frac{1}{2}}
    \left\{
    \begin{array}{@{\hspace{0mm}}c@{\hspace{0mm}}}
      \mathbf{+ 1}
      \left\{
      \begin{array}{@{\hspace{0mm}}c@{\hspace{0mm}}}
        e^{5 6}_{8 2},
        e^{5 7}_{8 3},
        e^{6 7}_{8 4}
      \end{array}
      \right\}
      \\
      - 1
      \left\{
      \begin{array}{@{\hspace{0mm}}c@{\hspace{0mm}}}
        e^{8 2}_{5 6},
        e^{8 3}_{5 7},
        e^{8 4}_{6 7}
      \end{array}
      \right\}
    \end{array}
    \right\},
    q^{\frac{1}{2}}
    \left\{
    \begin{array}{@{\hspace{0mm}}c@{\hspace{0mm}}}
      \mathbf{- 1}
      \left\{
      \begin{array}{@{\hspace{0mm}}c@{\hspace{0mm}}}
        e^{6 5}_{8 2},
        e^{7 5}_{8 3},
        e^{7 6}_{8 4}
      \end{array}
      \right\}
      \\
      + 1
      \left\{
      \begin{array}{@{\hspace{0mm}}c@{\hspace{0mm}}}
        e^{8 2}_{6 5},
        e^{8 3}_{7 5},
        e^{8 4}_{7 6}
      \end{array}
      \right\}
    \end{array}
    \right\}
  \right\}
  \\
  & &
  \hspace{-45pt}
  \Delta
  q^{3 \alpha + \frac{5}{2}}
  \left\{
    \overline{q}^{\frac{1}{2}}
    \left\{
    \begin{array}{@{\hspace{0mm}}c@{\hspace{0mm}}}
      \mathbf{+ 1}
      \left\{
      \begin{array}{@{\hspace{0mm}}c@{\hspace{0mm}}}
        e^{6 3}_{4 5},
        e^{7 2}_{3 6}
      \end{array}
      \right\}
      \\
      - 1
      \left\{
      \begin{array}{@{\hspace{0mm}}c@{\hspace{0mm}}}
        e^{4 5}_{6 3},
        e^{3 6}_{7 2}
      \end{array}
      \right\}
    \end{array}
    \right\},
    q^{\frac{1}{2}}
    \left\{
    \begin{array}{@{\hspace{0mm}}c@{\hspace{0mm}}}
      \mathbf{- 1}
      \left\{
      \begin{array}{@{\hspace{0mm}}c@{\hspace{0mm}}}
        e^{7 2}_{4 5}
      \end{array}
      \right\}
      \\
      + 1
      \left\{
      \begin{array}{@{\hspace{0mm}}c@{\hspace{0mm}}}
        e^{4 5}_{7 2}
      \end{array}
      \right\}
    \end{array}
    \right\}
  \right\}
  \\
  & &
  \hspace{-45pt}
  \Delta
  q^{3 \alpha + 3}
  A_{1}
  \left\{
    \overline{q}
    \left\{
    \begin{array}{@{\hspace{0mm}}c@{\hspace{0mm}}}
      e^{6 3}_{5 4}
      \\
      e^{5 4}_{6 3},
    \end{array}
    \right\},
    -
    \left\{
    \begin{array}{@{\hspace{0mm}}c@{\hspace{0mm}}}
      e^{7 2}_{5 4}
      \\
      e^{5 4}_{7 2}
    \end{array}
    \right\},
    q
    \left\{
    \begin{array}{@{\hspace{0mm}}c@{\hspace{0mm}}}
      e^{6 3}_{7 2}
      \\
      e^{7 2}_{6 3}
    \end{array}
    \right\}
  \right\}
  \\
  & &
  \hspace{-45pt}
  \Delta
  q^{3 \alpha + 2}
  A_{0}
  A_{2}
  \left\{
    -
    \overline{q}
    \left\{
    \begin{array}{@{\hspace{0mm}}c@{\hspace{0mm}}}
      e^{8 1}_{2 7}
      \\
      e^{2 7}_{8 1}
    \end{array}
    \right\},
    \left\{
    \begin{array}{@{\hspace{0mm}}c@{\hspace{0mm}}}
      e^{8 1}_{3 6}
      \\
      e^{3 6}_{8 1}
    \end{array}
    \right\}
    -
    q
    \left\{
    \begin{array}{@{\hspace{0mm}}c@{\hspace{0mm}}}
      e^{8 1}_{4 5}
      \\
      e^{4 5}_{8 1}
    \end{array}
    \right\}
  \right\}
  \\
  & &
  \hspace{-45pt}
  \Delta^2
  q^{3 \alpha + 3}
  A_{0}
  A_{1}
  A_{2}
  \left\{
    \overline{q}
    \left\{
    \begin{array}{@{\hspace{0mm}}c@{\hspace{0mm}}}
      \mathbf{- e^{5 4}_{8 1}}
      \\
              + e^{8 1}_{5 4}
    \end{array}
    \right\},
    -
    \left\{
    \begin{array}{@{\hspace{0mm}}c@{\hspace{0mm}}}
      \mathbf{- e^{6 3}_{8 1}}
      \\
              + e^{8 1}_{6 3}
    \end{array}
    \right\},
    q
    \left\{
    \begin{array}{@{\hspace{0mm}}c@{\hspace{0mm}}}
      \mathbf{- e^{7 2}_{8 1}}
      \\
              + e^{8 1}_{7 2}
    \end{array}
    \right\}
  \right\}.
\end{eqnarray*}

\normalsize